\documentclass[a4paper,
11pt, reqno]{amsart}
\usepackage[utf8x]{inputenc}
\usepackage{comment}
\usepackage{amssymb,amsmath,latexsym}
\usepackage{epsfig}
\usepackage{mathrsfs}
\usepackage{hyperref}
\usepackage{color}
\usepackage{epic,eepic,wrapfig,color, ifthen}
\usepackage{url}
\usepackage{hyperref}

\newtheorem{thm}{Theorem}[section]
\newtheorem{cor}[thm]{Corollary}
\newtheorem{lemma}[thm]{Lemma}
\newtheorem{prop}[thm]{Proposition}
\newtheorem{defn}[thm]{Definition}

\newtheorem{remark}[thm]{Remark}
\newtheorem{example}[thm]{Example}
\numberwithin{equation}{section}

\newcommand{\formula}[2][nolabel]
{\ifthenelse{\equal{#1}{nolabel}}
 {\begin{align*} #2 \end{align*}}
 {\ifthenelse{\equal{#1}{}}
  {\begin{align} #2 \end{align}}
  {\begin{align} \label{#1} #2 \end{align}}
 }
}

\def\pf{{\medskip\noindent {\bf Proof. }}}
\def\qed{{\hfill $\Box$ \bigskip}}
\def\R{{\mathbb R}}
\def\P{{\mathbb P}}
\def\E{{\mathbb E}}
\def\1{{\bf 1}}

\def\diam{{\mathrm{diam}}}

\DeclareMathOperator{\dist}{dist}
\DeclareMathOperator*{\esssup}{ess\,sup}

\newcommand{\cal}[1]{\mathcal{#1}}
 \def\sB {{\cal B}} 
  \def\FF{{\mathcal F}}

\def\R {{\mathbb R}}
\def\N {{\mathbb N}}

\def\HKBh {${\bf HK^h_B}$}

\def\HKUh {${\bf HK^h_U}$}

\def\nn{\nonumber}

\def\wt{\widetilde}
\def\wh{\widehat}

\def\E{{\mathbb E}}
\def\P{{\mathbb P}}
\def\eps{{\epsilon}}

\def\bea{\begin{align*}}
\def\eea{\end{align*}}
\def\bee{\begin{equation}}
\def\eee{\end{equation}}

\makeatletter
\@addtoreset{equation}{section}

\makeatother

\begin{document}
\bibliographystyle{plain}
\allowdisplaybreaks

\title[Heat kernel estimates for subordinate Markov processes]
{ \bf Heat kernel estimates for subordinate Markov processes and their applications}

\author{ Soobin Cho, \quad Panki Kim, \quad Renming Song \quad and \quad Zoran Vondra\v{c}ek}

\address[Cho]{Department of Mathematical Sciences,
	Seoul National University,
	Building 27, 1 Gwanak-ro, Gwanak-gu
	Seoul 08826, Republic of Korea}\thanks{The research of Soobin Cho is supported by the POSCO Science Fellowship of POSCO TJ Park Foundation.}
\curraddr{}
\email{soobin15@snu.ac.kr}

\address[Kim]{Department of Mathematical Sciences and Research Institute of Mathematics,
	Seoul National University,
	Building 27, 1 Gwanak-ro, Gwanak-gu
	Seoul 08826, Republic of Korea}\thanks{This research is  supported by the National Research Foundation of Korea(NRF) grant funded by the Korea government(MSIP) (No. NRF-2021R1A4A1027378).
}
\curraddr{}
\email{pkim@snu.ac.kr}

\address[Song]{
Department of Mathematics, University of Illinois, Urbana, IL 61801,
USA}
\curraddr{}\thanks{Research supported in part by a grant from
the Simons Foundation (\#429343, Renming Song)}
\email{rsong@math.uiuc.edu}

\address[Vondra\v{c}ek]
{
Department of Mathematics, Faculty of Science, University of Zagreb, Zagreb, Croatia,
}
\curraddr{}\thanks{ Research supported in part by the Croatian Science Foundation under the project 4197. 
(Zoran Vondra\v{c}ek)}
\email{vondra@math.hr}

 \date{}

\begin{abstract}
In this paper, we establish sharp two-sided estimates for transition densities of a large class of subordinate Markov processes. As applications, we show that 
the parabolic Harnack inequality and 
H\"older regularity hold for  parabolic functions of  such processes,
and derive sharp two-sided Green function estimates.
\end{abstract}

\maketitle



%%%%%%%%%%%%%%%%%%%%%%%%%%%%%%%%%%%%%%%%%%%%%%%%%%%%%%%%%%%%%%%%%%%%%%%%%%%%%%%%%%%%%%
%%%%%%%%%%%%%%%%                                                          Introduction                                                %%%%%%%%%%%%%%%%%%%%%%%%%%%%%%%
%%%%%%%%%%%%%%%%%%%%%%%%%%%%%%%%%%%%%%%%%%%%%%%%%%%%%%%%%%%%%%%%%%%%%%%%%%%%%%%%%%%%%% 

\allowdisplaybreaks
\section{Introduction}
Transition densities of Markov processes are of central importance in both probability
and analysis. The transition density 
$p(t, x, y)$ of a Markov process $X$ with generator $L$
is the fundamental solution of the equation $\partial_tu=Lu$. Hence the transition density $p(t, x, y)$ is also known as the heat kernel of $L$. The heat kernel is rarely known explicitly. Due to the importance
of heat kernels, there is a huge amount of literature devoted to estimates of heat kernels. 
 
The purpose of this paper is to study 
heat kernel estimates for subordinate Markov processes. 
The main motivation comes from \cite{KSV20}, where it was established that the jump kernels of subordinate killed L\'evy processes have 
an unusual form not observed before.
It is therefore plausible that the heat kernels of those processes will have some new features.
It turns out that this is indeed the case. To illustrate the new features, we explain below the motivating and also the simplest example covered by our results.

Let $D\subset \R^d$, 
$d \ge 1$,  be a 
bounded $C^{1,1}$ open set. 
For $x\in D$, let $\delta_D(x)$ denote the distance between $x$ and 
$D^c$. Let $Y$ be an isotropic $\alpha$-stable process in $\R^d$,
 $\alpha \in (0,2]$  and let $Y^D$ denote  the part process  of $Y$ 
killed upon exiting $D$. 
When $\alpha =2$, we further assume that $D$ is connected.
Sharp two-sided estimates of the heat kernel 
$p_D(t,x,y)$ of $Y^D$ were obtained in  \cite{Da, Zh}
(for $\alpha=2$)   and  \cite{CKS-jems} (for $\alpha<2$): 
there exist positive constants $c_i$, $i=1, \dots, 8,$ such that following estimates hold true.
 For $(t,x,y)\in (0,1]\times D\times D$,
$$
c_1h(t,x,y)\Big(t^{-d/\alpha}\wedge \frac{t}{|x-y|^{d+\alpha}}\Big) \le p_D(t,x,y)\le c_2
	h(t,x,y)\Big(t^{-d/\alpha}\wedge \frac{t}{|x-y|^{d+\alpha}}\Big), \quad\text{for }  \alpha<2,
$$
and 
$$
c_3h(t,x,y) \, t^{-d/2}e^{-c_4|x-y|^2/t} \le p_D(t,x,y)\le c_5
	h(t,x,y) \, t^{-d/2}e^{-c_6|x-y|^2/t} , \quad \text{for }  \alpha=2,
$$
where the boundary function $h(t,x,y)$ is 
given by
$$
h(t,x,y)=\Big(1\wedge \frac{\delta_D(x)^{\alpha} }{t}\Big)^{1/2}\Big(1\wedge \frac{\delta_D(y)^{\alpha} }{t}\Big)^{1/2}  =\Big(1\wedge \frac{\delta_D(x) }{t^{1/\alpha}}\Big)^{\alpha/2}\Big(1\wedge \frac{\delta_D(y) }{t^{1/\alpha}}\Big)^{\alpha/2}.$$ 
For $(t,x,y)\in [1,\infty)\times D \times D$, 
$$
c_7e^{-\lambda_1 t}\delta_D(x)^{\alpha/2}\delta_D(y)^{\alpha/2} \le p_D(t,x,y)\le c_8 e^{-\lambda_1 t}\delta_D(x)^{\alpha/2}\delta_D(y)^{\alpha/2},
$$
where $\lambda_1$ is the smallest eigenvalue of   the 
Dirichlet (fractional) Laplacian 
$(-\Delta)^{\alpha/2}\big |_D$. 

 Let $S=(S_t)_{t\ge 0}$ be a $\beta$-stable subordinator, $\beta\in (0,1)$,
independent of $Y^D$, and let $X=(X_t)_{t\ge 0}$ be the subordinate process: $X_t:=Y^D_{S_t}$. The generator of $X$ is equal to  
(the negative of)  
$\big((-\Delta)^{\alpha/2}\big |_D\big)^{\beta}$ -- the fractional power of the Dirichlet fractional Laplacian. 
In particular, when $\alpha=2$, $\big(-\Delta\big |_D\big)^{\beta}$ is 
called a spectral fractional Laplacian in the PDE literature, see, for instance,
\cite{BSV} and the references therein.
In this respect, the process $X$ bears some similarity with the isotropic $\alpha \beta $-stable process.
The heat kernel $q(t,x,y)$ of the subordinate process $X$ is given by
$$
q(t,x,y)=\int_0^{\infty} p_D(s,x,y)\P(S_t\in ds), \quad t>0, \ x,y\in D.
$$
Note that the distribution of $S_t$ is not explicitly known, making the above 
 integration rather delicate. 
To handle this integral, we establish some estimates on the distribution of $S_t$
 (in fact, for much more general subordinators than the stable ones).
Using these, we can
obtain sharp two-sided estimates of $q(t,x,y)$. 
To  present those estimates, we first introduce 
some notation. 
Denote $a\wedge b:=\min\{a,b\}$ and $a\vee b:=\max\{a,b\}$. The notation $f(s) \simeq  g(s)$ means that there exist comparison constants $c_1,c_2>0$ such that $c_1g(s)\leq f (s)\leq c_2 g(s)$ for 
specified range of the variable $s$.  
For $x,y\in D$, let
\begin{eqnarray*}
	\delta_{\vee}(x,y)=\delta_D(x)\vee \delta_D(y), &\quad & \delta_{\wedge}(x,y)=\delta_D(x)\wedge \delta_D(y).\end{eqnarray*}

Our main results, specialized to the present situation, are summarized below.

\begin{thm}\label{t:special}
 (1) For all $(t,x,y)\in (0,1]\times D\times D$, 
\begin{align}
q(t,x,y) &\simeq \Big(1 \wedge \frac{\delta_D(x)}{t^{1/(\alpha\beta)}} \Big)^{\alpha/2} \Big(1 \wedge \frac{\delta_D(y)}{t^{1/(\alpha\beta)}} \Big)^{\alpha/2} B^{\alpha,\beta}(t,x,y)\left( t^{-d/(\alpha \beta )} \wedge \frac{t}{|x-y|^{d+\alpha\beta}} \right)\label{e:nHKE1}\\
&\simeq \Big(1 \wedge \frac{\delta_D(x)}{t^{1/(\alpha\beta)}} \Big)^{\alpha/2} \Big(1 \wedge \frac{\delta_D(y)}{t^{1/(\alpha\beta)}} \Big)^{\alpha/2}\left( t^{-d/(\alpha \beta )} \wedge \frac{tB^{\alpha,\beta}(t,x,y)}{|x-y|^{d+\alpha\beta}} \right),\label{e:nHKE2}
\end{align}
where
$$
B^{2,\beta}(t,x,y):= \Big(1 \wedge \frac{\delta_D(x) \vee t^{1/(2\beta)}}{|x-y|}\Big)  \Big(1 \wedge \frac{\delta_D(y) \vee t^{1/(2\beta)}}{|x-y|}\Big)
$$
and for $\alpha<2$,
\begin{equation*}
	B^{\alpha,\beta}(t,x,y):= \begin{cases} \displaystyle
		\bigg(1 \wedge \frac{\delta_\wedge(x,y) \vee t^{1/(\alpha\beta)}}{|x-y|} \bigg)^{\alpha-\alpha\beta}, 	&\mbox{if }\, \beta>\displaystyle\frac12,\\[10pt]
	\displaystyle
	\bigg(1 \wedge \frac{\delta_\wedge(x,y) \vee t^{1/(\alpha\beta)}}{|x-y|} \bigg)^{\alpha/2} 	\bigg(1 \wedge \frac{\delta_{ \vee }(x,y) \vee t^{1/(\alpha\beta)}}{|x-y|} \bigg)^{\alpha/2-\alpha\beta}, &\mbox{if }\, \beta <\displaystyle\frac12, \\[10pt]
	\displaystyle	\bigg(1 \wedge \frac{\delta_\wedge(x,y)\vee t^{1/(\alpha\beta)}}{|x-y|} \bigg)^{\alpha/2} \log \bigg( e+ \frac{(\delta_\vee(x,y) \vee t^{1/(\alpha\beta)}) \wedge |x-y|}{ ( \delta_\wedge(x,y) \vee t^{1/(\alpha\beta)} ) \wedge |x-y|  }\bigg), &\mbox{if }\, \beta=\displaystyle\frac12.
	\end{cases} 
\end{equation*}

(2) For all $(t,x,y)\in [1,\infty)\times D \times D$,
$$
q(t,x,y)\simeq  e^{-t \lambda_1^\beta}
\delta_D(x)^{\alpha/2} \delta_D(y)^{\alpha/2},
$$
where $\lambda_1$ is the smallest eigenvalue of   the 
Dirichlet (fractional) Laplacian 
$(-\Delta)^{\alpha/2}\big |_D$.
\end{thm}

\begin{remark}\label{r:new}
{\rm
(1) 
 From  the forms of the heat kernel estimates \eqref{e:nHKE1}-\eqref{e:nHKE2} in Theorem \ref{t:special}(1), 
 one  can easily see the following: 
 (1-a) For  
$x, y$ away from the boundary (in the sense that 
$\delta_\wedge(x, y)
\ge |x-y|\vee t^{1/(\alpha\beta)}$), 
and for all $\beta\in (0,1)$, it holds that
\begin{align}\label{e:nqe}
q(t,x,y)\simeq t^{-d/(\alpha \beta )} \wedge \frac{t}{|x-y|^{d+\alpha \beta }}.
\end{align}
(1-b)
Dividing  
\eqref{e:nHKE1} by $t$ and letting $t\to 0$,  we can deduce that 
 the jump kernel is comparable with 
$$
\frac{B^{\alpha,\beta}(0,x,y)}{|x-y|^{d+\alpha\beta}}.
$$
 This indicates that the term $B^{\alpha, \beta}(t,x,y)$  in 
\eqref{e:nHKE1} comes from the boundary decay of the jump kernel.

\noindent	
(2) 
Using the identity
$$
	\bigg(1 \wedge \frac{r}{a} \bigg)\bigg(1 \wedge \frac{r \vee a}{c} \bigg)= \bigg(1 \wedge \frac{r}{a \vee c} \bigg), \quad a,c,r>0,$$
we see from Theorem \ref{t:special}(1) that, when $\alpha=2$, 
$$
	q(t,x,y)\simeq
	\Big(1 \wedge \frac{\delta_D(x)}{|x-y| \vee t^{1/\beta}} \Big)
		\Big(1 \wedge \frac{\delta_D(y)}{|x-y| \vee t^{1/\beta}} \Big)
	\left( t^{-d/(2 \beta )} \wedge \frac{t}{|x-y|^{d+2\beta}} \right).
$$	
}
\end{remark}

Recall that the two-sided estimates of the form \eqref{e:nqe} are valid for the heat kernel of 
the isotropic $\alpha \beta $-stable process in the whole space.
The novelty of the estimates for $q(t,x,y)$ is in the boundary term, which is quite unusual and 
involves interplays among $\delta_{\vee}(x,y)$, $\delta_{\wedge}(x,y)$ and time $t$ itself. 
In this respect, the form of the boundary term is very different 
from the boundary function $h(t,x,y)$ for the underlying process $Y^D$.

\iffalse
The Green function $G_D(x,y)$ of the process $X$ in $D$ is given by
\begin{equation*}
	G_D(x,y) = \int_0^\infty q(t,x,y)dt, \quad x,y\in D.
\end{equation*}
Then for all $x,y \in D$, 

\begin{equation}\label{e:Example_Green}
G_D(x,y)\simeq\left\{ \begin{array}{ll} 
\Big(1 \wedge \frac{\delta_D(x)^{\alpha/2} \delta_D(y)^{\alpha/2}}{|x-y|^\alpha} \right) \frac{1}{|x-y|^{d-\alpha \beta }},& d>\alpha \beta ,\\[10pt]
(\delta_D(x)\delta_D(y))^{(\alpha \beta -d)/2} \wedge \frac{\delta_D(x)^{\alpha/2}\delta_D(y)^{\alpha/2}}{|x-y|^{d+\alpha-\alpha \beta }} , & d<\alpha \beta , \\[10pt]
 \log \Big(1+ \frac{ \delta_D(x)^{\alpha/2}\delta_D(y)^{\alpha/2}}{|x-y|^\alpha}\right), & d=\alpha \beta .
\end{array}\right. 
\end{equation}
In case of $d \ge 2>\alpha \beta $, the above Green function estimates are obtained in \cite[Theorem 6.4]{KSV20}.  Letting $\beta \to  1$, \eqref{e:Example_Green} recovers the Green function estimates for killed isotropic $\alpha$-stable process in $D$; see  \cite[p. 182]{CZ} for $d=\alpha=2$, \cite{Zhao} for $d\ge 3$ and $\alpha=2$, and  \cite[Corollary 1.2]{CKS-jems} for $\alpha<2$.
\fi

All the estimates in Theorem \ref{t:special}
are consequences of  Theorems \ref{t:smallgeneral}, 
\ref{t:Slarge}
and Lemma \ref{l:asym-Bp}, see Example \ref{ex:basic}. 
Integrating the heat kernel estimates, we can obtain sharp two-sided estimates on the Green function of $X$, see Theorem \ref{t:Sgreen}  and Example \ref{ex:basic}.

In this paper, we obtain sharp two-sided  heat kernel estimates for subordinate Markov processes in a setting 
which is more general, in several directions, than 
that of the example above.
We allow (i) quite general subordinators, (ii) Markov processes with state space $D$ that is either a bounded or an unbounded subset of a locally compact separable metric space, and (iii) very general form of two-sided estimates of the heat kernel $p_D(t,x,y)$ of the underlying process. 
In the remaining part of the introduction,
we describe some of our assumptions and results,
 and lay out the structure of the paper.

In Section \ref{s:preliminaries}, we first introduce the main assumption on the subordinator $S=(S_t)_{t\ge 0}$. Let $\nu$ denote its L\'evy measure and 
$w(t):=\nu(t, \infty)$.
We assume that there exist constants $R_1\in (0,\infty]$, $c_1, c_2>0$ and $\beta_2\ge \beta_1>0$ such that
$$
c_1 \Big(\frac{R}{r}\Big)^{\beta_1} \le \frac{w(r)}{w(R)} \le c_2 \Big(\frac{R}{r}\Big)^{\beta_2}    \quad \text{for all} \;\; 0<r \le R  <R_1.
$$
This assumption is quite weak -- it implies that 
the Laplace exponent $\phi$ of the subordinator $S$ satisfies weak scaling conditions near infinity with lower index $\beta_1$ and upper index $\beta_2 \wedge 1$. 
Note that $\beta_2>1$ is allowed.
Building upon the results from \cite{CK1, CK20}, we show several auxiliary results leading to the important estimate \eqref{e:tw} saying that $\P(S_t \ge s) \simeq tw(s)$ 
for $2\phi^{-1}(1/t)^{-1}<s<R_1/2$.

 Sections \ref{s:setup}--\ref{s:key-estimates} are central to the paper. 
We start with the setup in Section \ref{s:setup}: 
the underlying space is a locally compact separable metric space $E$ with a Radon measure $m$ having full support and 
satisfying volume doubling conditions. 
The state space is a proper open subset $D$ of $E$, bounded or unbounded, and $Y^D$ is a Hunt process living on $D$. We assume that $Y^D$ admits a transition density $p_D(t,x,y)$. The main assumption on the transition density is given in Definition \ref{d:pD-estimates}. Depending on whether $D$ is bounded or unbounded, the assumptions are somewhat different. Roughly speaking, at least for small times, $p_D(t,x,y)$ is comparable to the product of two parts which we may call the boundary part and the interior part. The latter is specified in terms of the volume and two functions -- $\Psi\ge \Phi$ -- both enjoying the scaling property, and also includes a Gaussian part. We note that although in most examples it holds that $\Psi\simeq \Phi$, allowing for different functions enlarges the scope of examples. The boundary part is described through a boundary function $h(t,x,y)$ which is not specified but is required to satisfy certain assumptions -- see Definition \ref{d:df}. 
To justify  our assumptions on the heat kernel $p_D(t,x,y)$, we provide a number of examples from the literature satisfying those assumptions. 

The main object of our study is the subordinate process $X_t:=Y^D_{S_t}$. In Subsection \ref{ss:jump}, see Theorem \ref{t:jump-estimate}, we first establish sharp two-sided estimates of the jump kernel $J(x,y)$ of
 $X$, thus generalizing  
\cite[Theorem 8.4]{KSV20}.
Subsection \ref{ss:key} contains sharp two-sided estimates of the heat kernel $q(t,x,y)$ of 
$X$ which are the main results of the paper. The case of a bounded set $D$ and small time is given in Theorem \ref{t:smallgeneral}, and the case of unbounded $D$ and all time in Theorem \ref{t:largegeneral}. The near-diagonal estimates have a rather simple form, but the off-diagonal estimates are quite involved, containing four terms which cannot be compared under general assumptions. The form of the estimates is somewhat simplified in case  when 
the upper scaling index $\beta_2<1$ and $\Psi\simeq \Phi$, 
see Corollary \ref{c:interior-small}. Finally, Theorem \ref{t:Slarge} provides 
the large time estimates in case of bounded $D$.

In Section \ref{s:green} we apply our sharp two-sided heat kernel estimates to derive sharp two-sided estimates of the Green function $G(x,y)$ of 
$X$. The most general form of the estimates is given in Proposition \ref{p:Sgreen}. These can be simplified under additional assumptions on the boundary function $h(t,x,y)$. We obtain several forms of the estimates depending on the relationship between 
the parameters in the volume doubling condition, the scaling indices of the functions $w$ and $\Phi$, and the parameters coming from $h$. 
The main results of this section are Theorems \ref{t:Sgreen}, \ref{t:Sgreen_2} and \ref{t:Sgreen_3}.

In Section \ref{s:phi} we show that parabolic functions with respect to 
$X$ satisfy  H\"older regularity and 
the parabolic Harnack inequality.
By using the rough upper estimates and the interior estimates for the heat kernel from Proposition \ref{p:hke-upper-rough} and Corollary \ref{c:smallint} together with the jump kernel estimates from Theorem \ref{t:jump-estimate}, we establish that the process $X$ satisfies all the assumptions from \cite{CKK09, CKW20}  used in the proofs of those results.

Finally, in Section \ref{s:examples}, we 
provide a concrete and explicit example
which includes the motivating example from the beginning of this introduction, 
and, with help of Lemma \ref{l:asym-Bp}, derive the heat kernel estimates 
using the general  Theorems \ref{t:smallgeneral}, \ref{t:largegeneral}
and \ref{t:Slarge}. We also derive the jump kernel estimates and Green function 
estimates using
Theorem \ref{t:jump-estimate}, and Theorems \ref{t:Sgreen}, \ref{t:Sgreen_2}, \ref{t:Sgreen_3} respectively. 
As an application of these estimates, combined with the main results of \cite{KSV21},
we completely  determine the  region of the parameters where the boundary Harnack principle holds for the process $X_t=Y^D_{S_t}$, where $D$ is the upper half-space,
$Y^D$ is the process in Example \ref{ex:HK} (b-4) 
and $S_t$ is an independent
$\beta$-stable subordinator, $\beta\in (0, 1)$.
At the end we provide an interesting example in which the upper scaling index $\beta_2 >1$ and the two scaling function $\Phi$ and $\Psi$ are different.

{\bf Notations}:  Values of lower case letters with subscripts $c_i$, $i=0,1,2,...$ are fixed in each statement and proof, and the labeling of these constants starts anew in each proof.
Recall that  $a\wedge b:=\min\{a,b\}$, $a\vee b:=\max\{a,b\}$.
We use two notations 
for comparison of functions.
First, the notation $f(x) \simeq g(x)$ means that there exist  constants $c_1,c_2>0$ such that $c_1g(x)\leq f (x)\leq c_2 g(x)$ for 
specified range of $x$.
On the other hand, the notation $f(x) \asymp g_1(x) + g_2(x)h(cx)$ means that there exist constants $c_3,c_4,c_5,c_6>0$ such that $c_3(g_1(x)+g_2(x)h(c_4x))\le f(x) \le c_5(g_1(x)+g_2(x)h(c_6x))$ for 
specified range of $x$.
We use the convention $1/\infty=0$. 

%%%%%%%%%%%%%%%%%%%%%%%%%%%%%%%%%%%%%%%%%%%%%%%%%%%%%%%%%%%%%%%%%%%%%%%%%%%%%%%%%%%%%%
%%%%%%%%%%%%%%%%                                                             Preliminaries                                               %%%%%%%%%%%%%%%%%%%%%%%%%% %%
%%%%%%%%%%%%%%%%%%%%%%%%%%%%%%%%%%%%%%%%%%%%%%%%%%%%%%%%%%%%%%%%%%%%%%%%%%%%%%%%%%%%%% 
\section{Estimates on distributions of  subordinators }\label{s:preliminaries}

Let $S=(S_t)_{t\ge0}$ be a driftless subordinator (i.e., a non-decreasing pure-jump L\'evy process on $\R$ with $S_0=0$) with Laplace exponent $\phi$ given by
\begin{equation*}
\phi(\lambda) = -\log \E e^{-\lambda S_1}=\int_0^\infty (1-e^{-\lambda s}) \nu(ds).
\end{equation*}
Let $w(r):=\nu(r, \infty)$. 
 Using that $\phi(\lambda)=\lambda\int_0^{\infty}e^{-\lambda s}w(s)ds$, it is easy to  see (cf. the proof of \cite[Lemma 2.1]{CK1}) 
 that 
\begin{equation}\label{e:phi-w}
	e^{-1}\lambda \int_0^{1/\lambda}w(s)ds \le \phi(\lambda) \le 
	2\lambda \int_0^{1/\lambda}w(s)ds, \quad \lambda>0.
\end{equation} 
The following is our main assumption on the subordinator $S$.

\smallskip

\noindent{\bf (Poly-$R_1$)} There exist constants $R_1 \in (0, \infty]$, $c_1,c_2>0$ and $\beta_2 \ge \beta_1>0$ such that
\begin{align*}
c_1 \Big(\frac{R}{r}\Big)^{\beta_1} \le \frac{w(r)}{w(R)} \le c_2 \Big(\frac{R}{r}\Big)^{\beta_2}    \quad \text{for all} \;\; 0<r \le R  <R_1.
\end{align*}

Suppose 
that  {\bf (Poly-$R_1$)} holds. Then 
by \cite[Lemma 2.1(ii)]{CK1},
in case $R_1<\infty$, for every $r_0>0$, there exists 
$c_3=c_3(r_0)>0$ such that
\begin{equation}\label{e:phi-upper-scaling}
\frac{\phi(R)}{\phi(r)}\le 
c_3\Big(\frac{R}{r}\Big)^{\beta_2\wedge 1}, \quad r_0<r<R.
\end{equation}
On the other hand, by adapting the proof of \cite[Lemma 2.3(3)]{CK20}, 
we can get that, for every $r_0>0$, 
there exists $c_4=c_4(r_0)>0$ such that 
\begin{equation}\label{e:phi-lower-scaling}
\frac{\phi(R)}{\phi(r)}\ge 
c_4\Big(\frac{R}{r}\Big)^{\beta_1}, \quad r_0<r<R.
\end{equation}
As a consequence of \eqref{e:phi-upper-scaling} and \eqref{e:phi-lower-scaling}, we see that 
$\beta_1\le 1$ and that
$\phi^{-1}$ enjoys the following scaling: For every $t_0>0$, there exist 
$c_5, c_6>0$ depending on $t_0$ such that
\begin{equation}\label{e:phi-inv-scaling}
c_5\Big(\frac{t}{s}\Big)^{1/(\beta_2\wedge 1)}\le \frac{\phi^{-1}(t)}{\phi^{-1}(s)}\le c_6 \Big(\frac{t}{s}\Big)^{1/\beta_1}, \quad t_0<s<t.
\end{equation}
 In case when {\bf (Poly-$\infty$)} holds, \eqref{e:phi-upper-scaling} and  \eqref{e:phi-lower-scaling} are valid for all $0<r<R$, and \eqref{e:phi-inv-scaling} is valid for all $0<s<t$.

\begin{lemma}\label{l:phi-phi'}
Assume {\bf (Poly-$R_1$)} holds. For any $a>0$, there exists  $c_1=c_1(a)\in(0,1)$  such that
	\begin{equation}\label{e:phi-phi'}
c_1\phi(\lambda) \le \lambda\phi'(\lambda) \le \phi(\lambda), \quad \lambda>a.
	\end{equation}
	Further, if  {\bf (Poly-$\infty$)} holds,  then \eqref{e:phi-phi'} holds for all $\lambda>0$.
\end{lemma}
\pf  
The second inequality follows from the fact 
$1-e^{-u}-u e^{-u} \ge 0$ for  $u \ge 0$. 
The first  inequality follows from \cite[Lemma 1.3]{KM14} and  its proof. \qed

\begin{lemma}\label{l:phi-derivatives}
Assume {\bf (Poly-$R_1$)} holds. For any $a>0$, there exists  $c_1=c_1(a)>0$  such that
\begin{equation}\label{e:phi-derivatives}
|\phi''(\lambda)|\le c_1\lambda^{-1}\phi'(\lambda), \quad \lambda>a.
\end{equation}
Further, if  {\bf (Poly-$\infty$)} holds,  then \eqref{e:phi-derivatives} holds for all $\lambda>0$.
\end{lemma}
\pf 
The proof is similar to 
that of \cite[Lemma 2.1(3)]{CK20}, where the existence of L\'evy density  is assumed.
We give a detailed proof for the reader's convenience.

Since $e^{-x} \le x^{-2}$ for all $x>0$, we see that for all $\lambda > 1/R_1$,
\begin{equation}\label{e:sigma_2}
\lambda|\phi''(\lambda)| = \int_{0}^{1/\lambda} \lambda y^2e^{-\lambda y} \nu(dy) + \int_{1/\lambda}^\infty \lambda y^2 e^{-\lambda y} \nu(dy)  \le \int_0^{1/\lambda} y \nu(dy) + \lambda^{-1}w(1/\lambda).
\end{equation}
By {\bf (Poly-$R_1$)}, there exists $\eps \in (0,1/2)$ such that $w(\eps/\lambda) \ge 2w(1/\lambda)$ for all $\lambda>1/R_1$. 
Hence,
\begin{equation*}
\int_0^{1/\lambda} y \nu(dy) \ge \eps \lambda^{-1} \int_{\eps/\lambda}^{1/\lambda} \nu(dy) \ge \eps \lambda^{-1}w(1/\lambda) \quad \text{for all} \;\; \lambda>1/R_1.
\end{equation*}
It follows that 
\begin{equation}\label{e:sigma_3}
\phi'(\lambda) \ge e^{-1}\int_0^{1/\lambda} y \nu(dy) \ge (2e)^{-1} \eps \Big(\int_0^{1/\lambda} y \nu(dy) + \lambda^{-1}w(1/\lambda) \Big).
\end{equation}
Combining \eqref{e:sigma_3} with \eqref{e:sigma_2}, we get that 
 in the case $R_1=\infty$, \eqref{e:phi-derivatives} holds for all $\lambda>0$ with $c_1=\eps/(2e)$, and 
in the case $R_1<\infty$, \eqref{e:phi-derivatives} holds with
$$
c_1=\frac{\eps}{2e}\vee \sup_{ \lambda \in [a,1/R_1]}\big( \lambda|\phi''(\lambda)|/\phi'(\lambda)\big).
$$
\qed

Let $H:(0,\infty)\to (0,\infty)$ be defined by
$H(\lambda):=\phi(\lambda)-\lambda\phi'(\lambda)$, $\lambda>0$.
The function $H$ is strictly increasing, $H(0+)=0$, $\lim_{\lambda\to \infty}H(\lambda)=\int_0^{\infty}\nu(ds)=w(0+)$, 
and satisfies
\begin{equation}\label{e:H-upper-2}
\frac{H(\lambda)}{\lambda^2}=-\Big( \frac{\phi(\lambda)}{\lambda}\Big)'=\int_0^{\infty}e^{-\lambda s}s w(s)ds, \quad \lambda >0.
\end{equation}
Since $1-e^{-\lambda s} -\lambda se^{-\lambda}\ge 1-2/e$ when $s\ge 1/\lambda$, we see that
\begin{equation}\label{e:H-w}
 \phi(\lambda) \ge  H(\lambda)\ge \int_{1/\lambda}^{\infty}(1-e^{-\lambda s}-\lambda s e^{-\lambda s})\nu(ds)\ge \frac{e-2}{e}w(1/\lambda), \quad \lambda>0.
\end{equation}
Suppose that  {\bf (Poly-$R_1$)} holds. Then it follows from the proof of 
\cite[Lemma 2.3(3)]{CK20} that, for every $r_0>0$, there exists a constant
 $c=c(r_0)>0$ such that
\begin{equation}\label{e:H-lower-scaling}
\frac{H(R)}{H(r)}\ge c\Big(\frac{R}{r}\Big)^{\beta_1}, \quad r_0<r\le R.
\end{equation}
As a consequence of \eqref{e:H-lower-scaling},  
we have the following upper scaling for the inverse function $H^{-1}$:
\begin{equation}\label{e:H-inv-scaling}
\frac{H^{-1}(t)}{H^{-1}(s)}\le c^{-1/\beta_1}\Big(\frac{t}{s}\Big)^{1/\beta_1}, \quad H(r_0)<s\le t.
\end{equation}
In case when  {\bf (Poly-$\infty$)} holds, 
  \eqref{e:H-lower-scaling} and \eqref{e:H-inv-scaling} hold with $r_0=0$.  
Note 
that \eqref{e:H-lower-scaling} implies that $\lim_{\lambda\to \infty}H(\lambda)=+\infty$. 

Next we look at the function $b:(0,\infty)\to (0,\infty)$ defined by
$$
b(t):=(\phi'\circ H^{-1})(1/t)=\int_0^{\infty}s e^{-H^{-1}(1/t)s}\nu(ds), \quad t>0.
$$
The function $b$ is strictly increasing, $b(0+)=0$, and $\lim_{t\to \infty}b(t)=\int_0^{\infty}s\nu(ds)=\phi'(0+)$. This implies that $t\mapsto t b(t)$ is also strictly increasing and $\lim_{t\to \infty}t b(t)=+\infty$. Moreover, according to \cite[Lemma 2.4(ii)]{CK1}, cf.~also \cite[(2.13)]{CK20}, it holds that
\begin{equation}\label{e:comp-b-phi0}
\phi^{-1}(7/t)^{-1}\le t b(t) \le \phi^{-1}(1/t)^{-1} \quad \textrm{for all } t>0.
\end{equation}
Hence, under {\bf (Poly-$R_1$)}, we see from the scaling of $\phi^{-1}$ in \eqref{e:phi-inv-scaling} that, for every $t_0>0$,   there exists $c_1=c_1(t_0)>0$ such that
\begin{equation}\label{e:comp-b-phi}
c_1\phi^{-1}(1/t)^{-1}\le t b(t) \le \phi^{-1}(1/t)^{-1}\quad \textrm{for all } 0<t< t_0.
\end{equation}
Moreover, if  {\bf (Poly-$\infty$)}  holds, then \eqref{e:comp-b-phi} holds with $t_0=\infty$.

Finally, we introduce the function $\sigma=\sigma(t,s):(0,\infty) \times (0,\infty) \to [0,\infty)$ defined by 
$$
\sigma=\sigma(t,s):=
(\phi')^{-1}(s/t){\bf 1}_{(0, \phi'(0+))}(s/t), \quad s, t>0.
$$
Note that
$s\mapsto \sigma(t,s)$ is decreasing 
with $\lim_{s\to 0}\sigma(t,s)=\infty$ and $\lim_{s\to \infty}\sigma(t,s)=0$,
 while $t\mapsto \sigma(t,s)$ is increasing 
with $\lim_{t\to 0}\sigma(t,s)=0$ and $\lim_{t\to \infty}\sigma(t,s)=\infty$.
Further, by using the former and the fact that $H$ is increasing, we conclude that 
\begin{equation}\label{e:Hsigma}
t(H\circ\sigma)(t, tb(t))=1 \quad \textrm{ and }\quad t(H\circ\sigma)(t,s)<1 \quad \textrm{ for }s>tb(t).
\end{equation}

The function $\sigma$ plays a crucial role in estimating the left tail of the subordinator $S$. We first state a result which follows from \cite[Lemma 2.11]{CK1} and \cite[Lemma 5.2]{JP}.

\begin{prop}\label{p:lefttail}
	There exist constants $c_1,c_2>0$ such that for all $t>0$,
	$$
	c_1 \exp \big(-c_2 t(H \circ \sigma)(t,s) \big) \le \P(S_t \le s) \le  e \exp \big(- t(H \circ \sigma)(t,s)\big). 
	$$
\end{prop} 
\pf If $s \le tb(t)$, then it follows from \cite[Lemma 2.11]{CK1} and \cite[Lemma 5.2]{JP}
that there exist $c_1,c_2>0$ independent of $s$ and $t$ such that
$
	c_1 \exp \big(-c_2 t(H \circ \sigma)(t,s) \big) \le \P(S_t \le s) \le \exp \big(- t(H \circ \sigma)(t,s) \big).
	$
(Note that the function $b(t)$ in this paper is the same as $t^{-1}b(t)$ in \cite{CK1}.)  
In particular, taking $s=tb(t)$ and using \eqref{e:Hsigma}, we get
$
c_3:=c_1e^{-c_2}\le \P(S_t \le tb(t))\le e^{-1}.
$

If $s > tb(t)$, then by the second part of \eqref{e:Hsigma},
$\exp (-t(H \circ \sigma)(t,s))\ge e^{-1}$, and thus
$
\P(S_t\le s)\le 1 \le e \exp \big(-t (H \circ \sigma)(t,s)\big)
$
which gives the desired upper bound. For the desired lower bound, 
$$
\P(S_t\le s)\ge \P(S_t\le tb(t))  \ge c_3\ge
c_3 \exp  \big(-c_2t (H \circ \sigma)(t,s)\big).
$$
\qed

\begin{lemma}\label{l:sigma}
Suppose {\bf (Poly-$R_1$)}\ holds. Then, for any $a>0$, there exists $\delta=\delta(a)>0$  such that
	\begin{align}\label{e:sigmascale}
\frac{\sigma(t,u)}{\sigma(t,s)} \ge 2^{-\delta}\Big(\frac{s}{u}\Big)^{\delta}, \quad 0<u  \le s \le t\phi'(a).
	\end{align}
	Moreover, if  {\bf (Poly-$\infty$)}\ holds, then \eqref{e:sigmascale} holds for 
	all $0<u \le s <t \phi'(0+)$.
\end{lemma}
\pf  Let $a>0$. For all $0<2w \le t\phi'(a)$, it holds that $\sigma(t,2w) \ge a$. By the mean value theorem, the fact that both $|\phi''|$ and $s\mapsto \sigma(t,s)$ are decreasing,  and Lemma \ref{l:phi-derivatives}, we get 
\begin{align}\label{e:sigma1}
\frac{w}{t} &= (\phi'\circ \sigma)(t,2w) - (\phi'\circ \sigma)(t,w) \le |(\phi'' \circ \sigma)(t,2w)|(\sigma(t,w)-\sigma(t,2w))\nn\\
& \le c_1 \frac{(\phi'\circ \sigma)(t,2w)}{\sigma(t,2w)}(\sigma(t,w)-\sigma(t,2w)) = \frac{2c_1w}{t} \frac{(\sigma(t,w)-\sigma(t,2w))}{\sigma(t,2w)}.
\end{align}
Let  $\delta=\log_2(1+1/(2c_1))$. Then, we see from \eqref{e:sigma1} that for all $0<2w \le t\phi'(a)$,
\begin{equation}\label{e:sigmadoubling}
2^{\delta} \sigma(t,2w) \le \sigma(t,w).
\end{equation}
For any $0<u \le s \le t\phi'(a)$, 
let $n=n(u,s)$ be the largest integer such that $2^n u \le s$. Then, by \eqref{e:sigmadoubling},
 we obtain
\begin{align*}
\frac{\sigma(t,u)}{\sigma(t,s)} \ge 2^{\delta  n}\frac{\sigma(t,2^nu)}{\sigma(t,s)}  \ge 2^{\delta n} \ge 2^{\delta n} 2^{-\delta(n+1)}\Big(\frac{s}{u}\Big)^{\delta}=2^{-\delta}\Big(\frac{s}{u}\Big)^{\delta}.
\end{align*}
This proves the first assertion.

Assume now that $R_1= \infty$. Then \eqref{e:phi-derivatives} is valid for all 
$\lambda>0$ so that \eqref{e:sigma1} holds for all 
$0<2w<t\phi'(0+)$.
Hence,  \eqref{e:sigmadoubling} holds for all 
$0<2w <t\phi'(0+)$. 
We conclude the proof as in the first assertion.  \qed

\begin{lemma}\label{l:lefttail}
	Suppose that {\bf (Poly-$R_1$)}\ holds. Then, for all $\kappa, N>0$ and $T>0$, there exists a constant $C=C(T,\kappa, N)>0$  such that for all $0<t \le T$ and $0<s\le \phi^{-1}(1/t)^{-1}$,
	\begin{align}\label{e:leftexp}
	\exp \big( -\kappa t (H\circ \sigma)(t, s) \big) \le C(s \phi^{-1}(1/t))^{N}.
	\end{align}
	Moreover, if {\bf (Poly-$\infty$)}\ holds, then for all $\kappa ,N>0$, there exists a constant $C=C(\kappa, N)>0$ such that \eqref{e:leftexp} holds for all $0<s \le \phi^{-1}(1/t)^{-1}$.
\end{lemma} 
\pf  
Choose an arbitrary $t\in (0, T]$.
In view of \eqref{e:comp-b-phi}, since $e^{-x} \le 1$ for all $x\ge 0$, it suffices to prove \eqref{e:leftexp} for $0<s \le tb(t)$.  Recall that $t (H\circ \sigma)(t, tb(t)) = 1$. Hence, by \eqref{e:H-lower-scaling}, Lemma \ref{l:sigma} and  \eqref{e:comp-b-phi}, we have that, for all $0<s \le tb(t)$,
\begin{align*}
t (H\circ \sigma)(t, s)= \frac{(H\circ \sigma)(t, s)}{(H\circ \sigma)(t, tb(t))} \ge c_1 \Big(\frac{\sigma(t,s)}{\sigma(t,tb(t))}\Big)^{\beta_1} \ge c_2 \Big(\frac{\phi^{-1}(1/t)^{-1}}{s}\Big)^{\delta\beta_1},
\end{align*}
where $\delta=\delta(T)$ is the constant from Lemma \ref{l:sigma}.  Let $c_3:=\sup_{x>0}x^{N/(\delta \beta_1)} e^{-\kappa x}$. 
Then
\begin{align*}
\exp \big( -\kappa t (H\circ \sigma)(t, s) \big) \le c_3 \Big(\frac{1}{t (H\circ \sigma)(t, s)} \Big)^{N/(\delta\beta_1)} \le c_2 c_3\big(s \phi^{-1}(1/t)\big)^{N}.
\end{align*}
This proves the first assertion. Moreover, we can see that the second assertion is 
true by using that \eqref{e:H-inv-scaling} and \eqref{e:comp-b-phi} hold for
 $r_0=0$ and $t_0=\infty$, respectively and the second assertion of  Lemma \ref{l:sigma}. \qed

\begin{lemma}\label{l:leftint}
Let $f:(0, \infty) \to (0, \infty)$ be a given function. Assume that {\bf (Poly-$R_1$)}\ holds and there exist constants $c_1, p>0$ such that $s^p f(s) \le c_1t^pf(t)$ for all $0< s \le t$. Then for every $T>0$, there exists a constant $C=C(T,c_1, p)>0$ such that for any $t \in (0,  T]$,
\begin{align}\label{e:leftint}
\E[f(S_t): S_t \le r] \le C f(r) \exp \Big(- \frac{t}{2}(H \circ \sigma)(t, r) \Big), \quad 0<r \le \phi^{-1}(1/t)^{-1}.
\end{align}
Moreover, if {\bf (Poly-$\infty$)}\ holds, then there exists a constant $C=C(c_1, p)>0$ such that \eqref{e:leftint} holds for all $t>0$.
\end{lemma}
\pf
By using Proposition \ref{p:lefttail} in the second and Lemma \ref{l:lefttail} (with $\kappa=1/2$ and $N=p+1$) in the third inequality below, we get that
\begin{align*}
&\E[f(S_t): S_t \le r] = \sum_{j=0}^\infty \int_{2^{-j-1}r}^{2^{-j}r} f(s) \P(S_t \in ds) \le c_12^p\sum_{j=0}^\infty f(2^{-j}r) \P(S_t \le 2^{-j}r) \\
& \le c_1^22^pf(r)\sum_{j=0}^\infty 2^{jp}\exp \Big(- \frac{t}{2}(H \circ \sigma)(t,2^{-j}r) \Big)\exp \Big(- \frac{t}{2}(H \circ \sigma)(t,2^{-j}r) \Big)\\
& \le c_2f(r)\exp \Big(- \frac{t}{2}(H \circ \sigma)(t,r) \Big)\sum_{j=0}^\infty 2^{jp}\big(2^{-j}r \phi^{-1}(1/t)\big)^{p+1} \le 2c_2f(r)\exp \Big(- \frac{t}{2}(H \circ \sigma)(t,r) \Big).
\end{align*}
\qed

\begin{prop}\label{p:mode}
	Suppose that {\bf (Poly-$R_1$)}\ holds. Then for any  $T>0$, there exist constants $\delta \in (0,1)$ independent of $T$ and $\epsilon=\epsilon(T) \in (0,1)$ such that 
	\begin{align}\label{e:mode}
	\P\big(\epsilon\phi^{-1}(1/t)^{-1} \le S_t \le \phi^{-1}(1/t)^{-1}\big) \ge \delta, \quad t \in (0,T).
	\end{align}
	Moreover,  if {\bf (Poly-$\infty$)}\ holds, then there exist $\epsilon, \delta \in (0,1)$ such that \eqref{e:mode} holds with $T=\infty$.
\end{prop}
\pf
By Proposition \ref{p:lefttail}, \eqref{e:Hsigma} and \eqref{e:comp-b-phi}, there exists a constant $c_1 \in (0,1)$ such that  $\P(S_t \le \phi^{-1}(1/t)^{-1}) \ge c_1$ for all $t>0$. Let $c_2 = \log(2/c_1)$. Then,  using  Proposition \ref{p:lefttail} again, we get that for all $t>0$,
\begin{equation}\label{e:S-mode}
\P\big(tb(t/c_2) \le S_t \le \phi^{-1}(1/t)^{-1}\big) \ge c_1 - \P\big( S_t < tb(t/c_2)\big) \ge c_1 - e^{-c_2} = c_1/2.
\end{equation}
Since {\bf (Poly-$R_1$)} holds, by \eqref{e:comp-b-phi} and \eqref{e:phi-inv-scaling}, there exists $\epsilon \in (0,1)$ such that 
\begin{equation}\label{e:bdoubling}
tb(t/c_2) \ge \epsilon \phi^{-1}(1/t)^{-1}, \quad t \in (0, T).
\end{equation}
We also see that if {\bf (Poly-$\infty$)} holds, then \eqref{e:bdoubling} holds with $T=\infty$. 
Combining this with \eqref{e:S-mode}, we obtain \eqref{e:mode}.
\qed

\begin{lemma}\label{l:exppoly}
	Suppose that {\bf (Poly-$R_1$)}\ holds. Then, for any $\kappa>0$, there exists a constant $a_1=a_1(\kappa)>0$ such that for all $\phi^{-1}(1/t)^{-1} \le s<R_1/2$, we have
$	\exp\big(-\kappa s H^{-1}(1/t)\big) \le a_1 t w(s).$
\end{lemma}
\pf
According to \cite[Lemma 2.2]{CK1}, 
there exists $c_1>0$ such that 
$
H(1/s)^{1+\beta_2} \le c_1 \phi(1/s)^{ \beta_2} w(s)$, $
s\in (0, R_1/2].
$
Note that by \eqref{e:H-upper-2}, the map $\lambda \to \lambda^{-2}H(\lambda)$ is decreasing.
 Let $c_2:=\sup_{x>0} x^{2(1+\beta_2)}e^{-\kappa x}$. 
Then by using {\bf (Poly-$R_1$)} and the fact that $\phi$ is increasing, we get that for all $\phi^{-1}(1/t)^{-1} \le s < R_1/2$, since $1/s \le \phi^{-1}(1/t) \le H^{-1}(1/t)$,
\begin{align*}
&\exp \big( -\kappa sH^{-1}(1/t)\big) \le  c_2 \Big(\frac{1/s}{H^{-1}(1/t)}\Big)^{2(1+\beta_2)} \le c_2\Big(\frac{H(1/s)}{1/t}\Big)^{ 1+\beta_2} \nn\\
& \le c_1c_2 t^{1+\beta_2} \phi(1/s)^{\beta_2} w(s)  \le c_1 c_2 t^{1+\beta_2} \phi\big(\phi^{-1}(1/t)\big)^{\beta_2} w(s) = c_1 c_2 t w(s).
\end{align*}
This proves the lemma. \qed

\begin{prop}\label{p:righttail}
	 Suppose that {\bf (Poly-$R_1$)}\ holds. Then, for all $2\phi^{-1}(1/t)^{-1} < s < R_1/2$,
	\begin{align}\label{e:tw}
	\P(S_t \ge s) \simeq tw(s).
	\end{align}
	In particular, there exists a constant  $M>1$ such that for all $2\phi^{-1}(1/t)^{-1} < s < R_1/(2M)$,
	\begin{align*}
	\P(S_t \in [s,Ms]) \simeq tw(s).
	\end{align*}
\end{prop}
\pf The lower bound of \eqref{e:tw}  follows from \cite[Lemma 2.6]{CK1} (note that $t\phi(s^{-1})\le 1$). The upper bound of \eqref{e:tw} comes from the proof of \cite[Proposition 2.7]{CK1} with a bit of modification. We provide most of the proof for the reader's convenience. 

Pick an arbitrary $s\in (2\phi^{-1}(1/t)^{-1},  R_1/2)$.
Let $\eps = \log(5/4)/2 \in (0,1)$. We set
\begin{align*}
\mu^1:=\1_{(0, \eps/H^{-1}(1/t)]}\nu(dx), \quad \mu^2:=\1_{(\eps/H^{-1}(1/t), s]}\nu(dx) \;\; \text{and} \;\; \mu^3:=\1_{(s, \infty)}\nu(dx)
\end{align*}
and denote by $S^1, S^2$ and $S^3$ the independent driftless subordinators  with 
L\'evy measures $\mu^1, \mu^2$ and $\mu^3$,  respectively. Then  $S_t \le S^1_t + S^2_t + S^3_t$ (note that it may happen that $s<\epsilon/H^{-1}(1/t)$) and hence
\begin{align*}
\P(S_t \ge s) \le \P(S^1_t \ge {3s}/{4}) + \P(S^2_t \ge {s}/{4}) + \P(S^3_t>0).
\end{align*}

Since $S^3$ is a compound Poisson process, it holds that $\P(S^3_t>0) = 1-e^{-tw(s)} \le tw(s)$. Moreover, by following the proof of \cite[Proposition 2.7]{CK1}, one can obtain from \cite[Proposition 1 and Lemma 9]{KS} that $\P(S^2_t \ge s/4) \le ctw(s)$. Lastly, by using 
Markov's  inequality 
and \cite[Lemma 2.5]{CK1}, since $s > 2tb(t)$ due to \eqref{e:comp-b-phi0}, we have that
\begin{align*}
\P(S^1_t \ge {3s}/{4}) &\le \E \left[ \exp \left(-({3s}/{4})H^{-1}(1/t) + H^{-1}(1/t) S^1_t \right) \right] \\
&= \exp \big( -({3s}/{4})H^{-1}(1/t) + t \int_0^{\eps/H^{-1}(1/t)} (e^{H^{-1}(1/t) x}-1) \nu(dx)\big)\\
&\le \exp \big( -({3s}/{4})H^{-1}(1/t) + e^{2\eps}tH^{-1}(1/t) \int_0^{\eps/H^{-1}(1/t)} xe^{-H^{-1}(1/t) x} \nu(dx)\big)\\
&\le  \exp \big( -({3s}/{4})H^{-1}(1/t) + ({5}/{4})H^{-1}(1/t)tb(t)\big)\le\exp \big( -2^{-3}sH^{-1}(1/t)\big).
\end{align*}
We used the fact that $e^{y}-1 \le ye^{-y}e^{2y}$ for all $y \ge 0$ in the third line. Hence, by Lemma \ref{l:exppoly}, we get that $\P(S^1_t \ge 3s/4) \le ctw(s)$ and hence the first assertion holds.

The second assertion follows from 
{\bf (Poly-$R_1$)}.
\qed

%%%%%%%%%%%%%%%%%%%%%%%%%%%%%%%%%%%%%%%%%%%%%%%%%%%%%%%%%%%%%%%%%%%%%%%%%%%%%%%%%%%
%%%%%%%%%%%%%%%%%%%%%%%%%%            Key estimates and applications                       %%%%%%%%%%%%%%%%%%%%%%%%%%%%%%
%%%%%%%%%%%%%%%%%%%%%%%%%%%%%%%%%%%%%%%%%%%%%%%%%%%%%%%%%%%%%%%%%%%%%%%%%%%%%%%%%%%
 \section{Setup and main assumptions} \label{s:setup}
Let $(E,\rho)$ be a locally compact separable metric space such that all bounded closed sets are compact, and let $m$ a positive Radon measure on  $E$ with full support.
We use $B(x,r)$ to denote the open ball in $(E,\rho)$ of radius $r$ centered at $x$,
 and $V(x,r):=m(B(x,r))$ its volume.

Throughout the remainder of this paper, 
we assume the following volume doubling and reverse volume doubling properties 
with localization radius $R_E\in (0, \infty]$: 
There exist constants $d_2 \ge d_1>0$ such that, for every $a \ge 1$, there exists a constant $C_V=C_V(a) \ge 1$ satisfying
\begin{equation}\label{e:volume_doubling}
 C_V^{-1} \Big(\frac{R}{r}\Big)^{d_1} \le 	\frac{V(x,R)}{V(x,r)} \le C_V \Big(\frac{R}{r}\Big)^{d_2} \quad \text{for all} \; x \in E \;\text{ and }\; 0<r \le R<aR_E.
\end{equation}
As a consequence of \eqref{e:volume_doubling}, we see that for all $R_0, \epsilon,\eta>0$, there exists a constant $C=C(R_0, \epsilon, \eta)>0$ such that
\begin{equation}\label{e:compare_ball}
	V(x,r) \le CV(y,\eta r)  \quad \text{for all} \; x,y \in E \;\text{ and }\;  	\epsilon\rho(x,y)<r\le R_0.\end{equation}
Indeed, since $B(x,r) \subset B(y, r+\rho(x,y))$, we get from \eqref{e:volume_doubling} that  $V(x,r) \le V(y, r+ \rho(x,y)) \le V(y, (1+1/\epsilon)r) \le c_1V(y,\eta r)$.
Moreover, if the localization radius $R_E$ is infinite, 
then the above constant $C$ is independent of $R_0$ and  \eqref{e:compare_ball} holds for 
$\epsilon\rho(x,y)<r<\infty$.

Let $D$ be a proper  open subset of $E$. 
We use $\diam(D)$ to denote the diameter of $D$. If $\diam(D)<\infty$, i.e., $D$ is bounded, then by the assumption on $E$ it holds that $m(D)<\infty$. 
For $x\in D$, let 
$\delta_D(x)=\rho(x, E\setminus D)$. 
In most applications,  $D$  
will be an open subset of the Euclidean space $\R^d$, $d\ge 1$, 
and $m(dy)$ will be the Lebesgue measure. For simplicity we write $dy$ instead of $m(dy)$.
Let $Y^D=(Y^D_t, \P^x)$ be a  Hunt process in $D$.  
We assume that the semigroup of $Y^D$ admits a density 
$p_D(t,x,y)$, which we call 
the heat kernel of  $Y^D$.  
Thus, for any non-negative Borel function $f$ on $D$,
$$
\E^x[f(Y^D_t)]=\int_D f(y) p_D(t,x,y)\, dy.
$$
Let $S=(S_t)_{t\ge 0}$ be a driftless subordinator 
independent of 
$Y^D$.  We will be interested in the subordinate process $X_t:=Y^D_{S_t}$. 
It is well known 
(cf. \cite[p.67, pp. 73--75]{Bo84} 
and \cite{SV08}) that 
$X$   is also a Hunt process and admits a heat kernel $q(t,x,y)$ which is given by the formula
\begin{equation}\label{e:hk-X}
q(t,x,y)=\E[p_D(S_t,x,y)]=\int_0^{\infty}p_D(s,x,y)\P(S_t\in ds).
\end{equation}

Our goal is to find two-sided estimates of 
$q(t,x,y)$ under certain assumptions on the underlying heat kernel $p_D(t,x,y)$ and the subordinator $S$. On the subordinator we will impose the assumption 
{\bf (Poly-$R_1$)}.
Now we explain the assumptions we impose on 
$p_D(t,x,y)$. These assumptions are motivated by various examples from the literature.

We first introduce two functions $\Phi, \Psi:[0,\infty)\to [0,\infty)$, both strictly increasing and satisfying $\Psi(r)\ge \Phi(r)$ for all $r\ge 0$. Moreover, 
we always 
assume that both satisfy global scaling conditions: There exist constants  $\alpha_1,\alpha_2,\alpha_3,\alpha_4>0$ and $c_1,c_2,c_3,c_4>0$ such that for all $R \ge r >0$,
\begin{align}\label{e:psiphi}
 c_1 \Big(\frac{R}{r}\Big)^{\alpha_1} \le \frac{\Phi(R)}{\Phi(r)} \le c_2\Big(\frac{R}{r}\Big)^{\alpha_2}  \quad \text{and} \quad  c_3 \Big(\frac{R}{r}\Big)^{\alpha_3} \le \frac{\Psi(R)}{\Psi(r)} \le c_4 \Big(\frac{R}{r}\Big)^{\alpha_4}.
\end{align}
As an easy consequence we see that, for every $a\ge 1$, there exist two constants $c_1(a)>0$ and $c_2(a)>0$ such that, for all $r, R>0$ satisfying $0<r\le aR$, it holds that
\begin{equation}\label{e:phi-a-scaling}
 c_1(a) \Big(\frac{R}{r}\Big)^{\alpha_1} \le \frac{\Phi(R)}{\Phi(r)} \le c_2(a) \Big(\frac{R}{r}\Big)^{\alpha_2}.
\end{equation}

The following lemma shows that, without loss of generality, we may replace $\Phi$
by a nicer  function.

\begin{lemma}\label{l:derivative-inv-Phi}
There exists a strictly increasing differentiable functions $\wt{\Phi}$ 
satisfying the following two properties:

\noindent(P1) $\Phi(r) \simeq \wt{\Phi}(r)$ for all $r>0$ and $\wt{\Phi}$ satisfies \eqref{e:phi-a-scaling};

\noindent(P2) $\wt{\Phi}'(r)\simeq r^{-1}\wt{\Phi}(r)$ and $(\wt{\Phi}^{-1})'(t) \simeq t^{-1}\wt{\Phi}^{-1}(t)$ for $r,t>0$.
\end{lemma}
\pf  According to \cite[Lemmas 3.1 and 3.2]{CKKW}, for any $\alpha>\alpha_2$, there exists a complete Bernstein function $\varphi$ such that
$$
\Phi(r)\simeq \varphi(r^{-\alpha})^{-1} \quad \textrm{and} \quad \varphi'(r)\simeq r^{-1}\varphi(r) \quad \textrm{for all }r>0,
$$
and that $\varphi$ satisfies the weak scaling conditions with exponents $\alpha_1/\alpha$ and $\alpha_2/\alpha$. Let $\wt{\Phi}(r):=\varphi(r^{-\alpha})^{-1}$, $r>0$. It is straightforward to check that $\wt{\Phi}$ satisfies \eqref{e:phi-a-scaling} and also that $\wt{\Phi}'(r)\simeq r^{-1}\wt{\Phi}(r)$. Moreover, by the inverse function theorem, the second comparability in (P2) is also valid. \qed

\begin{lemma}\label{l:leftint-Phi}
	Let $f:(0, \infty) \to (0, \infty)$ be a given function. Assume that there exist constants $c_1, p>0$ such that $s^p f(s) \le c_1t^pf(t)$ for all $0< s \le t$. Then there exists a constant $c_2=c_2(c_1,p)>0$ such that for all $r, \kappa>0$, 
	\begin{equation}\label{e:leftint-Phi}
	\int_0^r f(s) \exp\Big( - \frac{\kappa^2}{ \Phi^{-1}(s)^2}\Big) ds \le \frac{c_2r^{p+1}f(r)}{\Phi(\kappa)^p}.
	\end{equation}
\end{lemma}
\pf Let $c_3:=\sup_{u>0}u^{p\alpha_2/2}e^{-u}$. Then by the scaling of $\Phi$, we have that
\begin{align*}
	&\int_0^r f(s) \exp\Big( - \frac{\kappa^2}{ \Phi^{-1}(s)^2}\Big) ds \le c_3 	\int_0^{r \wedge \Phi(\kappa)} f(s) \Big(\frac{ \Phi^{-1}(s)}{\kappa}\Big)^{p\alpha_2} ds + 	\int_{r \wedge \Phi(\kappa)}^r f(s) ds\\
	&\le c_4	\int_0^{r \wedge \Phi(\kappa)} f(s) \Big(\frac{ s}{\Phi(\kappa)}\Big)^{p} ds + 	\int_{r \wedge \Phi(\kappa)}^r s^{-p}s^pf(s) ds\\
	&\le \frac{c_1c_4r^pf(r)}{\Phi(\kappa)^p}  \int_0^r ds + \frac{c_1r^pf(r)}{(r \wedge \Phi(\kappa))^p}\int_{r \wedge \Phi(\kappa)}^r ds \le \frac{c_1(c_4+1)r^{p+1}f(r)}{\Phi(\kappa)^p}.
\end{align*}
\qed

\begin{defn}\label{d:df}
{\rm
We say that a function $h: (0,\infty) \times D \times D \to [0,1]$ is a \emph{boundary function} if it satisfies the following two properties:

 \noindent	(H1) For all fixed $x,y \in D$, the map $s \mapsto h(s,x,y)$ is non-increasing.
		
\noindent		(H2) There exist constants 
 $c_1>0$, $\gamma \ge 0$ 
such that
		\begin{align*}
		s^\gamma h(s,x,y) \le c_1 t^\gamma h(t, x, y), 
		 \quad \;\; 0<s \le t < 4\Phi(\diam(D))+1,  \; x,y \in D,
		\end{align*} 	
with $4\Phi(\diam(D))+1$ interpreted as $\infty$ when $D$ is unbounded.

A boundary function $h$ is said to be \emph{regular} if there exists  $c_2>0$ such that
$$
h(t,x,y)\ge c_2, \quad \;\;
 0<t<4\Phi(\diam(D))+1, \ x,y\in D \textrm{ with }
 \delta_\wedge(x, y)\ge \Phi^{-1}(t).
$$

 A regular boundary function $h$ is  said to be of \emph{Harnack-type} if  
 there exists  $c_3>0$ such that  for all $x,y,z\in D$ satisfying $\rho(x,z)\le (\rho(x,y)\wedge \delta_D(x))/2$,
\begin{equation}\label{e:H_3}
h(t,x,y)\le c_3 h(t,z,y), \quad \;\; 0<t< \Phi(\rho(x,y)).
\end{equation}
}
\end{defn}

{\it From now on, $h(t,x,y)$ always denotes a boundary function.}

\begin{remark}\label{r:regular_boundary}
{\rm	Suppose that $h$ is a regular boundary function. Then for every $\eps \in (0,1)$, there exists $c_1=c_1(\eps)>0$ such that 
	$$
	h(t,x,y)\ge c_1, \quad \;\;
	0<t<4\Phi(\diam(D))+1, \ x,y\in D \textrm{ with }
	\delta_\wedge(x, y)\ge \eps \Phi^{-1}(t).
	$$
Indeed, by \eqref{e:psiphi} and  (H2),  we see that for all $x,y\in D$ with 
$\delta_\wedge(x, y)\ge \eps \Phi^{-1}(t)$,
\begin{align*}
	h(t,x,y) \ge c_2 h( \Phi(\eps \Phi^{-1}(t)),x,y) \ge c_3.
\end{align*}}
\end{remark}

\begin{example}\label{ex:boundary-fn}
{\rm  
	(a) Let $p,q\ge0$. For $t>0$ and $x,y\in D$,
define  
\begin{equation}\label{e:def-hp}	h_{p,q}(t,x,y):=\Big(1\wedge \frac{\Phi(\delta_D(x))}{t}\Big)^{p} \Big(1\wedge \frac{\Phi(\delta_D( y ))}{t}\Big)^{q}, \qquad h_p(t,x,y):=h_{p,p}(t,x,y).\end{equation}
	Then  $h_{p,q}(t,x,y)$ is a typical example of a regular boundary function which is also of  Harnack-type.  
	Indeed, (H1) and the regularity is clear, while (H2) holds with $c_1=1$ and $\gamma=p+q$ since for all $0<s<t$, 
	\begin{align*}
		t^{p+q}h_{p,q}(t,x,y)= \Big(t\wedge \Phi(\delta_D(x))\Big)^{p}\Big(t\wedge \Phi(\delta_D(y))\Big)^{q} \ge s^{p+q}
	 h_{p,q}(s,x,y).
	\end{align*}
 Moreover, we see from \eqref{e:psiphi} that for all $x,y,z\in D$ satisfying $\rho(x,z)\le (\rho(x,y)\wedge \delta_D(x))/2$, since $\delta_D(z) \ge \delta_D(x)-\rho(x,z) \ge \delta_D(x)/2$, 
	\begin{equation*}
		\Big(1\wedge \frac{\Phi(\delta_D(x))}{t}\Big)^{p} \le \Big(1\wedge \frac{\Phi(2\delta_D(z))}{t}\Big)^{p} \le c_1 \Big(1\wedge \frac{\Phi(\delta_D(z))}{t}\Big)^{p}.
	\end{equation*}
Thus we conclude that $h_{p,q}(t,x,y)$ is of Harnack-type. 
	The boundary function $h_p(t,x,y)$ is very typical when $D$ is a 
	bounded smooth open subset of $\R^d$.

	\begin{comment}
	(a) Let $p>0$. For $t>0$ and $x,y\in D$, define  
	\begin{equation}\label{e:def-hp}
	h_p(t,x,y):= \left(\frac{\Phi(\delta_D(x))}{\Phi(\delta_D(x))+t} \Big)^{p}  \left( \frac{\Phi(\delta_D(y))}{\Phi(\delta_D(y))+t} \Big)^{p}.
	\end{equation} 
	Then  $h_p(t,x,y)$ is a typical example of a regular boundary function which is also of  Harnack-type.  
	Indeed, (H1) is clear, while (H2) holds with $c_1=1$ and $\gamma=2p$ since for all $0<s<t$, 
	\begin{align*}
	t^{2p}h_p(t,x,y)=\left(\frac{1}{\Phi(\delta_D(x))}  + \frac{1}{t}\Big)^{-p}  \left(\frac{1}{\Phi(\delta_D(y))}  + \frac{1}{t}\Big)^{-p} > s^{2p}h_p(s,x,y).
	\end{align*}
	Further, since $(1/2)(1\wedge u/t)\le u/(u+t)\le 1\wedge u/t$ for all $u,t>0$, it holds that 
	\begin{equation}\label{e:defa-alt}
	h_p(t,x,y)\simeq \Big(1\wedge \frac{\Phi(\delta_D(x))}{t}\Big)^{p}\Big(1\wedge \frac{\Phi(\delta_D(y))}{t}\Big)^{p}.
	\end{equation}
	This clearly shows that $h_p(t,x,y)$ is regular. Lastly, we see from \eqref{e:psiphi} that for all $x,y,z\in D$ satisfying $\rho(x,z)\le (\rho(x,y)\wedge \delta_D(x))/2$, since $\delta_D(z) \ge \delta_D(x)-\rho(x,z) \ge \delta_D(x)/2$, 
	\begin{equation*}
	\Big(1\wedge \frac{\Phi(\delta_D(x))}{t}\Big)^{p} \le \Big(1\wedge \frac{\Phi(2\delta_D(z))}{t}\Big)^{p} \le c_1 \Big(1\wedge \frac{\Phi(\delta_D(z))}{t}\Big)^{p}.
	\end{equation*}
	Combining these with \eqref{e:defa-alt}, we conclude that $h_p(t,x,y)$ is of Harnack-type. 
	The boundary function $h_p(t,x,y)$ is very typical when $D$ is a 
	bounded smooth open subset of $\R^d$.
	\end{comment}
\noindent (b)   
Let $h_p(t, x, y)$ be the function defined in \eqref{e:def-hp}. Then
$h_p(t\wedge 1, x, y)$ is also 
a regular boundary function of Harnack-type.
This is a typical boundary function for smooth exterior open sets.

\noindent (c)
A quite general example of a 
boundary function is obtained as follows. Suppose that $Y^D$ admits a dual process $\wh{Y}^D$.
Let $\zeta$ and $\wh{\zeta}$ be the lifetimes of $Y^D$ and $\wh{Y}^D$ respectively. 
Assume  that the survival probabilities $\P^x(\zeta>t)$ and $\P^y(\wh{\zeta}>t)$ satisfy the following doubling property: $\P^x(\zeta>t/2) \simeq \P^x(\zeta>t)$ and $\P^y(\wh{\zeta}>t/2) \simeq \P^y(\wh{\zeta}>t)$ for all $0<t<4\Phi(\diam(D))+1$ and $x,y \in D$. Then $h(t,x,y):=\P^x(\zeta>t)\P^y(\wh{\zeta}>t)$ is a boundary function. Indeed, (H1) is clear, while (H2) follows from 
the doubling property of survival probabilities assumed above.
The survival probabilities usually satisfy the doubling property, see for instance,
\cite[Lemma 2.21]{CKSV} and its proof. In fact, by \cite[Lemma 2.21]{CKSV} and its proof, one can see that the boundary function $h$ above is often regular.

Moreover, the above $h(t,x,y)$ is of Harnack-type if, in addition,
(1) it is regular; (2)  $Y^D$ satisfies the  (interior elliptic) Harnack inequality 
and (3)  there is 
$c_1>0$ such that for all $x \in D$ and $ \Phi(\delta_D(x))<t< \Phi(\diam(D))$,
\begin{align}\label{e:sufficient_H3}
\P^x(\zeta>t) &\simeq \P^x(\zeta> \tau_{U(x,t)} )= \P^x(Y^D_{\tau_{U(x,t)}},
\in D),
\end{align}
where $U(x,t):=B(x,c_1\Phi^{-1}(t)) \cap D$ and
$\tau_{V}=\inf\{t>0: Y^D_t\notin V\}$.

To see this, we fix $x,z \in D$ satisfying $\rho(x,z) \le \delta_D(x)/2$. If $\delta_D(x) \vee \delta_D(z)\ge (c_1 \wedge 2^{-1})\Phi^{-1}(t)$, then we have $\delta_D(x) \wedge \delta_D(z) \ge \delta_D(x) \vee \delta_D(z)-2^{-1}\delta_D(x) \ge 2^{-1}(c_1\wedge 2^{-1})\Phi^{-1}(t)$. By Remark \ref{r:regular_boundary}, it follows that $
1 \ge \P^x(\zeta>t) \wedge \P^z(\zeta>t) \ge h(t,x,x)  \wedge h(t,z,z) \ge c_2.$
Hence, we obtain $h(t,x,y)/h(t,z,y) = \P^x(\zeta>t)/\P^z(\zeta>t) \le 1/c_2$. If $\delta_D(x)\vee \delta_D(z)<(c_1\wedge 2^{-1})\Phi^{-1}(t)$, then $B(x, \delta_D(x)) \subset U(x,t)$ so that  $v \mapsto \P^{v}(Y^D_{\tau_{U(x,t)}} \in D)$ is harmonic in $B(x,\delta_D(x))$ with respect to $Y^D$. Using \eqref{e:sufficient_H3} twice, we see from the  Harnack inequality, \eqref{e:psiphi} and (H2) that
\begin{align*}
	\P^x(\zeta>t) &\le c_3\P^z(Y^D_{\tau_{U(x,t)}} \in D) \le c_3\P^z(Y^D_{\tau_{U(z, \Phi(\Phi^{-1}(t)/2))}} \in D) \le c_4 \P^z(\zeta> \Phi(\Phi^{-1}(t)/2)) \le c_5 \P^z(\zeta>t).
\end{align*}
The second inequality above is valid since $U(z, \Phi(\Phi^{-1}(t)/2)) \subset U(x,t)$. 
Therefore, we obtain \eqref{e:H_3}.

Under the setting and assumptions in   \cite[Section 2]{CKSV} (Assumptions {\bf A} and {\bf U} in \cite{CKSV}),  for the Hunt process $Y$ defined right below \cite[(2.27)]{CKSV} on 
a $\kappa$-fat open set  $D$ with a critical killing potential  $\mu \in \mathbf{K}_1(D)$, by \cite[Lemma 2.21]{CKSV}, we know that the boundary function $h(t,x,y)=\P^x(\zeta>t)\P^y(\wh{\zeta}>t)$ is 
of Harnack-type. (See  \cite[Definition 2.19]{CKSV}  and \cite[Definition 2.12]{CKSV} for the definitions of a $\kappa$-fat open set and the class $\mathbf{K}_1(D)$, respectively.) 
See \cite{BGR, BGR2, CKS} for related work. 
}
\end{example}
For later use, we record the following simple consequence of (H1) and (H2): Let $k>1$ and $s,t>0$ satisfy $k^{-1}s\le t\le ks\le 4\Phi(\mathrm{diam} (D))$. Then for all $x,y\in D$, 
\begin{equation}\label{e:boundary-function-c}
c_1^{-1}k^{-\gamma}h(s,x,y)\le h(t,x,y)\le c_1 k^{\gamma}h(s,x,y),
\end{equation}
where $c_1$  is the constant from (H2).

\begin{defn}\label{d:pD-estimates}
{\rm
 Let $h(t,x,y)$ be a boundary function. 
	
\noindent (a)  We say that  \HKBh \ holds,  if $D$ is bounded and the following estimates hold: (i) there exist $C_0 \in \{0,1\}$ and $c_1, c_2, c_3, c_4>0$  such that for all $(t,x,y) \in (0,1] \times D \times D$,
\begin{align}\label{e:small}
	&	c_1h(t,x,y) \left[\frac{1}{V(x,\Phi^{-1}(t))} \wedge \left( \frac{C_0 t}{V(x,\rho(x,y))\Psi(\rho(x,y))} + \frac{1}{V(x,\Phi^{-1}(t))} \exp\Big(-\frac{c_2\rho(x,y)^2}{\Phi^{-1}(t)^2}\Big)  \right) \right]\nn\\
	& \le	p_D(t,x,y)   \nn\\ &\le c_3 h(t,x,y) \left[ \frac{1}{V(x,\Phi^{-1}(t))} \wedge \left( \frac{C_0t}{V(x,\rho(x,y))\Psi(\rho(x,y))} + \frac{1}{V(x,\Phi^{-1}(t))} \exp\Big(-\frac{c_4\rho(x,y)^2}{\Phi^{-1}(t)^2}\Big)  \right) \right],
\end{align}
and (ii) there exists a constant $\lambda_D>0$ such that for all $(t,x,y) \in [1, \infty) \times D \times D$,
\begin{align}\label{e:bdd}
		p_D(t,x,y) \simeq e^{-\lambda_D t}h(1,x,y).
\end{align}
		
\noindent (b) We say that \HKUh \  holds, if  
the constant $R_E$ in \eqref{e:volume_doubling} is infinite 
and 
\eqref{e:small} holds 
for all $(t,x,y) \in (0,\infty) \times D \times D$. 
}
\end{defn}

By using the function $(1 \wedge \frac{R_1}{10\Phi(\diam(D))})\Phi(r)$ instead of $\Phi(r)$, we may and do assume that $\Phi(\diam(D)) < R_1/8$ whenever  {\bf (Poly-$R_1$)} and \HKBh \  
hold.

\begin{remark}\label{r:HKBh}
{\rm 
One can easily see that if \HKBh \ holds, then for every $T>0$, there exist constants $c_1,c_2,c_3,c_4>0$ such that \eqref{e:small} holds for all $(t,x,y) \in (0, T] \times D \times D$, and \eqref{e:bdd} holds for all $(t,x,y) \in [T, \infty) \times D \times D$.
}
\end{remark}

\begin{remark}\label{r:HKBh2}
	{\rm 
	Note that $a \wedge (b+c) \le (a\wedge b) + (a \wedge c) \le 2(a \wedge (b+c))$ for all $a,b,c>0$. Hence \eqref{e:small} is 
equivalent to the statement that 	
	for all $(t,x,y) \in (0,1] \times D \times D$,
	\begin{align}\label{e:small_mixed}
		p_D(t,x,y) \asymp h(t,x,y) \left[\left(\frac{1}{V(x,\Phi^{-1}(t))} \wedge  \frac{C_0t}{V(x,\rho(x,y))\Psi(\rho(x,y))}\right)+ \frac{1}{V(x,\Phi^{-1}(t))} \exp\Big(-\frac{c\rho(x,y)^2}{\Phi^{-1}(t)^2}\Big)  \right].
	\end{align}
	}
\end{remark}

\begin{example}\label{ex:HK}{\rm
Here are several examples of processes satisfying \HKBh \ 
or \HKUh . 
We will not try to give the most 
 general examples but the reader 
will see from examples below that our setup is general enough to cover almost all known cases.
 In all examples below, the boundary functions are of Harnack type. 

\smallskip

\noindent  (a) Suppose that $D$ is a bounded $C^{1,1}$ open subset of $\R^d$.

 (1) If $D$ is connected and $Y^D$ is 
 the killed Brownian motion in $D$, then  \HKBh\  is satisfied with $C_0=0$, $\Phi(r)=r^2$ and boundary function $h_{1/2}$.
 See \cite{CKP} for a more general example.

(2) 
If $\alpha\in (0, 2)$ and  $Y^D$ is a killed isotropic $\alpha$-stable process in $D$, then \HKBh\ is satisfied
with $\Phi(r)=\Psi(r)=r^\alpha$ and boundary function $h_{1/2}$, cf. \cite{CKS-jems}. More generally,
suppose $\chi$ is a complete Bernstein function satisfying 
global weak scaling conditions with indices
$\beta_1, \beta_2\in (0, 1)$,  $Y$ is a
subordinate Brownian motion in $\R^d$ via an independent subordinator with
Laplace exponent $\chi$, $Y^D$ is 
the part process of $Y$ in $D$.
Then \HKBh\ is satisfied
with $\Phi(r)=\Psi(r)=1/\chi(r^{-2})$ and boundary function $h_{1/2}$, cf. \cite{CKS}.
See \cite{BGR3, GKK, KK} for more general examples.

(3)
If $D$ is connected and $Y$ is the independent sum of isotropic $\alpha$-stable process and Brownian motion, then its part process $Y^D$ in $D$
satisfies \HKBh\ 
with $\Phi(r)=r^2 \wedge r^\alpha$, $\Psi(r)=r^\alpha$ and boundary function $h_{1/2}$, cf. \cite{CKS50}. 
More generally,
suppose $\chi$ is a complete Bernstein function satisfying the conditions in the paragraph above and
$Y$ is the independent sum of Brownian motion and a subordinate Brownian motion 
via a subordinator with Laplace exponent $\chi$, 
then its part process $Y^D$ in 
$D$
satisfies \HKBh\ 
with $\Phi(r)=\Phi_{\chi}(r):=r^2 \wedge (1/\chi(r^{-2}))$, $\Psi(r)=1/\chi(r^{-2})$ and boundary function $h_{1/2}$, cf. \cite{CKS16}.
Note that since $\lim_{\lambda \to \infty} \chi(\lambda)/\lambda=0$ (see \eqref{e:phi-w}), for every $a>0$, there are comparability constants depending on $a$ such that $\Phi_\chi(r) \simeq r^2$ for $r \in (0,a)$. We remark here that the estimates in  \cite[(1.4)]{CKS50} and \cite[(1.14)]{CKS16} 
are comparable to \eqref{e:small_mixed}   since $t \le 1$.

(4) 
Suppose that
$\chi$ is a complete Bernstein function such that the function $\lambda\mapsto\chi(\lambda)-\lambda \chi'(\lambda)$ satisfies 
weak scaling conditions for $\lambda\ge a>0$ with upper index $\delta<2$
and lower index $\gamma>2^{-1}\1_{\{\delta\ge 1\}}$. 
Suppose that $Y$ is a
subordinate Brownian motion in $\R^d$ via an independent subordinator with
Laplace exponent $\chi$, 
$Y^D$ is the part process of $Y$ in $D$.
Then  \HKBh\ is satisfied
with $\Phi(r)=1/\chi(r^{-2})$, $\Psi(r)=1/(\chi(r^{-2})- r^{-2}\chi'(r^{-2}))$ and boundary function $h_{1/2}$, cf. \cite{KM18}.

(5)
Let $\alpha \in (1,2)$ and $Y^D$ be a censored $\alpha$-stable process in $D$. 
Then it follows from \cite{CKS-ptrf} that \HKBh\ is satisfied
with $\Phi(r)=\Psi(r)=r^\alpha$ and boundary function $h_{(\alpha-1)/\alpha}$.

(6)
Let $\alpha \in (0,2)$ and  $Z^D$ be the part process, in $D$, of a reflected  isotropic  $\alpha$-stable process in $\overline{D}$.
For any $q\in [\alpha-1, \alpha)\cap (0, \alpha)$,
let $Y^D$ be the process on $D$ corresponding to the Feynman-Kac semigroup of $Z^D$ via the multiplicative functional $\exp(-\int^t_0C(d, \alpha, q)
\dist(Z^D_s, \partial D)^{-\alpha}ds)$, where the positive constant $C(d, \alpha, q)$ is defined on
\cite[p. 233]{CKSV}. It follows from \cite[Theorem 3.2]{CKSV} that the small
time estimates \eqref{e:small} holds with 
$\Phi(r)=\Psi(r)=r^\alpha$ and $h_{q/\alpha}$.
Using the small time estimates and the argument in \cite[Section 4]{CS}, one can easily show that the semigroup of $Y^D$ is intrinsically ultracontractive. With this, one 
can easily check that the large time estimates in Definition 
\ref{d:pD-estimates}(a)(ii) holds. Thus \HKBh\ holds.

(7) Suppose that $D$ is connected, $d \ge 3$ and $\kappa\ge -\frac14$.
Let
$Y^D$ be the process corresponding to $\Delta|_D-\kappa \delta_D(x)^{-2}$, the Dirichlet Laplacian in $D$ with critical
potential $\kappa \delta_D(x)^{-2}$. It follows from \cite[(6)]{DD} and \cite[Corollary 1.8]{FMT} that the heat kernel of $Y^D$ satisfies  \HKBh\ with $C_0=0$, $\Phi(r)=r^2$ and boundary function $h_{p}$, where
$p=\frac12(\frac12+\sqrt{\frac14+\kappa})$.

(8) 
Suppose that $\alpha \in (1,2)$ and $d \ge 2$.  Let $b:\R^d\to \R^d$ such that $|b|$ is in the Kato class 
$\mathbb{K}_{d, \alpha-1}$ (see \cite[Definition 1.1]{CKS5} for definition). 
Let $Y$ be an $\alpha$-stable process with drift $b$ in $\R^d$, that is, a process
with generator $-(-\Delta)^{\alpha/2}+b\cdot\nabla$, and let $Y^D$ be the part process of $Y$ in  $D$.
By \cite[Theorem 1.3]{CKS5},
\HKBh\ holds with $\Phi(r)=\Psi(r)=r^\alpha$ and $h_{1/2}$. See also
\cite{KSo}.

(9)
For general setups in which \HKBh\ is satisfied, see \cite[Section 2]{CKSV} and \cite{GS}.

\smallskip

\noindent  (b) Suppose that $D$ is an unbounded $C^{1,1}$ open subset of $\R^d$.

(1) 
If $D$ is the domain above the graph of a bounded Lipschitz function in $\R^{d-1}$,
then the killed Brownian motion in $D$ satisfies  \HKUh\ with $C_0=0$, $\Phi(r)=r^2$ and  a boundary function defined in terms of  survival probabilities like in  Example \ref{ex:boundary-fn}(b), 
which is of Harnack type  (cf. \cite{Var}).

(2)
Suppose that $D$ is a half-space-like $C^{1, 1}$ open set in $\R^d$ and  $\alpha\in (0, 2)$.  Let $Y^D$ be the part process in $D$ of an isotropic $\alpha$-stable process.
Then by
\cite[Theorem 1.2]{CT},  \HKUh\ is satisfied
with $\Phi(r)=\Psi(r)=r^\alpha$ and boundary function $h_{1/2}$. 
More generally, let $Y^D$ be the part process  in $D$ of the independent sum of Brownian motion and an  isotropic $\alpha$-stable process.  By  \cite[Theorem 1.4  and Remark 1.5(ii)]{CKS-ejp}, \HKUh\ is satisfied with 
$\Phi(r)=r^2 \wedge r^\alpha$, $\Psi(r)=r^\alpha$ and boundary function $h_{1/2}$.
When $D$ is an exterior $C^{1, 1}$ open set in $\R^d$ with $d>\alpha$ and $Y^D$ 
is part process in $D$ of an isotropic $\alpha$-stable process, it follows from 
\cite[Theorem 1.2]{CT} that \HKUh\ is satisfied
with $\Phi(r)=\Psi(r)=r^\alpha$ and boundary function $h_{1/2}(t\wedge 1, x, y)$.
See \cite{K} for a more general example.

(3)
Suppose $D$ is the upper half space in $\R^d$. Let $\chi$ be a complete Bernstein 
function satisfying  global weak scaling conditions 
with indices $\alpha_1, \alpha_2\in (0, 1)$, $Y$ be a
subordinate Brownian motion in $\R^d$ via an independent subordinator with
Laplace exponent $\chi$, 
$Y^D$ be the part process of $Y$ in $D$.
It follows from \cite[Theorem 5.10]{KSV-spa14} that  \HKUh\ is satisfied
with $\Phi(r)=\Psi(r)=1/\chi(r^{-2})$ and boundary function $h_{1/2}$.
See \cite{CK} for a more general example.

(4)
Suppose that $D$ is the upper half space in  $\R^d$ and 
$\alpha\in (0, 2)$. 
Let $Z^D$ be the part process, in $D$, of a  reflected  isotropic $\alpha$-stable process in $\overline{D}$.
For any 
$q\in [\alpha-1, \alpha)\cap (0, \alpha)$, 
let $Y^D$ be the process on $D$ corresponding to the Feynman-Kac semigroup of $Z^D$ via the multiplicative functional $\exp(-\int^t_0C(d, \alpha, q)
\delta_D (Z^D_s)^{-\alpha}ds)$, 
where $C(d, \alpha, q)$ is defined on
\cite[p. 233]{CKSV}. It follows from \cite[Theorem 3.2]{CKSV} that \HKUh\ is satisfied
with $\Phi(r)=\Psi(r)=r^\alpha$ and 
boundary function $h_{q/\alpha}$.

(5)
Suppose that $D=\R^d\setminus\{0\}$ and $\alpha\in (0, 2)$. Let $Z$ be an isotropic $\alpha$-stable process in $\R^d$.  
For any $q\in (0, \alpha)$, 
let $Y^D$ be the process on $D$ corresponding to the Feynman-Kac semigroup of $Z^D$ via the multiplicative functional $\exp(-\int^t_0\widetilde{C}(d, \alpha, q)
|Z^D_s|^{-\alpha}ds)$, where $\widetilde{C}(d, \alpha, q)$ is defined on
\cite[p. 250]{CKSV}.  It follows from 
\cite[Theorem 3.9]{CKSV}
and \cite[Theorem 1.1]{JW}
 that \HKUh\ is satisfied
with $\Phi(r)=\Psi(r)=r^\alpha$ and 
boundary function $h_{q/\alpha}$.

(6)
Suppose that $D=\R^d \setminus \{0\}$, $d \ge 2$ or $D=(0,\infty)$.  Let $Y^D$ be a  process with generator $\Delta + (a-1) |x|^{-2}\sum_{i,j=1}^d x_ix_j 
\partial_{ij}+ \kappa  |x|^{-2}  \cdot \nabla - b|x|^{-2}$ 
for some $a>0$, $\kappa,b \in \R$ such that
$$
\Lambda:= \frac{1}{2}\sqrt{ \frac{b}{a} + \Big( \frac{d-1+\kappa-a}{2a}\Big)^2 } \ge \frac{1}{4a} \big( (d-1+\kappa-a) \vee ((2a-1)d+1-\kappa-3a)\big).
$$
Note that when $a=1$ and $\kappa ,b\ge 0$,  the above inequality is always true.  It follows from \cite[Proposition 4.14, Theorem 6.2, Corollary 6.4]{MNS18}   that \HKUh\ is satisfied
with $C_0=0$, $\Phi(r)=r^2$ and 
boundary function $h_{p,q}$ where $p=\Lambda- (d-1+\kappa-a)/(4a)$ and $q=\Lambda- ((2a-1)d+1-\kappa-3a)/(4a)$.

(7)
Suppose that $\alpha \in (1,2)$ and  $D=\R^d \setminus \{0\}$, $d \ge 3$. Let $Y^D$ be a  process with generator $-(-\Delta)^{-\alpha/2} + \kappa |x|^{-\alpha} x \cdot \nabla$ for some $\kappa \in (0,\infty)$. It follows from \cite[Theorems 4 and 5]{KSe20}   that \HKUh\ is satisfied
with $\Phi(r)=\Psi(r)=r^\alpha$ and 
boundary function $h=h_{0,\beta/\alpha}$ for $\beta \in (0,\alpha)$ determined by the equation 
at the beginning of  \cite[Section 3.2]{KSe20}.

}
\end{example}

We now briefly discuss the term
$$
I(t,x,y, C_0):= \frac{1}{V(x,\Phi^{-1}(t))}  \wedge \left( \frac{C_0 t}{V(x,\rho(x,y))\Psi(\rho(x,y))} + \frac{1}{V(x,\Phi^{-1}(t))}  
\exp\Big(-\frac{c_1\rho(x,y)^2}{\Phi^{-1}(t)^2}\Big)  \right)
$$
appearing in \eqref{e:small}. If $C_0=0$, then clearly
\begin{equation}\label{e:interior-estimate-0}
I(t,x,y,0)=\frac{1}{V(x,\Phi^{-1}(t))}  
\exp\Big(-\frac{c_1\rho(x,y)^2}{\Phi^{-1}(t)^2}\Big).
\end{equation}
Suppose now that $C_0=1$. 
\begin{lemma}
For any $a\ge 1$, there are comparability constants depending on $a$ such that 
	\begin{align}\label{e:interior-estimate-1}
	&I(t,x,y,1)\asymp 
	\begin{cases}
	\displaystyle\frac{1}{V(x,\Phi^{-1}(t))} , & t\ge a^{-1}\Phi(\rho(x,y)),\\
	\displaystyle \frac{ t}{V(x,\rho(x,y))\Psi(\rho(x,y))} + \frac{1}{V(x,\Phi^{-1}(t))}  
	\exp\Big(-\frac{c\rho(x,y)^2}{\Phi^{-1}(t)^2}\Big),
	&  t<a\Phi(\rho(x,y)). 
	\end{cases}	
	\end{align}
In particular, if $\Psi(r) \simeq \Phi(r)$ for $r \in (0,R_1)$, then
\begin{equation}\label{e:Phi-approx-Psi}
I(t,x,y,1)\simeq\frac{1}{V(x,\Phi^{-1}(t))}  \wedge   \frac{t}{V(x,\rho(x,y))\Phi(\rho(x,y))}, \quad t>0, \; x,y \in D, \; \rho(x,y)<R_1.
\end{equation}
\end{lemma}
\pf 
If $t \ge a^{-1}\Phi(\rho(x,y))$, then by \eqref{e:phi-a-scaling}, 
\begin{equation*}
	\frac{1}{V(x,\Phi^{-1}(t))} \ge I(t,x,y,1) \ge \frac{1}{V(x,\Phi^{-1}(t))} \exp \Big(- \frac{c_1\rho(x,y)^2}{\Phi^{-1}(t)^2} \Big) \ge  \frac{c_2(a,c_1)}{V(x,\Phi^{-1}(t))}.
\end{equation*}
Assume that $t<a\Phi(\rho(x,y))$. Set
$$
g(t, x, y):=\frac{ t}{V(x,\rho(x,y))\Psi(\rho(x,y))} + \frac{1}{V(x,\Phi^{-1}(t))} \exp\Big(-\frac{c_1\rho(x,y)^2}{\Phi^{-1}(t)^2}\Big).
$$
Clearly, $I(t,x,y,1) \le g(t,x,y)$. Further, by using that $\Psi\ge \Phi$ and \eqref{e:phi-a-scaling}, we have
\begin{equation*}
	g(t,x,y) \le \frac{a}{V(x, \Phi^{-1}(t/a))}+ \frac{1}{V(x, \Phi^{-1}(t))}\le  \frac{c_3(a)+1}{V(x, \Phi^{-1}(t))}.
\end{equation*}
Hence, $I(t,x,y,1) = V(x,\Phi^{-1}(t))^{-1} \wedge g(t,x,y) \ge (c_3(a)+1)^{-1}g(t,x,y)$. Thus,  \eqref{e:interior-estimate-1} holds.

Now, we assume that $\Psi(r) \simeq \Phi(r)$ for $r\in(0, R_1)$. Using \eqref{e:volume_doubling},  \eqref{e:phi-a-scaling} and the fact that  $e^{-u} \le k^ku^{-k}$ for all $u,k>0$,  we get that for all $t>0$ and $x,y\in D$ satisfying $t<\Phi(\rho(x,y))$ and $\rho(x,y)<R_1$,
\begin{align*}
\frac{1}{V(x,\Phi^{-1}(t))} \exp \Big(- \frac{c_1\rho(x,y)^2}{\Phi^{-1}(t)^2}\Big) &\le \frac{c_4}{V(x,\Phi^{-1}(t))} \Big( \frac{\Phi^{-1}(t)^2}{c_1\rho(x,y)^2}\Big)^{(d_2+\alpha_1)/2}\\
& \le \frac{c_5t}{V(x,\rho(x,y))\Phi(\rho(x,y))}  \le \frac{c_6t}{V(x,\rho(x,y))\Psi(\rho(x,y))}.
\end{align*}
Thus, we can deduce \eqref{e:Phi-approx-Psi} from \eqref{e:interior-estimate-1}. \qed

\section{Jump kernel and heat kernel estimates}\label{s:key-estimates}

For a given boundary function $h$, we define for $(t,x,y) \in [0,\infty) \times D \times D$,
\begin{equation}\label{e:defsBstar}
\sB^*_h(x,y):= \int_{0}^{\Phi(\rho(x,y))} h(s,x,y)w(s)ds 
\end{equation} 
and if $\phi^{-1}(1/t)^{-1}\le \Phi(\rho(x,y))$,
\begin{equation}\label{e:defsB}
\sB_h(t,x,y):= \int_{2\phi^{-1}(1/t)^{-1}}^{4\Phi(\rho(x,y))} h(s,x,y)w(s)ds.
\end{equation} 
Since $\int_0^rw(s)ds<\infty$ for all $r>0$ (see \eqref{e:phi-w}) and $h \le 1$, 
the integral in \eqref{e:defsBstar} converges.
Note that, 
by  (H1),
$\sB_h^*(x,y) \simeq \sB_h(0,x,y)$ for all $(x,y) \in D \times D$. 

 \subsection{Jump kernel estimates}\label{ss:jump}

The jump kernel of the subordinate process $X$ is given by
\begin{equation}\label{e:jumping-J}
J(x,y)=\int_0^{\infty}p_D(s,x,y) \nu(ds)\,,  \quad x,y\in D.
\end{equation}
 See \cite[p.74]{Bo84} and also \cite{Miy}. 

\begin{thm}\label{t:jump-estimate}
Suppose that  either 
(1) {\bf (Poly-$R_1$)} and \HKBh \ hold, or  (2)  {\bf (Poly-$\infty$)} and \HKUh \ hold. 
	Then, for $(x,y) \in D\times D$ with $x \neq y$, 
	\begin{equation}\label{e:jump-estimate}
	J(x,y) \simeq \frac{C_0\sB_h^*(x,y)}{V(x,\rho(x,y))\Psi(\rho(x,y))} +  	h(\Phi(\rho(x,y)), x, y)\frac{ w\big(\Phi(\rho(x,y))\big)}{V(x,\rho(x,y))}.
	\end{equation}
\end{thm} 
\pf  Since the proofs are similar, we only give the proof 
of the case (1), 
which is more complicated.
 Fix $x,y\in D$ with $x\neq y$ and let $r:=\rho(x,y)>0$. 
By Remark \ref{r:HKBh}, \eqref{e:small} and \eqref{e:bdd} 
hold with  $T:=\Phi(2\mathrm{diam}(D))$. Then by \eqref{e:jumping-J} and  \eqref{e:interior-estimate-1},
\begin{align*}
	J(x,y)&
	\asymp
\frac{C_0}{V(x,r)\Psi(r)} \int_0^{\Phi(r)} sh(s,x,y) \nu(ds)   + \int_0^{\Phi(r)}\frac{h(s,x,y)}{V(x,\Phi^{-1}(s))}\exp\Big(-\frac{
	cr^2}{\Phi^{-1}(s)^2}\Big)\nu(ds)\\
	&\quad  +\int_{\Phi(r)}^T\frac{h(s,x,y)}{V(x,\Phi^{-1}(s))} \nu(ds) +h(1,x,y)\int_T^{\infty}e^{-\lambda_D s}\nu(ds)=:C_0J_1+J_2+J_3+J_4.
\end{align*}

Since  {\bf (Poly-$R_1$)} holds, there exists a constant $a>1$ such that $w(s/a) \ge 2w(s)$ for all $s <R_1$. Therefore, by \eqref{e:boundary-function-c}, since we assumed $\Phi(\diam(D))<R_1/8$,
\begin{align*}
&V(x,r)\Psi(r)J_1 = \sum_{i \in \N}\int_{a^{-i}\Phi(r)}^{a^{-i+1}\Phi(r)} sh(s,x,y)\nu(ds) \\
&\simeq  \sum_{i \in \N} a^{-i}\Phi(r) h(a^{-i}\Phi(r),x,y) \big(w ( a^{-i}\Phi(r)) - w ( a^{-i+1}\Phi(r) )  \big)\\
&\simeq  \sum_{i \in \N} a^{-i}\Phi(r) h(a^{-i}\Phi(r),x,y) w ( a^{-i}\Phi(r) )  \simeq \sum_{i \in \N}\int_{a^{-i}\Phi(r)}^{a^{-i+1}\Phi(r)} h(s,x,y)w(s) ds = \sB_h^*(x,y).
\end{align*}
Next, by (H1), the scaling  and monotonicity of $\Phi$,
we get that
\begin{align*}
J_2 &\ge \frac{h(\Phi(r),x,y)}{V(x,r)}\int_{\Phi(r)/a}^{\Phi(r)} \exp\Big(-\frac{c_1r^2}{\Phi^{-1}(s)^2}\Big)\nu(ds)\ge \frac{c_2 h(\Phi(r),x,y)}{V(x,r)} \int_{\Phi(r)/a}^{\Phi(r)}\nu(ds) \\
&= \frac{c_2 h(\Phi(r),x,y)}{V(x,r)} \big( w(\Phi(r)/a) - w(\Phi(r))\big) \ge \frac{c_2 h(\Phi(r),x,y)w (\Phi(r))}{V(x,r)}.
\end{align*} 
Hence, we obtain the lower bound in \eqref{e:jump-estimate}.

Now, we prove the upper bound in \eqref{e:jump-estimate}.  Let $\wt{\Phi}$ be the function in Lemma \ref{l:derivative-inv-Phi}.  Since $s \mapsto V(x,\wt\Phi^{-1}(s))^{-1}$ and $s \mapsto h(s,x,y)$ are non-increasing,
using the Leibniz rule for product, integration by parts 
and the property (P2) of $\wt\Phi^{-1}$ in Lemma \ref{l:derivative-inv-Phi}, we obtain
\begin{align}\label{e:estimate-J-1}
J_2  & \le c\int_{0}^{\Phi(r)} \frac{h(s,x,y)}{V(x,\wt\Phi^{-1}(s))} \exp\Big(-\frac{c_3r^2}{\wt\Phi^{-1}(s)^2}\Big)  \Big(-\frac{d}{ds}w(s)\Big) \nn\\
& \le c\int_{0}^{\Phi(r)} \frac{h(s,x,y)w(s)}{V(x,\wt\Phi^{-1}(s))} \bigg(\frac{d}{ds}\exp\Big(-\frac{c_3r^2}{\wt\Phi^{-1}(s)^2}\Big) \bigg) ds \nn\\
& \le c \int_0^{\Phi(r)}  \frac{h(s,x,y)w(s)}{V(x,\wt\Phi^{-1}(s))} \frac{r^2}{s\wt \Phi^{-1}(s)^2}\exp\Big(-\frac{c_3r^2}{\wt\Phi^{-1}(s)^2}\Big) ds.
\end{align}
In the second inequality above, we used the following: Since $h \le 1$, $e^{-x} \le k^kx^{-k}$ for all $x,k>0$ and $\lim_{s \to 0}sw(s)=0$ (because $w$ is the tail of the L\'evy mesure $\nu$), 
by using \eqref{e:volume_doubling} and the scaling of $\wt \Phi^{-1}$, we have that
\begin{align*}
&\lim_{s \to 0} \frac{h(s,x,y)w(s)}{V(x,\wt\Phi^{-1}(s))} \exp\Big(-\frac{c_3r^2}{\wt\Phi^{-1}(s)^2}\Big) \le c\,  \lim_{s \to 0} \frac{w(s)}{V(x,\wt\Phi^{-1}(s))}\bigg(\frac{\wt\Phi^{-1}(s)^2}{r^2}\bigg)^{(d_2+\alpha_2)/2}  \\
&\le \frac{c}{r^{d_2+\alpha_2}V(x,\wt \Phi^{-1}(1))} \lim_{s \to 0} w(s) \wt \Phi^{-1}(s)^{\alpha_2}\le  \frac{c\,\wt \Phi^{-1}(1)^{\alpha_2}}{r^{d_2+\alpha_2}V(x,\wt \Phi^{-1}(1))}  \lim_{s \to 0} sw(s) =0.  
\end{align*}
By {\bf (Poly-$R_1$)}, (H2), \eqref{e:volume_doubling},  \eqref{e:phi-a-scaling} and the fact that $\Phi \simeq \wt \Phi$, we can use Lemma \ref{l:leftint-Phi} with $f(s)=h(s,x,y)w(s)$ $V(x,\wt\Phi^{-1}(s))^{-1}s^{-1}\wt\Phi^{-1}(s))^{-2}$ and $p=\gamma+\beta_2+1+(d_2+2)/\alpha_1$ to deduce from  \eqref{e:estimate-J-1} that
\begin{align}\label{e:estimate-J-11}
J_2 &\le c  \frac{\Phi(r)^{\gamma+\beta_2+1+(d_2+2)/\alpha_1+1}}{\Phi(r)^{\gamma+\beta_2+1+(d_2+2)/\alpha_1}} 
\frac{ h(\Phi(r),x,y) w (\Phi(r))r^2}{V(x, r) \Phi(r)r^2}
=  \frac{ch(\Phi(r),x,y)w(\Phi(r))}{V(x,r)}.
\end{align}

For $J_3$ and $J_4$, 
since $s \mapsto V(x,\Phi^{-1}(s))^{-1}$, $s \mapsto h(s,x,y)$ and $s \mapsto w(s)$ are non-increasing, 
we have by the boundedness of $D$ that 
\begin{align*}
J_3 + J_4 &\le \frac{h(\Phi(r),x,y)w(\Phi(r))}{ V(x, r)}  + h(1,x,y)w(T)\le \frac{ch(\Phi(r),x,y)w(\Phi(r))}{ V(x, r)}.
\end{align*}
This completes the proof. \qed

Suppose that $\Psi\simeq \Phi$ and $C_0=1$. Then the first term in \eqref{e:jump-estimate} dominates the second. Indeed, by (H1), 
\begin{eqnarray*}
\sB_h^*(x,y)&=&\int_0^{\Phi(\rho(x,y))}h(s,x,y)w(s)\, ds \ge h(\Phi(\rho(x,y)),x,y)\int_0^{\Phi(\rho(x,y))}w(s)\, ds\\
&\ge &   h(\Phi(\rho(x,y)),x,y)w(\Phi(\rho(x,y)))\Phi(\rho(x,y)).
\end{eqnarray*}
Moreover, if $\beta_2<1$, then according to \cite[Lemma 2.6, Proposition 2.9]{Mi} and  \eqref{e:phi-upper-scaling}, we get that  $w(s) \simeq \phi(1/s)$ for all $0<s<R_1/2$. Therefore, it holds that
$$
J(x,y)\simeq \frac{1}{V(x, \rho(x,y))\Phi(\rho(x,y))}\int_0^{\Phi(\rho(x,y))}h(s,x,y)\phi(1/s)\, ds.
$$
In case the boundary function is equal to $h_{1/2}$, the integral above can be estimated in the same way as in \cite[Lemma 8.1]{KSV20}, 
cf. \cite[(8.4)]{KSV20}. 

Suppose that $C_0=0$. Then 
\begin{equation}\label{e:J-special}
J(x,y)\simeq h(\Phi(\rho(x,y))\frac{w(h(\Phi(\rho(x,y)),x,y)}{V(x, \rho(x,y))}.
\end{equation}
In particular, in the context of   Example \ref{ex:HK}(b-1), 
and assuming $\beta_2<1$, the above formula reduces to \cite[Theorem 4.4.(1)]{KSV19}. Similarly, 
if $D$ is an exterior $C^{1,1}$ domain in $\R^d$, the boundary function is equal to $h_{1/2}(t\wedge1, x,y)$ and $\beta_2<1$, then \eqref{e:J-special} reduces to  \cite[Theorem 4.4.(2)]{KSV19}.

 \subsection{Heat kernel estimates}\label{ss:key}
 
 Let
 \begin{equation}\label{e:defpsi}
 \psi(r):=\frac{1}{\phi(1/\Phi(r))}, \quad r>0.
 \end{equation}
 Since  $\phi$ and $\Phi$ are strictly increasing, $\psi$ is also strictly increasing. Moreover, it follows from \eqref{e:phi-upper-scaling},  \eqref{e:phi-lower-scaling} and \eqref{e:psiphi} that, for every $R_0>0$, there exist $c_1,c_2>0$ such that
 \begin{equation}\label{e:upper-scaling-psi}
 c_1 \Big(\frac{R}{r}\Big)^{\alpha_1 \beta_1} \le  \frac{\psi(R)}{\psi(r)}\le c_2 \Big(\frac{R}{r}\Big)^{\alpha_2(\beta_2\wedge 1)}, \quad 0<r<R< R_0.
 \end{equation}
 In case when {\bf (Poly-$\infty$)} holds, \eqref{e:upper-scaling-psi} is valid with $R_0=\infty$. We note that
 \begin{equation}\label{e:psi-inverse}
 \psi^{-1}(t) = \Phi^{-1}\big(\phi^{-1}(1/t)^{-1}\big), \quad t>0. 
 \end{equation}
 
 Recall the definition of the function $\sB_h(t,x,y)$  from \eqref{e:defsB}.
 
\begin{thm}\label{t:smallgeneral}
	Suppose that {\bf (Poly-$R_1$)}\ and \HKBh \ hold. Then for every $T>0$, the following estimates are valid for all $(t,x,y) \in (0, T] \times D \times D$:

	\noindent	{\rm (i)} If $\psi(\rho(x,y)) \le t$, then
	\begin{align}\label{e:on}
	q(t,x,y) \simeq \frac{h(\phi^{-1}(1/t)^{-1},x,y)}{V(x,\psi^{-1}(t))}.
	\end{align}

	\noindent	{\rm (ii)} If $\psi(\rho(x,y)) \ge t$, then
	\begin{align}\label{e:general-off}
	q(t,x,y)  &
	\asymp
  \frac{C_0}{V(x,\rho(x,y)) \Psi(\rho(x,y))}\bigg(t\sB_h(t,x,y) +  \frac{h(\phi^{-1}(1/t)^{-1},x,y)}{\phi^{-1}(1/t)} \bigg)  \nn\\
	&\quad +  \frac{h(\phi^{-1}(1/t)^{-1},x,y)}{V(x,\psi^{-1}(t))} \exp\Big(-\frac{c\,\rho(x,y)^2}{\psi^{-1}(t)^2}\Big) + h(\Phi(\rho(x,y)), x, y)  \frac{t w\big(\Phi(\rho(x,y))\big)}{V(x,\rho(x,y))}.
	\end{align}
\end{thm}
\pf Take $x,y \in D$ and let $r:=\rho(x,y)$. We start by establishing some relations valid for all $t\in (0,T]$. By Proposition \ref{p:mode}, there exist constants $\delta, \epsilon \in (0,1)$ such that
\begin{align}\label{e:qlower}
q(t,x,y) \ge \delta \inf_{s \in [\epsilon \phi^{-1}(1/t)^{-1}, \phi^{-1}(1/t)^{-1}]}p_D(s,x,y), \quad t \in (0,T].
\end{align} 
On the other hand, by Remark \ref{r:HKBh} (with $T=\Phi(\diam(D))$), \eqref{e:interior-estimate-0}, \eqref{e:interior-estimate-1}, \eqref{e:bdd} and the fact  that $\exp(-cr^2/\Phi^{-1}(s)^2)\simeq 1$ 
when $s>\Phi(r)$,   we see that 
\begin{align}\label{e:decompq}
&q(t,x,y)\nn\\
&	\asymp
 \int_0^{\Phi(r)} h(s,x,y)\bigg(\frac{C_0s}{V(x,r)\Psi(r)} + \frac{1}{V(x,\Phi^{-1}(s))}\exp\Big(-\frac{cr^2}{\Phi^{-1}(s)^2}\Big) \bigg) \P(S_t \in ds) \nn\\
&\quad +  \int_{\Phi(r)}^{\Phi(\diam(D))} \frac{h(s,x,y)}{V(x,\Phi^{-1}(s))}  \P(S_t \in ds)  + h(1,x,y)\int_{\Phi(\diam(D))}^\infty e^{-\lambda_D s}\, \P(S_t \in ds) \nn\\
&	\asymp
 C_0\int_0^{\Phi(r)} \frac{sh(s,x,y)}{V(x,r)\Psi(r)}\P(S_t \in ds) + \int_0^{\Phi(\diam(D))} \frac{h(s,x,y)}{V(x,\Phi^{-1}(s))} \exp\Big(-\frac{cr^2}{\Phi^{-1}(s)^2}\Big)  \P(S_t \in ds)  \nn\\
&\quad  + h(1,x,y)\int_{\Phi(\diam(D))}^\infty e^{-\lambda_D s}\P(S_t \in ds) =:C_0I_1+I_2+I_3.
\end{align}

  (i) Assume that $\psi(r) \le t$.  By Remark \ref{r:HKBh}, \eqref{e:interior-estimate-1}, 
  (H1),
  \eqref{e:volume_doubling}, the scaling of $\Phi^{-1}$ and \eqref{e:psi-inverse},  there exists a constant 
$c_1>0$ such that 
\begin{equation*}
\inf_{s \in [\epsilon \phi^{-1}(1/t)^{-1}, \phi^{-1}(1/t)^{-1}]}p_D(s,x,y) \ge 
c_1\frac{ h(\phi^{-1}(1/t)^{-1},x,y)}{V(x,\psi^{-1}(t))}, \quad t \in (0,T].
\end{equation*}
Hence, the lower bound in \eqref{e:on} follows from \eqref{e:qlower}.

Now, we prove the upper bound in \eqref{e:on}. 
First, using 
Lemma \ref{l:leftint} in the first inequality below, the assumption $\Psi \ge \Phi$ and
Lemma \ref{l:lefttail} with $N=\gamma+d_2/\alpha_1$ in the second, 
(H2), \eqref{e:psi-inverse} and \eqref{e:psiphi} in the third, and \eqref{e:volume_doubling} in the last, we get that
\begin{align*}
I_1 &\le \frac{\Phi(r)h(\Phi(r),x,y)}{V(x,r) \Psi(r)} \exp \Big(- \frac{t}{2} (H \circ \sigma) \big(t, \Phi(r)\big) \Big) \\
& \le c_2
\frac{h(\Phi(r),x,y)}{V(x,r) }\big( \Phi(r) \phi^{-1}(1/t)\big)^{\gamma} \big( \Phi(r) \phi^{-1}(1/t)\big)^{d_2/\alpha_1}\\
& \le c_3
\frac{h(\phi^{-1}(1/t)^{-1},x,y)}{V(x,r) } \Big(\frac{r}{\psi^{-1}(t)}\Big)^{d_2} \le c_4 \frac{h(\phi^{-1}(1/t)^{-1},x,y)}{V(x,\psi^{-1}(t))}.
\end{align*}
Next, we observe that
\begin{align*}
I_2 &\le \int_{0}^{\phi^{-1}(1/t)^{-1}} \frac{h(s,x,y)}{V(x,\Phi^{-1}(s))}\P(S_t \in ds) + \int_{\phi^{-1}(1/t)^{-1}}^{\infty} \frac{h(s,x,y)}{V(x,\Phi^{-1}(s))}\P(S_t \in ds)=:I_{2,1}+I_{2,2}.
\end{align*}
By (H2), \eqref{e:volume_doubling} and \eqref{e:psiphi},  we can apply Lemma \ref{l:leftint} with $f(s)= h(s,x,y)V(x,\Phi^{-1}(s))^{-1}$ and 
$p=\gamma+d_2/\alpha_1$ to get that 
\begin{align*}
I_{2,1} & \le c_5 \frac{h(\phi^{-1}(1/t)^{-1},x,y)}{V(x,\psi^{-1}(t))}\exp \Big(- \frac{t}{2} (H \circ \sigma) (t, \phi^{-1}(1/t)^{-1}) \Big) \le c_5 \frac{h(\phi^{-1}(1/t)^{-1},x,y)}{V(x,\psi^{-1}(t))}.
\end{align*} 
Moreover, we see from (H1), \eqref{e:psi-inverse} and the monotonicity of $\psi^{-1}$ that 
$$I_{2,2} \le   \frac{h(\phi^{-1}(1/t)^{-1},x,y)}{V(x,\psi^{-1}(t))}\P(S_t \ge \phi^{-1}(1/t)^{-1}) \le  \frac{h(\phi^{-1}(1/t)^{-1},x,y)}{V(x,\psi^{-1}(t))}.$$ 
Lastly, by using (H1) and  (H2), since $\phi$ and $\psi$ are increasing and $t \le T$, we have that 
$$
I_3 \le h(1,x,y) \le  c_6  \frac{h(\phi^{-1}(1/t)^{-1},x,y)}{V(x,\psi^{-1}(t))}.
$$
Hence, we obtain the upper bound in \eqref{e:on} from \eqref{e:decompq}.

\medskip

 (ii) Assume that $\psi(r) \ge t$. First we establish the lower bound. From \eqref{e:qlower},  Remark \ref{r:HKBh}, \eqref{e:interior-estimate-1}, \eqref{e:psi-inverse},  (H1), and the scaling and monotonicity of $\psi^{-1}$, we get that
\begin{align}\label{e:q-lower-off}
q(t,x,y) \ge c_7   h(\phi^{-1}(1/t)^{-1},x,y)\left[
\frac{C_0\phi^{-1}(1/t)^{-1} }{V(x,r) \Psi(r)} 
+ \frac{1}{V(x,\psi^{-1}(t))} \exp\Big(-\frac{c_8r^2}{\psi^{-1}(t)^{2}}\Big) \right].
\end{align}
We also see from Remark \ref{r:HKBh} that
\begin{align}\label{e:q-term-sB}
q(t,x,y)\ge\int_{2\phi^{-1}(1/t)^{-1}}^{4\Phi(r)}p_D(s,x,y)\P(S_t \in ds) \ge  \frac{c_9C_0}{V(x,r)\Psi(r)}
\int_{2\phi^{-1}(1/t)^{-1}}^{4\Phi(r)} sh(s,x,y)\P(S_t \in ds),
\end{align}
where in the second inequality we used \eqref{e:interior-estimate-1} and neglected the second term.  Let $M>1$ be the constant in Proposition \ref{p:righttail}. If $\Phi(r)>M\phi^{-1}(1/t)^{-1}$,  then  by (H1), (H2) and Proposition \ref{p:righttail}, since we assumed $\Phi(\diam(D))<R_1/8$, 
\begin{align*}
&\int_{2\phi^{-1}(1/t)^{-1}}^{4\Phi(r)} sh(s,x,y)\P(S_t \in ds) \ge \sum_{\substack{i \in \N\\M^i \le 2\Phi(r)\phi^{-1}(1/t)}}\int_{2M^{i-1}\phi^{-1}(1/t)^{-1}}^{2M^i \phi^{-1}(1/t)^{-1}} sh(s,x,y)\P(S_t \in ds)\\
& \ge c_{10}M^{-1}  t \sum_{\substack{i \in \N\\M^i \le 2\Phi(r)\phi^{-1}(1/t)}} 2M^i \phi^{-1}(1/t)^{-1}  h(2M^{i-1}\phi^{-1}(1/t)^{-1},x,y)  w \big( 2M^{i-1}\phi^{-1}(1/t)^{-1} \big) \\
&\ge c_{10}M^{-1}t \sum_{\substack{i \in \N\\M^i \le 2\Phi(r)\phi^{-1}(1/t)}}\int_{2M^{i-1}\phi^{-1}(1/t)^{-1}}^{2M^i \phi^{-1}(1/t)^{-1}} h(s,x,y)w(s)ds\\
& \ge c_{10}M^{-1}t \int_{2\phi^{-1}(1/t)^{-1}}^{4\Phi(r)/M} h(s,x,y) w(s)ds \ge 
c_{10}M^{-1}t \int_{2\Phi(r)/M}^{4\Phi(r)/M} h(s,x,y) w(s)ds \\
& \ge 2c_{10}M^{-2} t \Phi(r) h(4\Phi(r)/M, x,y) w(4\Phi(r)/M) \ge 
c_{10}2^{-1}M^{-2} t \int_{4\Phi(r)/M}^{4\Phi(r)} h(s,x,y) w(s)ds.
\end{align*}
By the fourth and the last inequalities above, we deduce from \eqref{e:q-term-sB} that
$$
q(t,x,y) \ge \frac{c_9C_0}{2V(x,r)\Psi(r)}
2\int_{2\phi^{-1}(1/t)^{-1}}^{4\Phi(r)} sh(s,x,y)\P(S_t \in ds) \ge  \frac{c_9c_{10}}{2M^2}\frac{C_0t\sB_h(t,x,y)}{V(x,r)\Psi(r)}.
$$
In case when $\Phi(r) \le M\phi^{-1}(1/t)^{-1}$, we see from (H1), \eqref{e:H-w} and \eqref{e:q-lower-off} that
\begin{align*}
\frac{C_0t\sB_h(t,x,y)}{V(x,r)\Psi(r)} &\le \frac{C_0t}{V(x,r)\Psi(r)} \int_{2\phi^{-1}(1/t)^{-1}}^{4M\phi^{-1}(1/t)^{-1}}h(s,x,y) w(s) ds \\
&\le \frac{4eMC_0}{e-2}
\frac{t \phi^{-1}(1/t)^{-1} }{V(x,r)\Psi(r)}
h( \phi^{-1}(1/t)^{-1},x,y) \phi\big( \phi^{-1}(1/t) \big) \le c_{11}\,q(t,x,y).
\end{align*}

Hence, it remains to prove that there exists a constant 
$c_{12}>0$ such that
\begin{align*}
q(t,x,y) \ge c_{12}
h(\Phi(r),x, y) \frac{tw(\Phi(r))}{V(x,r)}.
\end{align*}
Recall that $M>1$ is the constant in Proposition \ref{p:righttail}. By using Proposition \ref{p:righttail}, \eqref{e:interior-estimate-1}, (H1), \eqref{e:volume_doubling} and the scaling of $\Phi$, we get that, if $4\Phi(r)/M>2\phi^{-1}(1/t)^{-1}$, then
\begin{align*}
q(t,x,y) \ge c_{13}\int_{4\Phi(r)/M}^{4\Phi(r)} \frac{h(s,x,y)}{V(x,\Phi^{-1}(s))} \P(S_t \in ds)\ge 
c_{14}h(\Phi(r),x, y)\frac{tw(\Phi(r))}{V(x,r)}.
\end{align*}
If $4\Phi(r)/M \le 2\phi^{-1}(1/t)^{-1}$, then  by \eqref{e:psi-inverse}, the scaling of $\Phi$ and the assumption that $\psi(r) \ge t$, we get $\psi^{-1}(t) \le r\le 
c_{15}\psi^{-1}(t)$ for some $c_{15}>1$. 
By  \eqref{e:psi-inverse}  and \eqref{e:H-w}, since $w$ is non-increasing,
\begin{equation}\label{e:twphi-inv-bdd}
tw(\Phi(r)) \le tw( \phi^{-1}(1/t)^{-1}) \le  e(e-2)^{-1}t\phi\big(\phi^{-1}(1/t)\big) = e(e-2)^{-1}.
\end{equation}
Therefore, by (H1),   \eqref{e:q-lower-off} (neglecting the first term) and \eqref{e:psiphi}, we obtain 
\begin{align*}
q(t,x,y) \ge 
c_{16}\frac{h(\Phi(r),x,y)}{V(x,r)} \exp \Big( - \frac{c_8c_{15}\psi^{-1}(t)^2}{\psi^{-1}(t)^2} \Big)  \ge  c_{17} h(\Phi(r),x, y) \frac{t w(\Phi(r))}{V(x,r)}.
\end{align*}
This completes the proof of the lower bound.

\smallskip

Now we prove the upper bound. Recall \eqref{e:decompq}. Observe that
\begin{align*}
V(x,r)\Psi(r)I_1 &\le \int_{0}^{\phi^{-1}(1/t)^{-1}} sh(s,x,y)\P(S_t \in ds)+ \int_{\phi^{-1}(1/t)^{-1}}^{2\phi^{-1}(1/t)^{-1}} sh(s,x,y) \P(S_t \in ds) \\
& \quad + \int_{2\phi^{-1}(1/t)^{-1}}^{4\Phi(r)} sh(s,x,y) \P(S_t \in ds)=:K_1+K_2+K_3.
\end{align*}
We get from Lemma \ref{l:leftint} that  
$K_1 \le c_{18}\phi^{-1}(1/t)^{-1}h(\phi^{-1}(1/t)^{-1},x,y).$
Next, by (H1), we have
$K_2 \le 2\phi^{-1}(1/t)^{-1}$ $h(\phi^{-1}(1/t)^{-1},x,y)$.
To bound $K_3$, we use integration by parts and Proposition \ref{p:righttail} to obtain
\begin{align*}
K_3 &=\int_{2\phi^{-1}(1/t)^{-1}}^{4\Phi(r)} sh(s,x,y) \frac{d}{ds}\big(-\P(S_t \ge s) \big)\\
&\le 2\phi^{-1}(1/t)^{-1} h(\phi^{-1}(1/t)^{-1},x,y)  + \int_{2\phi^{-1}(1/t)^{-1}}^{4\Phi(r)} h(s,x,y)\P(S_t \ge s) ds\\
& \quad + \int_{2\phi^{-1}(1/t)^{-1}}^{4\Phi(r)} s\P(S_t \ge s) \frac{dh(s,x,y)}{ds} \le 
c_{19}\big(\phi^{-1}(1/t)^{-1} h(\phi^{-1}(1/t)^{-1},x,y)  + t \sB_h(t,x,y) \big).
\end{align*} 
In the second inequality above, we used the fact that $s \mapsto h(s,x,y)$ is non-increasing (so that 
$s \mapsto \frac{d}{ds} h(s,x,y)\le  0$ a.e.).

Now, we estimate $I_2$. We have 
\begin{align*}
I_2 &\le \int_0^{2\phi^{-1}(1/t)^{-1}} \frac{h(s,x,y)}{V(x,\Phi^{-1}(s))} \exp\Big(-
\frac{cr^2}{\Phi^{-1}(s)^2}\Big)  \P(S_t \in ds)\\
& \quad  + \int_{2\phi^{-1}(1/t)^{-1}}^{4\Phi(r)} \frac{h(s,x,y)}{V(x,\Phi^{-1}(s))}\exp\Big(-\frac{
cr^2}{\Phi^{-1}(s)^2}\Big)  \P(S_t \in ds)  + \int_{4\Phi(r)}^{\infty}
 \frac{h(s,x,y)}{V(x,\Phi^{-1}(s))}  \P(S_t \in ds)\\
& =: L_1+ L_2+L_3.
\end{align*}
By applying Lemma \ref{l:leftint}, 
we get from (H1),  (H2), \eqref{e:volume_doubling}, the scaling of $\Phi$ and \eqref{e:psi-inverse}  that 
\begin{align*}
L_1& \le \exp\Big(-\frac{
cr^2}{\Phi^{-1}(2\phi^{-1}(1/t)^{-1})^2}\Big)  
\int_0^{2\phi^{-1}(1/t)^{-1}} \frac{h(s,x,y)}{V(x,\Phi^{-1}(s))}   \P(S_t \in ds)\\
& \le 
c_{20} \frac{h(\phi^{-1}(1/t)^{-1},x,y)}{V(x,\psi^{-1}(t))}\exp\Big(-\frac{c_{21}r^2}{\psi^{-1}(t)^2}\Big).
\end{align*}
 Let $\wt{\Phi}$ be the function in Lemma \ref{l:derivative-inv-Phi}. By using integration by parts and similar calculations to \eqref{e:estimate-J-1}, we get that
\begin{align*}
L_2  & \le 
c_{22}
\int_{2\phi^{-1}(1/t)^{-1}}^{4\Phi(r)} \frac{h(s,x,y)}{V(x,\wt\Phi^{-1}(s))} \exp\Big(
-\frac{c_{23}r^2}{\wt\Phi^{-1}(s)^2}\Big)  
\frac{d}{ds}\big(-\P(S_t \ge s) \big)\\
& \le c _{24}\, 
\bigg[ \frac{h(2\phi^{-1}(1/t)^{-1},x,y)}{V(x,\wt\Phi^{-1}(2\phi^{-1}(1/t)^{-1}))} \exp\Big(
-\frac{c_{23}r^2}{\wt\Phi^{-1}(2\phi^{-1}(1/t)^{-1})^2}\Big)
\P\big(S_t \ge 2\phi^{-1}(1/t)^{-1}\big)\\
& \qquad+\int_{2\phi^{-1}(1/t)^{-1}}^{4\Phi(r)} \frac{h(s,x,y)}{V(x,\wt\Phi^{-1}(s))}\P(S_t \ge s)\frac{r^2}{s\wt\Phi^{-1}(s)^2} \exp
\Big(-\frac{c_{23}r^2}{\wt\Phi^{-1}(s)^2}\Big) 
ds \bigg] =:
c_{25}\,(L_{2,1} + L_{2,2}).
\end{align*}
By (H1), \eqref{e:volume_doubling}, the scaling of $\Phi$ and \eqref{e:psi-inverse}, since $\Phi \simeq \wt\Phi$, we see that
$$
L_{2,1} \le  c_{26}h(\phi^{-1}(1/t)^{-1},x,y)V(x,\psi^{-1}(t))^{-1}\exp\big(-c_{27}r^2/\psi^{-1}(t)^2\big).
$$
Also, by using  Proposition \ref{p:righttail} and repeating the calculation in \eqref{e:estimate-J-11}, we get that
\begin{align*}
L_{2,2}& \le c_{28}
t\int_{2\phi^{-1}(1/t)^{-1}}^{4\Phi(r)} \frac{h(s,x,y)w(s)}{V(x,\wt\Phi^{-1}(s))} \frac{r^2}{s\wt\Phi^{-1}(s)^2}\exp
\Big(-\frac{c_{23}r^2}{\wt\Phi^{-1}(s)^2}\Big) ds\le 
\frac{c_{29}t h(\Phi(r),x,y) w(\Phi(r))}{V(x,r)}.
\end{align*}
By (H1) and  Proposition \ref{p:righttail}, we obtain
\begin{align*}
L_3 \le  h(\Phi(r),x,y)V(x,r)^{-1} \P(S_t \ge 4\Phi(r))\le  
c_{30}th(\Phi(r),x,y) V(x,r)^{-1}w(\Phi(r)).
\end{align*}

Finally, 
we estimate $I_3$. By Proposition \ref{p:righttail}, since $D$ is bounded,  we get from (H1) and (H2) that 
$$
I_3 \le 
c_{31} th(1,x,y)w\big(\Phi(\diam(D))\big) \le c_{32} 
th(\Phi(r),x,y) V(x,r)^{-1}w(\Phi(r)).
$$
This completes the proof. \qed

 If one assumes  \HKUh 
 , instead of \HKBh 
 ,
together with {\bf (Poly-$\infty$)}, the results of Theorem \ref{t:smallgeneral}(i) and (ii)  are valid for all time. The proof is analogous to 
the proof of Theorem \ref{t:smallgeneral}  and hence omitted.
\begin{thm}\label{t:largegeneral} 
Suppose that {\bf (Poly-$\infty$)}  and \HKUh \ hold. Then the 
assertions in Theorem \ref{t:smallgeneral}(i)--(ii)
hold for all $(t,x,y) \in (0, \infty) \times D \times D$.
\end{thm}

When the upper scaling index $\beta_2$ in  {\bf (Poly-$R_1$)} is strictly less than $1$,  we can obtain the following simpler 
form of off-diagonal estimate.

\begin{cor}\label{c:interior-small}
	Suppose that {\bf (Poly-$R_1$)}\ holds with $\beta_2<1$ and  $\Phi(r) \simeq \Psi(r)$ for $r \in (0,R_1)$.

\noindent (i) If \HKBh \ holds, then for every $T>0$, the following estimates hold for all $(t,x,y) \in (0, T] \times D \times D$:
	
		(1) If $\psi(\rho(x,y)) \le t$, then \eqref{e:on} holds.
	
	 (2) 
If $\psi(\rho(x,y)) \ge t$, then 
\begin{align}\label{e:simple-off0}
q(t,x,y) &\simeq 
\frac{t}{V(x,\rho(x,y))} \times
  \begin{cases} \displaystyle 
\frac{h(\Phi(\rho(x,y)),x,y)}{\psi(\rho(x,y))}	
&\text{ when } C_0=0,\\[10pt]
\displaystyle \frac{\sB_h(t,x,y)}{ \Phi(\rho(x,y))}
&\text{ when } C_0=1.
\end{cases}
\end{align}

	\noindent (ii) If $R_1=\infty$ and \HKUh \ holds, then (1) and (2) above hold for all $(t,x,y) \in (0,\infty) \times D \times D$.
\end{cor} 
\pf According to \cite[Lemma 2.6, Proposition 2.9]{Mi} and  \eqref{e:phi-upper-scaling}, since $\beta_2<1$, we get that  $w(s) \simeq \phi(1/s)$ for all $0<s<R_1/2$.

(i)  We only need to deal with the case (2), i.e., the case $\psi(\rho(x,y)) \ge t$. 
We first assume that $C_0=0$. 
Using $w(s) \simeq \phi(1/s)$, we get that $w(\Phi(\rho(x,y))) \simeq 1/\psi(\rho(x,y))$. 
Thus by Theorem \ref{t:smallgeneral}, 
it remains to show that for any given $c_1>0$, there exists $c_2>0$ such that  for all $(t, x,y) \in (0,T] \times D \times D$ satisfying $\psi(\rho(x,y)) \ge t$,
\begin{equation}\label{e:new}
\frac{h(\phi^{-1}(1/t)^{-1},x,y)}{V(x,\psi^{-1}(t))} \exp\Big(-\frac{c_1 \rho(x,y)^2}{\psi^{-1}(t)^2}\Big)\le c_2 \frac{th(\Phi(\rho(x,y)),x,y)}{V(x,\rho(x,y))\psi(\rho(x,y))}.
\end{equation}
Let $c_3:=\sup_{u>0}u^{(d_2+\alpha_2\gamma+\alpha_2\beta_2)/2}e^{-u}$. By \eqref{e:volume_doubling}, \eqref{e:psiphi}, \eqref{e:upper-scaling-psi}, \eqref{e:psi-inverse} and (H2), we obtain that
\begin{align*}
&\frac{h(\phi^{-1}(1/t)^{-1},x,y)}{V(x,\psi^{-1}(t))} \exp\Big(-\frac{c_1 \rho(x,y)^2}{\psi^{-1}(t)^2}\Big) \le c_3\frac{h(\phi^{-1}(1/t)^{-1},x,y)}{V(x,\psi^{-1}(t))} \Big(\frac{\psi^{-1}(t)^2}{c_1 \rho(x,y)^2}\Big)^{(d_2+\alpha_2\gamma + \alpha_2\beta_2)/2}\\
& \le c_4 \frac{h(\phi^{-1}(1/t)^{-1},x,y)}{V(x,\rho(x,y))} \Big(\frac{\Phi(\psi^{-1}(t))}{ \Phi(\rho(x,y))}\Big)^{\gamma}\frac{t}{\psi(\rho(x,y))} \le  c_5 \frac{th(\Phi(\rho(x,y)),x,y)}{V(x,\rho(x,y))\psi(\rho(x,y))}.
\end{align*}

Now, let $C_0=1$. Since $\Phi(r) \simeq \Psi(r)$ for $r \in (0,R_1)$, 
the first term on the right hand side of \eqref{e:general-off} is comparable with
$$
\frac{ t\sB_h(t,x,y)}{V(x,\rho(x,y)) \Phi(\rho(x,y))}+\frac{h(\phi^{-1}(1/t)^{-1},x,y)}{V(x,\rho(x,y)) \Phi(\rho(x,y)){\phi^{-1}(1/t)}}.
$$
By (H2) and {\bf (Poly-$R_1$)}, it holds that for all $(t, x,y) \in (0,T] \times D \times D$ satisfying $\psi(\rho(x,y)) \ge t$,
\begin{align*}
	t\sB_h(t,x,y) &\ge t\int_{2\Phi(\rho(x,y))}^{4\Phi(\rho(x,y))} h(s,x,y)w(s)ds\ge  c_6t h(\Phi(\rho(x,y)),x,y) w \big(\Phi(\rho(x,y))\big)  \Phi(\rho(x,y)).
\end{align*}
Combining this with \eqref{e:new}, using $w(s) \simeq \phi(1/s)$, one can see that the first term on the right hand side of \eqref{e:general-off} dominates the other two terms. Further, by (H2), \eqref{e:H-w}  and \eqref{e:phi-upper-scaling}, 
we see that for all $(t, x,y) \in (0,T] \times D \times D$ satisfying $\psi(\rho(x,y)) \ge t$,
\begin{align*}
t\sB_h(t,x,y) &\ge t\int_{2\phi^{-1}(1/t)^{-1}}^{4\phi^{-1}(1/t)^{-1}} h(s,x,y)w(s)ds \ge c_7t\int_{2\phi^{-1}(1/t)^{-1}}^{4\phi^{-1}(1/t)^{-1}} h(s,x,y)\phi(1/s)ds \\
& \ge  c_8 \phi^{-1}(1/t)^{-1} h(\phi^{-1}(1/t)^{-1},x,y).
\end{align*}
This yields the desired conclusion.

(ii) 
This
can be proved by the same argument as that of (i). We omit the details here. \qed

In the case when $D$ is a bounded $C^{1, 1}$ domain, $Y^D$ is a killed Brownian motion in $D$ and $S$ is an $(\alpha/2)$-stable subordinator, part (i) of the corollary above is equivalent to \cite[Theorem 4.7]{So}.
In the case when $D$ is an exterior $C^{1, 1}$ domain, $Y^D$ is a killed Brownian motion in $D$ and $S$ is an $(\alpha/2)$-stable subordinator, part (ii) of the corollary above corrects \cite[Theorem 4.6]{So}.

For future use, 
we note the following rough upper  estimates on $q(t,x,y)$.
\begin{prop}\label{p:hke-upper-rough}
	(i) Suppose that {\bf (Poly-$R_1$)} and \HKBh \ hold. Then for every $T>0$, there exists a constant $C>0$ such that for all $(t, x, y)\in (0,T]\times D\times D$,
	\begin{equation}\label{e:hke-upper-rough}
	q(t,x,y) \le Ch(\phi^{-1}(1/t)^{-1},x,y) \bigg( \frac{1}{V(x,\psi^{-1}(t))} \wedge \frac{t}{V(x,\rho(x,y))\psi(\rho(x,y))}\bigg).
	\end{equation}
	
	\noindent (ii) Suppose that  {\bf (Poly-$\infty$)} and \HKUh \ hold. Then, there exists a constant $C>0$ such that \eqref{e:hke-upper-rough} holds for all $(t, x, y)\in (0,\infty)\times D\times D$.
\end{prop} 
\pf (i) Take $x,y \in D$ and let $r:=\rho(x,y)$. If $\psi(r) \le t$, then \eqref{e:hke-upper-rough} follows from Theorem \ref{t:smallgeneral}(i). Hence, we assume that $\psi(r) \ge t$ and estimate each term in Theorem \ref{t:smallgeneral}(ii) separately. 

First, by using (H1), the fact that $\Psi \ge \Phi$ and \eqref{e:phi-w}, we get
\begin{align*}
&\frac{t}{V(x,r) \Psi(r)}\int_{2/\phi^{-1}(1/t)}^{4\Phi(r)} h(s,x,y)w(s)ds \le \frac{th(1/\phi^{-1}(1/t)^{-1},x,y)}{V(x,r)\Phi(r)}  \int_0^{4\Phi(r)}w(s)ds\\
& \le 4e\frac{th(1/\phi^{-1}(1/t)^{-1},x,y)}{V(x,r)}  \phi\big(1/(4\Phi(r))\big)\le 4e
\frac{th(\phi^{-1}(1/t)^{-1}, x, y)}{V(x,r)\psi(r)}.
\end{align*}

Next, we note that, since $\phi$ is a Bernstein function, the map $u \mapsto \phi(u)/u$ is decreasing so that $\Phi(r)\phi(1/\Phi(r)) \ge \phi^{-1}(1/t)^{-1}\phi(\phi^{-1}(1/t))$. Hence, since $\Psi \ge \Phi$, it holds that
\begin{align*}
  \frac{h(\phi^{-1}(1/t)^{-1},x,y)}{V(x,r) \Psi(r)\phi^{-1}(1/t)}\le   \frac{h(\phi^{-1}(1/t)^{-1},x,y)}{V(x,r)\psi(r) \phi(\phi^{-1}(1/t))}= 
\frac{th(\phi^{-1}(1/t)^{-1}, x, y)}{V(x,r)\psi(r)}.
\end{align*}

Thirdly, by \eqref{e:volume_doubling} and  \eqref{e:upper-scaling-psi}, we see that for $c_2:=\sup_{u>0} u^{(d_2+ \alpha_2(\beta_2 \wedge 1))/2}e^{-u}$,
\begin{align*}
\frac{h(\phi^{-1}(1/t)^{-1},x,y)}{V(x,\psi^{-1}(t))} \exp\Big(-\frac{c_1r^2}{\psi^{-1}(t)^2}\Big) &\le c_2\frac{h(\phi^{-1}(1/t)^{-1},x,y)}{V(x,\psi^{-1}(t))} \Big(\frac{\psi^{-1}(t)^2}{c_1r^2}\Big)^{(d_2+ \alpha_2(\beta_2 \wedge 1))/2}\\
&\le c_3\frac{th(\phi^{-1}(1/t)^{-1},x,y)}{V(x,r)\psi(r)}.
\end{align*}

Lastly, by (H1) and \eqref{e:H-w},
\begin{align*}
\frac{t h(\Phi(r), x, y)  w(\Phi(r))}{V(x,r)} \le \frac{e}{e-2}\frac{th(\phi^{-1}(1/t)^{-1},x,y)}{V(x,r)\psi(r)}.
\end{align*}

(ii) By using Theorem \ref{t:largegeneral} instead of Theorem \ref{t:smallgeneral}, we obtain the result by repeating the proof of (i). \qed

As a corollary to Theorems \ref{t:smallgeneral} and \ref{t:largegeneral}, we obtain the following interior estimates on $q(t,x,y)$ in case of a regular boundary function.
\begin{cor}\label{c:smallint} Suppose that 
$h(t,x,y)$ is a regular 	boundary function.
	
\noindent
(i) If {\bf (Poly-$R_1$)}\ and \HKBh \ hold, then  for every $T>0$, the following estimates hold for all $(t,x,y) \in (0, T] \times D \times D$ satisfying  
$\delta_\wedge(x, y)\ge \rho(x,y) \vee \psi^{-1}(t)$.
	
{\rm (1)} If  $\psi(\rho(x,y)) \le t$, then \quad
$
\displaystyle	q(t,x,y) \simeq \frac{1}{V(x,\psi^{-1}(t))}.
$
	
{\rm (2)} If $\psi(\rho(x,y)) \ge t$, then 
	\begin{align*}
	q(t,x,y) & 	\asymp
	\frac{C_0}{V(x,\rho(x,y)) \Psi(\rho(x,y))} \bigg( t\int_{2\phi^{-1}(1/t)^{-1}}^{4\Phi(\rho(x,y))} w(s)ds  + \frac{1}{\phi^{-1}(1/t)} \bigg)\nn\\
	&\quad\;\; +  \frac{1}{V(x,\psi^{-1}(t))} \exp\Big(-\frac{c\,\rho(x,y)^2}{\psi^{-1}(t)^2}\Big) + \frac{t w\big(\Phi(\rho(x,y))\big)}{V(x,\rho(x,y))} .
	\end{align*}

\noindent	(ii)  If {\bf (Poly-$\infty$)}\ and \HKUh \ hold,
then (1) and (2) above hold for all $(t,x,y) \in (0, \infty) \times D \times D$ satisfying 
$\delta_\wedge(x, y) \ge \rho(x,y) \vee \psi^{-1}(t)$.
\end{cor}

Now we give the large time estimates for $q(t, x, y)$ under  \HKBh.

\begin{thm}\label{t:Slarge}
	Suppose that {\bf (Poly-$R_1$)}\ and \HKBh \ hold. Then for every $T>0$, 
	\begin{align}\label{e:Slarge}
q(t,x,y) \simeq e^{-t\phi(\lambda_D)}h(1,x,y), \quad (t, x, y)\in [T, \infty)\times D\times D.
	\end{align}
\end{thm}
\pf Fix 
$x,y \in D$ and $s_0\in(0,1)$ 
such that  $(H \circ \sigma)(T,s_0) \ge 2\phi(\lambda_D)+1/T$. 
Since $\lim_{s \to 0} (H \circ \sigma)(T,s) = \infty$, 
such an $s_0$ always exists.
 Then, since $H$ is non-decreasing and $\phi'$ is non-increasing, we see that
 \begin{equation}\label{e:Slarge_Hsigma}
	(H\circ \sigma)(t,s_0) \ge 2\phi(\lambda_D) + 1/T, \quad t \ge T. 
 \end{equation}
 
By (H2), \eqref{e:volume_doubling} and \eqref{e:psiphi},  we can apply Lemma \ref{l:leftint} with $f(s)= h(s,x,y)V(x,\Phi^{-1}(s))^{-1}$.  Using
 Remark \ref{r:HKBh} 
(with $T=s_0$),
Lemma \ref{l:leftint} and \eqref{e:Slarge_Hsigma}, since $\phi$ is the Laplace 
exponent of $S$, we get that, for all $t \ge T$, 
\begin{align*}
q(t,x,y) &\le c_1\int_0^{s_0} \frac{h(s,x,y)}{V(x,\Phi^{-1}(s))} \P(S_t \in ds) + c_1h(1,x,y) \int_{s_0}^\infty e^{-\lambda_D s} \P(S_t \in ds)\\
& \le c_2 \frac{h(s_0,x,y)}{V(x,\Phi^{-1}(s_0))} \exp \Big(-\frac{t}{2}(H \circ \sigma)(t, s_0) \Big) + c_1h(1,x,y) \E[e^{-\lambda_D S_t}] \le c_3 h(1,x,y) e^{-t\phi(\lambda_D)}.
\end{align*}

On the other hand, we also see from Remark \ref{r:HKBh}  that
\begin{align*}
q(t,x,y) &\ge c_4 h(1,x,y)
\int_{s_0}^\infty e^{-\lambda_D s}\P(S_t \in ds)=c_4h(1,x,y)\Big( e^{-t \phi(\lambda_D)} - \int_0^{s_0} e^{-\lambda_D s}\P(S_t \in ds) \Big).
\end{align*}
According to Proposition \ref{p:lefttail} and \eqref{e:Slarge_Hsigma}, it holds that for all $t \ge T$,
\begin{equation*}
\int_0^{s_0}e^{-\lambda_D s}\P(S_t \in ds) \le \P(S_t \le s_0) \le \exp \big( - 2t\phi(\lambda_D) \big).
\end{equation*}
Therefore, we conclude that  for all $t \ge T$,
$
q(t,x,y) \ge c_4(1-e^{-T\phi(\lambda_D)})h(1,x,y)e^{-t \phi(\lambda_D)}.
$
\qed

\section{Green function estimates}\label{ss:green}\label{s:green}
In this section,
we always assume that either (1)  {\bf (Poly-$R_1$)}\ and \HKBh \ hold, or (2) {\bf (Poly-$\infty$)}  and 	\HKUh \ hold. 
The Green function $G_D$ of $X$ is given by
$
G_D(x,y):= \int_{0}^\infty q(t,x,y) dt.
$

As an application of the heat kernel estimates obtained in 
the previous section, 
we can obtain two-sided estimates on the Green function. To this end, we prove a
simple lemma first.

\begin{lemma}\label{l:intscale}
	Let  $f:I\to [0,\infty)$ be a function defined on an interval  $I \subset [0,\infty)$. 
	Assume that there exist constants $c_1, c_2>0$, $p_1, p_2 \in \R$ such that
	\begin{equation}\label{e:scaling-of-f}
		c_1\Big(\frac{r_2}{r_1}\Big)^{p_1}
		\le \frac{f(r_2)}{f(r_1)} \le c_2\Big(\frac{r_2}{r_1}\Big)^{p_2} \quad \text{for all} \;\;  r_1, r_2 \in I, \; 0<r_1 \le r_2.
	\end{equation} 
	For any $a>1$, there exists a constant $c_3>0$ such that for all $r, R \in I$, $ar \le R$,
	\begin{equation}\label{e:intscale_lower}
		\int_{r}^{R} s^{-1}f(s)ds \ge c_3 \big( f(r) + f(R)\big).
	\end{equation}
	
	\noindent	(i) If we assume $p_1>0$, 
	then, for any $a>1$,  $\int_{r}^{R} s^{-1}f(s)ds \simeq f(R)$ for all $r, R \in I$, $ar \le R$, with comparison constants depending on $a$.
	
	\noindent	(ii) If we assume $p_2<0$,
	then, for any $a>1$, $\int_{r}^{R} s^{-1}f(s)ds \simeq f(r)$ for all $r, R \in I$, $0<ar \le R$, with comparison constants depending on $a$.
\end{lemma}
\pf Suppose $a>1$. For all $r, R \in I$, $ar \le R$, 
by \eqref{e:scaling-of-f},
\begin{align*}
	\int_{r}^{R} s^{-1}f(s)ds \ge 
	\int_{r}^{ar} s^{-1}f(s)ds \vee \int_{R/a}^{R} s^{-1}f(s)ds 
	\ge c \big( f(r) \vee f(R) \big).
\end{align*}
Thus we only need to prove the upper bounds in (i) and (ii).

(i) 
By \eqref{e:scaling-of-f}, 
since $p_1>0$, we have that
\begin{align*}
	\int_{r}^{R} s^{-1}f(s)ds = f(R) \int_r^R \frac{f(s)}{sf(R)} ds \le c_1^{-1}f(R) \int_r^R \frac{s^{p_1-1}}{R^{p_1}}ds \le c_1^{-1}p_1^{-1}f(R).
\end{align*}

(ii) Similarly, 
by \eqref{e:scaling-of-f}, 
since $p_2<0$, we have that \begin{align*}
	\int_{r}^{R} s^{-1}f(s)ds = f(r) \int_r^R \frac{f(s)}{sf(r)} ds \le c_2f(r) \int_r^R \frac{s^{p_2-1}}{r^{p_2}}ds \le -c_2p_2^{-1} f(r). 
\end{align*}
\qed

Recall the definition of $\psi$ in \eqref{e:defpsi}. 
The following
 simple observation will be used in the proof of the next proposition and also later.

\begin{lemma}\label{l:newrs}
If $D$ is bounded, then there  exists  $c>0$ such that for any $x, y\in D$,
$$
 \int_{\Phi(\rho(x,y))}^{2\Phi(\diam(D))}  \frac{h(s,x,y)}{sV(x,\Phi^{-1}(s))\phi(1/s)}ds\ge c\left( \frac{h(\Phi(\rho(x, y)),x,y)\psi(\rho(x, y))}{V(x,\rho(x, y))} + h(1,x,y)\right).
$$
\end{lemma}

\pf 
This follows easily from \eqref{e:intscale_lower}, (H1), (H2) and \eqref{e:compare_ball}. 
\qed

The following proposition provides the first and most general estimate of the Green function.

\begin{prop}\label{p:Sgreen}
 It holds that for $x,y \in D$,
	\begin{equation}\label{e:Green_integral_form}
G_D(x,y) \simeq \frac{C_0}{V(x, \rho(x,y))\Psi(\rho(x,y))} \int_0^{\Phi(\rho(x,y))} \frac{h(s,x,y)}{\phi(1/s)}ds+  \int_{\Phi(\rho(x,y))}^{2\Phi(\diam(D))}  \frac{h(s,x,y)}{sV(x,\Phi^{-1}(s))\phi(1/s)}ds.
	\end{equation} 
\end{prop}

\pf  Since the proofs are similar, we only give the proof when \HKBh \ holds, which is more complicated.

Take $x,y \in D$ and let $r:=\rho(x,y)$. Set $T_D:=1/\phi\big(1/(2\Phi(\diam(D)))\big)$. 
By a change of variables and Lemma \ref{l:phi-phi'}, we have that 
\begin{align}\label{e:green_G4}
	\int_{0}^{\psi(r)} \frac{h(\phi^{-1}(1/t)^{-1},x,y)}{\phi^{-1}(1/t)}  dt = \int_{0}^{\Phi(r)} \frac{h(s,x,y)\phi'(1/s)}{s \phi(1/s)^2}ds \simeq \int_{0}^{\Phi(r)} \frac{h(s,x,y)}{\phi(1/s)}ds
\end{align}
and
\begin{align}\label{e:green_G5}
\int_{\psi(r)}^{T_D} \frac{h(\phi^{-1}(1/t)^{-1},x,y)}{V(x,\psi^{-1}(t))}  dt &= \int_{\Phi(r)}^{2\Phi(\diam(D))} \frac{h(s,x,y)}{V(x,\Phi^{-1}(s))} \frac{\phi'(1/s)}{s^2 \phi(1/s)^2}  ds\nn\\
&\simeq \int_{\Phi(r)}^{2\Phi(\diam(D))}  \frac{h(s,x,y)}{sV(x,\Phi^{-1}(s))\phi(1/s)}ds.
\end{align}
Combining  with Theorem \ref{t:smallgeneral} (with $T=T_D$), we arrive at  the lower bound in \eqref{e:Green_integral_form}.

 By Theorems \ref{t:smallgeneral} and \ref{t:Slarge} (with $T=T_D$), 
 we have that
\begin{align*}
G_D(x,y)&
\le 
c_0\int_0^{\psi(r)} \frac{h(\phi^{-1}(1/t)^{-1},x,y)}{V(x,\psi^{-1}(t))} \exp \Big(-\frac{c_1r^2}{\psi^{-1}(t)^2} \Big) dt  +c_0h(\Phi(r),x,y) \frac{w(\Phi(r))}{V(x,r)} \int_0^{\psi(r)}tdt \\
&\quad +c_0 \frac{C_0}{V(x,r)\Psi(r)} \int_0^{\psi(r)} t \sB_h(t,x,y)  dt  + c_0\frac{C_0}{V(x,r)\Psi(r)} \int_0^{\psi(r)}   \frac{h(\phi^{-1}(1/t)^{-1},x,y)}{\phi^{-1}(1/t)} dt  \\
&\quad +c_0 \int_{\psi(r)}^{T_D} \frac{h(\phi^{-1}(1/t)^{-1},x,y)}{V(x,\psi^{-1}(t))}  dt +c_0h(1,x,y) \int_{T_D}^\infty e^{-t\phi(\lambda_D)}dt\\
&=: c_0(G_1+G_2+C_0G_3+C_0G_4+G_5+G_6).
\end{align*} 
First note that, following the proof of  Lemma \ref{l:leftint-Phi} and with help of \eqref{e:upper-scaling-psi}, one can see that under the assumptions of Lemma \ref{l:leftint-Phi}, there exists 
$c_1>0$ such that for all $r, \kappa \in (0, T_D)$,
$$
	\int_0^r f(s) \exp\Big( - \frac{\kappa^2}{ \psi^{-1}(s)^2}\Big) ds \le 
	\frac{c_1r^{p+1}f(r)}{\psi(\kappa)^p}.
$$
Applying this inequality 
with $f(t)=h(\phi^{-1}(1/t)^{-1},x,y) V(x,\psi^{-1}(t))^{-1}$ and $p:=\gamma/\beta_1+d_2/(\alpha_1\beta_1)$, we get from (H2), 
\eqref{e:phi-inv-scaling}, \eqref{e:volume_doubling} and
\eqref{e:upper-scaling-psi} that
\begin{align*}
G_1 &\le 
c_2 \frac{h(\Phi(r),x,y)}{V(x,r)} \frac{\psi(r)^{p+1}}{\psi(r)^p}=c_2
h(\Phi(r),x,y)\frac{\psi(r)}{V(x,r)}.
\end{align*}

For $G_2$,  we see from \eqref{e:H-w} that 
\begin{align*}
G_2 \le \frac{e}{2(e-2)}h(\Phi(r),x,y)\frac{\psi(r)^2}{V(x,r) \psi(r)}=\frac{e}{2(e-2)}h(\Phi(r),x,y)\frac{\psi(r)}{V(x,r)}.
\end{align*}

For $G_3$, we use Fubini's theorem to get that
\begin{align}\label{e:G1decom}
 V(x,r)\Psi(r) G_3  &= \int_0^{\psi(r)} t \int_{2\phi^{-1}(1/t)^{-1}}^{4\Phi(r)} h(s,x,y)w(s) ds dt \nn\\
&= \int_0^{2\Phi(r)} h(s,x,y)w(s)\int_0^{\phi(2/s)^{-1}} t dt  ds +  \int_{2\Phi(r)}^{4\Phi(r)} h(s,x,y)w(s)\int_0^{\psi(r)} t dt  ds \nn\\
&=:G_{3,1} + G_{3,2}.
\end{align}
By \eqref{e:H-w}, a change of variables and (H2), we get
\begin{align*}
	G_{3,1} &\le  
c_3 \int_0^{2\Phi(r)} \frac{h(s,x,y)\phi(2/s)}{\phi(2/s)^2} ds  =   2c_3 \int_0^{\Phi(r)} \frac{h(2s,x,y)}{\phi(1/s)} ds \le c_4 \int_0^{\Phi(r)} \frac{h(s,x,y)}{\phi(1/s)} ds.
\end{align*}
On the other hand,  we also get from   \eqref{e:H-w}, (H1) and  \eqref{e:phi-upper-scaling} that
\begin{align*}
G_{3,2} &\le \frac{e}{2(e-2)}h(\Phi(r),x,y) \phi(1/\Phi(r)) \psi(r)^2 \int_{2\Phi(r)}^{4\Phi(r)}ds =\frac{e}{(e-2)}h(\Phi(r),x,y)  \psi(r) \Phi(r) \\
& \le \frac{2e}{(e-2)} \frac{\psi(r)}{\phi(2/\Phi(r))^{-1}} \int_{\Phi(r)/2}^{\Phi(r)} \frac{h(s,x,y)}{\phi(1/s)} ds \le 
c_5 \int_{\Phi(r)/2}^{\Phi(r)} \frac{h(s,x,y)}{\phi(1/s)} ds \le c_5 \int_{0}^{\Phi(r)} \frac{h(s,x,y)}{\phi(1/s)} ds.
 \end{align*}

Clearly, $G_6\le \phi(\lambda_D)^{-1}h(1,x,y)$.

Recall from \eqref{e:green_G4} and \eqref{e:green_G5} that
$$G_4+G_5 \simeq \frac{C_0}{V(x, r)\Psi(r)} \int_0^{\Phi(r)} \frac{h(s,x,y)}{\phi(1/s)}ds+  \int_{\Phi(r)}^{2\Phi(\diam(D))}  \frac{h(s,x,y)}{sV(x,\Phi^{-1}(s))\phi(1/s)}ds.$$

It follows from Lemma \ref{l:newrs} and the upper bounds above  on $G_1, G_2, G_6$ 
that $G_5$ dominates $G_1+G_2+G_6$. 
Since $G_4$ dominates $G_3$, the proof is complete.
\qed

In the remainder of this section, under some additional assumptions on the boundary function, we will obtain 
Green function estimates in 
simpler forms.  
Lemma \ref{l:intscale} will be a useful tool in all simplifications.

We start with the following condition which is a counterpart of 
(H2).

\medskip

\noindent (H2*) There exist constants $c_1,\gamma_* >0$ such that for all 
$x,y \in D$ and $s,t\ge 0$ with 
$\Phi(\delta_\vee(x, y))
\le s \le t < 2\Phi(\diam(D))$,
\begin{align*}
s^{\gamma_*} h(s,x,y) \ge c_1 t^{\gamma_*} h(t, x, y).
\end{align*}

Note that the $\gamma_*$ above is less than or equal to $\gamma$.

\begin{remark}\label{r:H2*_epsilon}
	{\rm Suppose that 
		the boundary function $h(t,x,y)$ satisfies (H2*). Then for every $\epsilon \in (0,1)$, there exists $c_2=c_2(\epsilon)>0$ 
such that for all 
$x,y \in D$ and $s,t\ge 0$ with $\epsilon
\Phi(\delta_\vee(x, y))
\le s \le t < 2\Phi(\diam(D))$,
		\begin{align*}
s^{\gamma_*} h(s,x,y) \ge c_2 t^{\gamma_*} h(t, x, y).
		\end{align*}  
Indeed, let $\epsilon
\Phi(\delta_\vee(x, y))
\le s \le 
\Phi(\delta_\vee(x, y))$
and $s\le t < 2\Phi(\diam(D))$.  
If $t \le 
\Phi(\delta_\vee(x, y))$,
then $\epsilon t \le s$ so that by (H1), 
	\begin{align*}
	s^{\gamma_*}h(s,x,y) \ge  	s^{\gamma_*}h(t,x,y)  \ge  \eps^{\gamma_*} t^{\gamma_*}h(t,x,y).
\end{align*}
If $t > 
\Phi(\delta_\vee(x, y))$,
then  by using (H1) in the first inequality below, (H2*) in the second, and the condition that $s \ge \eps
\Phi(\delta_\vee(x, y))$
 in the last inequality, we see that
		\begin{align*}
		s^{\gamma_*}h(s,x,y) \ge s^{\gamma_*}h( 
		\Phi(\delta_\vee(x, y)), x, y)
		\ge c_1\Big (\frac{s}{
		\Phi(\delta_\vee(x, y))
		}\Big)^{\gamma_*}t^{\gamma_*}h(t, x,y) \ge c_1 \epsilon^{\gamma_*} t^{\gamma_*}h(t,x,y).
		\end{align*}
	}
\end{remark}

\begin{example}\label{ex:H2_opposite}
	{\rm Let $p,q\ge0$, $p+q>0$. Recall that the boundary function $h_{p,q}(t,x,y)$ defined 
   in \eqref{e:def-hp} satisfies (H2) with $\gamma=p+q$. 
    We claim that
	$h_{p,q}(t,x,y)$ also satisfies (H2*) with $\gamma_*=\gamma=p+q$. Indeed, 
	 for all $x,y \in D$ and $ 
	 \Phi(\delta_\vee(x, y))
	 <s<t$, 
	\begin{align*}
		s^{p+q}h_{p,q}(s,x,y)&=\Phi(\delta_D(x))^p \Phi(\delta_D(y))^q  = t^{p+q}
h_{p, q}(t,x,y).	
	\end{align*}	
}
\end{example}

In the remainder of this section,  
 we let $d_1,d_2, \gamma,\gamma_*,\beta_1,\beta_2$ and $\alpha_1, \alpha_2$ be the constants in  \eqref{e:volume_doubling}, (H2), (H2*),  {\bf (Poly-$R_1$)} and the scaling indices of  $\Phi$ in \eqref{e:psiphi}, respectively.

Let 
\begin{equation}\label{e:gReenpart}
\widetilde{G}_D(x, y):=\int_{\Phi(\rho(x,y))}^{2\Phi(\diam(D))}  \frac{h(s,x,y)}{sV(x,\Phi^{-1}(s))\phi(1/s)}ds
\end{equation}
denote the second term on the right-hand side of the estimate \eqref{e:Green_integral_form}.

\begin{lemma}\label{l:newrs2}
The following estimates hold for all $x,y \in D$.	
	
	\noindent (i) If  $d_1 > \alpha_2(\beta_2 \wedge 1)$, then 
	\begin{equation*}
		\widetilde{G}_D(x, y)
		\simeq  h\big(\Phi(\rho(x,y)),x,y\big) \frac{\psi(\rho(x,y))}{V(x,\rho(x,y)) }.
	\end{equation*}

\noindent (ii) If 
$d_2<\alpha_1(\beta_1 - \gamma)$, then
\begin{equation*}
\widetilde{G}_D(x, y)
	\simeq  \begin{cases}
		h(1,x,y), &\mbox{ when \ \HKBh \ holds,}\\
		\infty , &\mbox{ when \ \HKUh \ holds.}
	\end{cases}
\end{equation*} 

	Below, we also assume  that  $h(t,x,y)$ is  regular  and (H2*) holds.
	
	\smallskip

	\noindent (iii)
If $\alpha_1\beta_1>d_2 \ge d_1 >  \alpha_2((\beta_2 \wedge 1)-\gamma_*)$, then 
\begin{equation*}
	\widetilde{G}_D(x, y)
	\simeq 
	h\big(\Phi(\rho(x,y)),x,y\big)\frac{\psi(\rho(x,y) \vee \delta_\vee(x,y))}{V(x,\rho(x,y) \vee \delta_\vee(x,y)) }. 
\end{equation*}

\noindent (iv) If  
$d_1=d_2=\alpha_1\beta_1=\alpha_2\beta_2$,
then
\begin{equation*}
\widetilde{G}_D(x, y)
	\simeq  
	h\big(\Phi(\rho(x,y)),x,y\big) \log \Big(e+\frac{\delta_\vee(x,y)}{\rho(x,y)}\Big).
\end{equation*}  

\noindent (v) If  $\alpha_1=\alpha_2$, $\beta_1=\beta_2$, $\gamma=\gamma_*$ and  
$d_1=d_2=\alpha_1(\beta_1-\gamma)$, then
\begin{equation*}
\widetilde{G}_D(x, y)
	\simeq  \begin{cases}
h(1,x,y)	\log \left(e+{\diam(D)}(\rho(x,y) \vee \delta_\vee(x,y))^{-1}\right), &\mbox{ when \ \HKBh \ holds,}\\
		\infty , &\mbox{ when \ \HKUh \ holds.}
	\end{cases} 
\end{equation*}  
\end{lemma}

\pf Take $x,y \in D$. 
Let $\delta_{\wedge}:=\delta_{\wedge}(x,y)$ and $\delta_{\vee}:=\delta_{\vee}(x,y)$.  
Define 
	\begin{equation*}
	g(s):= \frac{h(s,x,y)}{V(x,\Phi^{-1}(s))\phi(1/s)}, \quad s>0.
\end{equation*} 
Then
$$
\widetilde{G}_D(x,y) = \int_{\Phi(\rho(x,y))}^{2\Phi(\diam(D))}  \frac{g(s)}{s}ds.
$$

By (H1), (H2),  \eqref{e:volume_doubling}, \eqref{e:psiphi}, \eqref{e:phi-upper-scaling} and \eqref{e:phi-lower-scaling},  there exist $c_1,c_2>0$ such that
\begin{equation}\label{e:green_integrand_0}
	c_1\Big(\frac{r}{s}\Big)^{-\gamma-d_2/\alpha_1+\beta_1} \le \frac{g(r)}{g(s)}\le c_2\Big(\frac{r}{s}\Big)^{-d_1/\alpha_2+(\beta_2\wedge 1)}, \quad 0<s \le r<2\Phi(\diam(D)).
\end{equation} 
If $h(t,x,y)$ is regular, then by Remark \ref{r:regular_boundary}, for every $a>0$, there exists $c_3=c_3(a)>0$ such that
\begin{equation}\label{e:green_integrand_1}
	c_3\Big(\frac{r}{s}\Big)^{-d_2/\alpha_1+\beta_1} \le  \frac{g(r)}{g(s)} \le c_2\Big(\frac{r}{s}\Big)^{-d_1/\alpha_2+(\beta_2\wedge 1)}, \;\; 
	0<s \le r<\Phi(a\delta_{\wedge}
	) \wedge 2\Phi(\diam(D));
\end{equation}
if furthermore (H2*) further holds, then by Remark \ref{r:H2*_epsilon},  there exists $c_4>0$ such that 
\begin{equation}\label{e:green_integrand_2}
	c_1\Big(\frac{r}{s}\Big)^{-\gamma-d_2/\alpha_1+\beta_1} \le\frac{g(r)}{g(s)}\le c_4\Big(\frac{r}{s}\Big)^{-\gamma_*-d_1/\alpha_2+(\beta_2\wedge 1)}, \;\; \Phi(\delta_{\vee}/2)  <s \le r<2\Phi(\diam(D)).
\end{equation}

(i) By  \eqref{e:green_integrand_0}, since $-d_1/\alpha_2 + (\beta_2 \wedge 1)<0$, the result follows from Lemma \ref{l:intscale}(ii).

(ii) If $D$ is bounded, then by  \eqref{e:green_integrand_0} and Lemma \ref{l:intscale}(i), since $-\gamma-d_2/\alpha_1 + \beta_1>0$, it holds that $
\widetilde{G}_D(x,y)
\simeq g(\Phi(\diam(D)))$. By \eqref{e:compare_ball}, there exists a constant 
$c_5>1$ such that $c_5^{-1} \le V(z, \diam(D)) \le c_5$ for all $z \in D$.
Hence, 
by using (H1), (H2) and the definition of $g$,
we get that $G_D(x,y) \simeq h(1,x,y)$. If $D$ is unbounded, then we see from \eqref{e:green_integrand_0} and Lemma \ref{l:intscale}(i) that
\begin{equation*}
\widetilde{G}_D(x,y)
\simeq \lim_{r \to \infty} \int_{\Phi(\rho(x,y))}^r \frac{g(s)}{s}ds \simeq \lim_{r \to \infty} g(r) \ge c_1g(1) \lim_{r \to \infty} r^{-\gamma-d_2/\alpha_1+\beta_1}=\infty.
\end{equation*}

(iii)  Suppose that $\delta_{\vee} \le 2\rho(x,y)$. 
Since
$-\gamma_*-d_1/\alpha_2+(\beta_2 \wedge 1)<0$,  by \eqref{e:green_integrand_2} and Lemma \ref{l:intscale}(ii), 
\begin{align*}
\widetilde{G}_D(x,y)
\simeq g(\Phi(\rho(x,y))) = \frac{h(\Phi(\rho(x,y)),x,y)\psi(\rho(x,y))}{V(x, \rho(x,y))}.
\end{align*}
Hence the  result follows from \eqref{e:volume_doubling} and \eqref{e:upper-scaling-psi}.

Suppose now that $\delta_{\vee}>2\rho(x,y)$. 
Then    $\delta_{\wedge} \ge \delta_{\vee}-\rho(x,y)> \delta_{\vee}/2>\rho(x,y)$. Since $h$ is regular,  we get $h(\Phi(\delta_{\vee}),x,y)  \simeq  h(\Phi(\rho(x,y)),x,y) \simeq 1$. Further,
 since $-d_2/\alpha_1+\beta_1>0$ and $-\gamma_*-d_1/\alpha_2+(\beta_2 \wedge 1)<0$,  by the scaling of $\Phi$, \eqref{e:green_integrand_1}, \eqref{e:green_integrand_2}  and Lemma \ref{l:intscale}(i)-(ii), we get 
 \begin{align*}
\widetilde{G}_D(x,y)
 	&\simeq \int_{\Phi(\rho(x,y))}^{2\Phi(\delta_{\wedge})} \frac{g(s)}{s}ds + \int_{\Phi(\delta_{\vee})}^{2\Phi(\diam(D))} \frac{g(s)}{s}ds \simeq g(\Phi(\delta_{\wedge})) + g(\Phi(\delta_{\vee})) \\
 	& \simeq  g(\Phi(\delta_{\vee}))= \frac{h(\Phi(\delta_{\vee}),x,y)\psi(\delta_{\vee})}{V(x, \delta_{\vee})}  \simeq \frac{\psi(\delta_{\vee})}{V(x, \delta_{\vee})}  \simeq \frac{h(\Phi(\rho(x,y)),x,y)\psi(\delta_{\vee})}{V(x, \delta_{\vee})}.
 \end{align*}
This finishes the proof for (iii).

(iv) Since $d_1=d_2$, by \eqref{e:volume_doubling} and \eqref{e:compare_ball}, we see that for every $a>0$, there are comparability constants depending on $a$  such that for all $w,z \in D$ and $0<r<a\,\diam(D)$,
\begin{equation}\label{e:uniform_volume}
V(w,r) \simeq\Big (\frac{r}{\rho(w,z)}\Big)^{d_1}V(w, \rho(w,z))  \simeq \Big(\frac{r}{\rho(w,z)}\Big)^{d_1}V(z, \rho(w,z)) \simeq V(z,r) \simeq r^{d_1}V(z,1).
\end{equation}
Moreover, since $\beta_1=\beta_2$ and $\alpha_1=\alpha_2$, by  \eqref{e:phi-upper-scaling}, \eqref{e:phi-lower-scaling} and \eqref{e:psiphi}, we get that
\begin{equation}\label{e:phi_strong_stable}
\phi(1/s)^{-1} \simeq s^{\beta_1},  \quad  0<s<2\Phi(\diam(D)) \quad \text{ and } \;\quad 	\Phi^{-1}(s)\simeq s^{1/\alpha_1},  \quad  s>0,
\end{equation}
so that $g(s) \simeq h(s,x,y)$ for all   $0<s<2\Phi(\diam(D))$. In particular, since $h$ is regular, we see from Remark \ref{r:regular_boundary} that  
\begin{equation}\label{e:Green_iv}
	g(s) \simeq 1, \quad 0<s<
	2\Phi(\delta_{\wedge}).
\end{equation}

 If $\delta_{\vee} \le 2\rho(x,y)$, then by \eqref{e:green_integrand_2} and Lemma \ref{l:intscale}(ii),
$$
\widetilde{G}_D(x,y)
\simeq g(\Phi(\rho(x,y))) \simeq h(\Phi(\rho(x,y)),x,y)\simeq h(\Phi(\rho(x,y)),x,y)  \log \Big(e+\frac{\delta_{\vee}}{\rho(x,y)}\Big).
$$

If $\delta_{\vee}>2\rho(x,y)$, then we get $\delta_{\wedge}>\delta_{\vee}/2> \rho(x,y)$ as in (iii), and by \eqref{e:Green_iv}, \eqref{e:green_integrand_2} and Lemma \ref{l:intscale}(ii),
 \begin{align*}
	\widetilde{G}_D(x,y)
	&\simeq \int_{\Phi(\rho(x,y))}^{2\Phi(\delta_{\wedge})} \frac{g(s)}{s}ds + \int_{\Phi(\delta_{\vee})}^{2\Phi(\diam(D))} \frac{g(s)}{s}ds \simeq \int_{\Phi(\rho(x,y))}^{2\Phi(\delta_{\wedge})} \frac{ds}{s} + g(\Phi(\delta_{\vee})).
\end{align*}
Note that since $\Phi(s) \simeq s^{\alpha_1}$ for $s>0$, we have $\int_{\Phi(\rho(x,y))}^{2\Phi(\delta_{\wedge})} s^{-1}ds \simeq \log \Big(e+\frac{\delta_{\wedge}}{\rho(x,y)}\Big)$ 
so that
$$ 
g(\Phi(\delta_{\vee})) \simeq h(\Phi(\delta_{\vee}),x,y) \le 1 \le \log \Big(e+\frac{\delta_{\wedge}}{\rho(x,y)}\Big) \simeq \int_{\Phi(\rho(x,y))}^{2\Phi(\delta_{\wedge})} \frac{ds}{s}.
$$
Eventually, since $\delta_{\wedge} \simeq \delta_{\vee}$ and $h(\Phi(\rho(x,y)),x,y) \simeq 1$ in this case, we obtain that
$$
\widetilde{G}_D(x,y)
\simeq  \log \Big(e+\frac{\delta_{\wedge}}{\rho(x,y)}\Big) \simeq  h(\Phi(\rho(x,y)),x,y)\log \Big(e+\frac{\delta_{\vee}}{\rho(x,y)}\Big).
$$

(v) By \eqref{e:uniform_volume},  \eqref{e:phi_strong_stable}, the regularity of $h$,  (H2), Remark \ref{r:H2*_epsilon} and \eqref{e:boundary-function-c}, we have
\begin{equation}\label{e:gamma=gamma0_0}
g(s) \simeq s^{\gamma}, \quad 0<s<\Phi(\delta_{\wedge})
\end{equation}
and
\begin{equation}\label{e:gamma=gamma0}
g(s) \simeq s^{\gamma}h(s,x,y) \simeq t^\gamma h(t,x,y), \quad \;\; \Phi(\delta_{\vee}/2)<s \le t<2\Phi(\diam(D))+1.
\end{equation}

If $\delta_{\vee} \le 2\rho(x,y)$, then since $\Phi(s) \simeq s^{\alpha_1}$ for $s>0$ in this case, we get from \eqref{e:gamma=gamma0} that 
\begin{align*}
\widetilde{G}_D(x,y)	
	\simeq  h(1,x,y)\int_{\Phi(\rho(x,y))}^{2\Phi(\diam(D))} \frac{ds}{s} \simeq h(1,x,y) 	\log \Big(e+\frac{\diam(D)}{\rho(x,y)}\Big).
\end{align*}

If $\delta_{\vee} >2\rho(x,y)$, then $\delta_{\wedge}>\delta_{\vee}/2> \rho(x,y)$ as in (iii) and hence by \eqref{e:gamma=gamma0_0} and \eqref{e:gamma=gamma0},
 \begin{align*}
\widetilde{G}_D(x,y)	
	&\simeq \int_{\Phi(\rho(x,y))}^{2\Phi(\delta_{\wedge})} s^{\gamma-1}ds +h(1,x,y) \int_{\Phi(\delta_{\vee})}^{2\Phi(\diam(D))} \frac{ds}{s}. 
\end{align*}
Since $\Phi(s) \simeq s^{\alpha_1}$ for $s>0$, $\delta_{\wedge} \le \delta_{\vee} \le 2\delta_{\wedge}$ and $h$ is regular, by \eqref{e:gamma=gamma0},
\begin{align*}
&h(1,x,y) \int_{\Phi(\delta_{\vee})}^{2\Phi(\diam(D))} \frac{ds}{s} \simeq h(1,x,y) 	\log \Big(e+\frac{\diam(D)}{\delta_{\vee}}\Big)\\
& \ge h(1,x,y) \simeq \Phi(\delta_{\wedge})^{\gamma} h(\Phi(\delta_{\wedge}),x,y) \simeq \Phi(\delta_{\wedge})^{\gamma} \ge  \gamma^{-1} \int_{\Phi(\rho(x,y))}^{2\Phi(\delta_{\wedge})} s^{\gamma-1}ds.
\end{align*}
This completes the proof. \qed

In the next lemma we show that under the additional assumption that $\gamma<\beta_1+1$, the first term on the right-hand side of \eqref{e:Green_integral_form} is dominated by $\wt{G}_D(x,y)$.

\begin{lemma}\label{l:Sgreen}
If either $C_0=0$ or $\gamma<\beta_1+1$, then
$
G_D(x, y)\simeq 
\widetilde{G}_D(x, y)  \text{ on }D\times D.
$
\end{lemma}

\pf When $C_0=0$, the assertion follows from Proposition \ref{p:Sgreen}. So we now assume $\gamma<\beta_1+1$. According to \eqref{e:phi-lower-scaling} and (H2), we have
\begin{equation*}\frac{h(t,x,y)\phi(1/t)^{-1}}{h(s,x,y)\phi(1/s)^{-1}}  \ge c\Big(\frac{t}{s}\Big)^{\beta_1-\gamma} \quad \text{for all} \;\;0<s \le t < \Phi(\diam(D)).
\end{equation*}
Thus, since $\Psi \ge \Phi$, by Lemma \ref{l:intscale}(i), (H1) and  \eqref{e:phi-upper-scaling}, we get
\begin{equation*}
\frac{C_0}{V(x,r)\Psi(r)}\int_0^{\Phi(r)}\frac{h(s,x,y)}{\phi(1/s)}ds \le  c\frac{h(\Phi(r),x,y) \psi(r)\Phi(r)}{V(x,r)\Psi(r)} \le c\frac{h(\Phi(r),x,y) \psi(r)}{V(x,r)}.
\end{equation*}
Combining the above with 
Lemma \ref{l:newrs} and Proposition \ref{p:Sgreen}, we get the assertion. \qed

Define
\begin{equation}\label{e:newrs}
	g_0(x,y)=
	\begin{cases} \displaystyle
	\frac{\psi(\rho(x,y))}{V(x,\rho(x,y))}
		,& \mbox{ if } d_1>\alpha_2(\beta_2 \wedge 1),\\[7pt]
 \displaystyle	\log\Big( e +  \frac{\delta_\vee(x,y)}{\rho(x,y)}\Big), &\mbox{ if } d_1=d_2=\alpha_1\beta_1=\alpha_2\beta_2, \\[7pt]
	 \displaystyle	\frac{\psi(\rho(x,y) \vee \delta_\vee(x,y))}{V(x,\rho(x,y) \vee \delta_\vee(x,y))}
		, & \mbox{ if } d_2<\alpha_1\beta_1.	
	\end{cases}
	\end{equation}

By combining Proposition \ref{p:Sgreen}, Lemma \ref{l:newrs2} and Lemma \ref{l:Sgreen} we arrive at the following result.
\begin{thm}\label{t:Sgreen}
Suppose that $C_0=0$ or $\gamma<\beta_1+1$, $h(t,x,y)$ is regular and (H2*) holds. 
\begin{itemize}
	\item[(a)] Suppose also that one of the following holds: 
 (1) $ d_1>\alpha_2(\beta_2 \wedge 1)$ or (2) $d_1=d_2=\alpha_1\beta_1=\alpha_2\beta_2$ or (3) $d_2<\alpha_1\beta_1$. Then it holds that
\begin{equation}\label{e:Sgreen}
G_D(x,y)\simeq   h(\Phi(\rho(x,y)),x,y) g_0(x,y).
\end{equation}
	\item[(b)] If $d_2<\alpha_1(\beta_1 - \gamma)$, then
\begin{equation*}
	G_D(x,y) 
	\simeq  \begin{cases}
		h(1,x,y), &\mbox{ when \ \HKBh \ holds,}\\
		\infty , &\mbox{ when \ \HKUh \ holds.}
	\end{cases}
\end{equation*} 
	\item[(c)] If  $\alpha_1=\alpha_2$, $\beta_1=\beta_2$, $\gamma=\gamma_*$ and  
$d_1=d_2=\alpha_1(\beta_1-\gamma)$, then
\begin{equation*}
	G_D(x,y) 
	\simeq  \begin{cases}
h(1,x,y)	\log \Big(e+{\diam(D)}(\rho(x,y) \vee \delta_\vee(x,y))^{-1}\Big), &\mbox{ when \ \HKBh \ holds,}\\
		\infty , &\mbox{ when \ \HKUh \ holds.}
	\end{cases} 
\end{equation*}  
\end{itemize}
\end{thm}

\medskip

When $C_0=1$, Theorem \ref{t:Sgreen} only deals with the case $\gamma<\beta_1+1$.
To cover the case when $\gamma$ is large, we assume the following condition.

\medskip

\noindent (H2**) There exist constants  $c_1>0$, $\gamma_{**} \in (0, 1+\beta_1)$ such that for all 
$x,y \in D$ and $s,t\ge 0$ with 
$\Phi(\delta_\wedge(x, y))\le s \le t < \Phi(\delta_\vee(x, y))$,
\begin{align*}
	s^{\gamma_{**}} h(s,x,y) \le c_1 t^{\gamma_{**}} h(t, x, y).
\end{align*}

\begin{example}\label{ex:H2**}
	{\rm For $p,q\ge0$, let $h_{p,q}(t,x,y)$ be the boundary function defined in \eqref{e:def-hp}. If $p \vee q <1+\beta_1$, then  $h_{p,q}(t,x,y)$ satisfies (H2**). Indeed, we see that  for all $x,y \in D$ and 
$\Phi(\delta_\wedge(x, y))\le s \le t < \Phi(\delta_\vee(x, y))$,
		\begin{equation*}
			s^{p \vee q}h_{p,q}(s,x,y) = \begin{cases}
			\Phi(\delta_D(x))^p s^{p\vee q - p}, &\mbox{ if } \delta_D(x) < \delta_D(y)\\
				\Phi(\delta_D(y))^q s^{p \vee q - q}, &\mbox{ if } \delta_D(x) > \delta_D(y) 
			\end{cases}\,
		\le t^{p \vee q} h_{p,q}(t,x,y).
		\end{equation*}
		}
\end{example}

\medskip

For a given boundary function $h$, we define for $x,y \in D$,
\begin{align*}	
	[h](x,y)&:=h\big(\Phi( \rho(x,y) \wedge \delta_\vee(x,y)),x,y\big)\left( 1 \wedge  \frac{\Phi(\delta_\vee(x,y))\psi( \delta_\vee(x,y))}{\Phi(\rho(x,y))\psi(\rho(x,y))}\right).
\end{align*}
One can see that for all $p,q \ge 0$,
\begin{align}\label{e:bracket_hpq}
[h_{p,q}](x,y) &\simeq \Big(1 \wedge \frac{\Phi(\delta_D(x))}{\Phi(\rho(x,y))}\Big)^p\Big(1 \wedge \frac{\Phi(\delta_D(y))}{\Phi(\rho(x,y))}\Big)^q\Big(1 \wedge \frac{\Phi(\delta_\vee(x,y)) }{\Phi(\rho(x,y))}\Big)^{1-p-q}\Big(1 \wedge \frac{ \psi(\delta_\vee(x,y))}{ \psi(\rho(x,y))}\Big).
\end{align}
Indeed, for $x,y \in D$, if $\delta_\wedge(x,y) \ge \rho(x,y)$, then $[h_{p,q}](x,y) = 1$ and the right-hand side of \eqref{e:bracket_hpq} is also equal to $1$. If $\delta_\wedge(x,y)< \rho(x,y)$, then $\delta_\vee(x,y) < \rho(x,y) + \delta_\wedge(x,y) \le \rho(x,y)+ \rho(x,y) \le 2\rho(x,y)$ and hence by \eqref{e:boundary-function-c}  and the scaling properties of $\Phi$ and $\psi$, 
\begin{align*}
\frac{\Phi(\rho(x,y))^{1-p-q} \psi(\rho(x,y))}{\Phi(\delta_\vee(x,y))^{1-p-q} \psi(\delta_\vee(x,y))}	[h_{p,q}](x,y)&\simeq   \frac{\Phi(\delta_\vee(x,y))^{p+q}}{\Phi(\rho(x,y))^{p+q} } h_{p,q}(\Phi(\delta_\vee(x,y),x,y)= \frac{\Phi(\delta_D(x))^p\Phi(\delta_D(y))^q}{\Phi(\rho(x,y))^{p+q}}\\
& \simeq \Big(1 \wedge \frac{\Phi(\delta_D(x))}{\Phi(\rho(x,y))}\Big)^p\Big(1 \wedge \frac{\Phi(\delta_D(y))}{\Phi(\rho(x,y))}\Big)^q.
\end{align*}

Recall that $g_0(x, y)$ is defined in \eqref{e:newrs}.

\begin{thm}\label{t:Sgreen_2}
 Suppose that $C_0=1$,  $h(t,x,y)$ is  regular, (H2*) holds with $\gamma_*>(\beta_2 \wedge 1) +1$ and (H2**) holds.   Suppose also that 
one of the following holds: 
 (1) $ d_1>\alpha_2(\beta_2 \wedge 1)$ or (2) $d_1=d_2=\alpha_1\beta_1=\alpha_2\beta_2$ or (3) $d_2<\alpha_1\beta_1$. Then it holds that
\begin{equation}\label{e:Sgreen_2}
G_D(x,y) \simeq  [h](x,y)  \frac{\Phi(\rho(x,y))}{\Psi(\rho(x,y))}\frac{\psi(\rho(x,y))}{V(x,\rho(x,y))}+ h(\Phi(\rho(x,y)),x,y) g_0(x,y).
\end{equation}
In particular, if $\Psi \simeq \Phi$, then
\begin{equation}\label{e:Sgreen_20}
	G_D(x,y) \simeq  [h](x,y)  g_0(x,y).
\end{equation}
	
\end{thm}
\pf 
Take $x,y \in D$ and let $r:=\rho(x, y)$ and
 $\delta_{\vee}:=\delta_{\vee}(x,y)$.
Observe that by \eqref{e:phi-upper-scaling}, \eqref{e:phi-lower-scaling}, (H1), (H2**) and the regularity of $h$, 
\begin{equation}\label{e:Green_0_1}
c_1\Big(\frac{t}{s}\Big)^{\beta_1-\gamma_{**}} \le \frac{h(t,x,y)/\phi(1/t)}{h(s,x,y)/\phi(1/s)} \le c_2\Big(\frac{t}{s}\Big)^{\beta_2 \wedge 1}, \quad 0<s \le t < \Phi(\delta_{\vee}).
\end{equation}
Note also 
that by  \eqref{e:phi-upper-scaling}, \eqref{e:phi-lower-scaling}, (H2), (H2*) and Remark \ref{r:H2*_epsilon},
\begin{equation}\label{e:Green_0_2}
	c_3\Big(\frac{t}{s}\Big)^{\beta_1-\gamma} \le \frac{h(t,x,y)/\phi(1/t)}{h(s,x,y)/\phi(1/s)} \le c_4\Big(\frac{t}{s}\Big)^{\beta_2 \wedge 1-\gamma_*}, \quad \Phi(\delta_{\vee})/2 \le s \le t  < \Phi(\diam(D)).
\end{equation}

If $\delta_{\vee}>2r$, then by Lemma \ref{l:intscale}(i) and  \eqref{e:Green_0_1},
$
\int_0^{\Phi(r)} h(s,x,y)\phi(1/s)^{-1}ds \simeq h(\Phi(r),x,y) \Phi(r) \psi(r).$

If $\delta_{\vee} \le 2r$, then by Lemma \ref{l:intscale}(i)-(ii), \eqref{e:Green_0_1},  \eqref{e:Green_0_2} and the scaling property of $\phi$, since $\beta_1-\gamma_{**}>-1$ and   $\beta_2 \wedge 1 -\gamma_*<-1$, we get that
\begin{align*}
\int_0^{\Phi(r)} \frac{h(s,x,y)}{\phi(1/s)}ds = \int_0^{\Phi(\delta_{\vee})/2} \frac{h(s,x,y)}{\phi(1/s)}ds + \int_{\Phi(\delta_{\vee})/2}^{\Phi(r)} \frac{h(s,x,y)}{\phi(1/s)}ds \simeq h(\Phi(\delta_{\vee}),x,y)\Phi(\delta_{\vee}) \psi(\delta_{\vee}).
\end{align*}

Therefore, in either case, it holds that
\begin{equation}\label{e:Sgreen_0}
	\int_0^{\Phi(r)} \frac{h(s,x,y)}{\phi(1/s)}ds\simeq h(\Phi(r \wedge \delta_{\vee}),x,y)   \left( 1 \wedge  \frac{\Phi(\delta_{\vee})\psi( \delta_{\vee})}{\Phi(r)\psi(r)}\right) \Phi(r)\psi(r) = [h](x,y) \Phi(r)\psi(r).
\end{equation}
Combining this with Proposition \ref{p:Sgreen} and 
Lemma \ref{l:newrs2},
we get \eqref{e:Sgreen_2}.

Now we also assume that $\Psi \simeq \Phi$. If $\delta_{\vee}>r$, then $[h](x,y) = h(\Phi(r),x,y)$. Hence, we see from  
Lemma \ref{l:newrs} 
that in  \eqref{e:Sgreen_2}, the second term dominates the first one so that \eqref{e:Sgreen_20} holds.  If $\delta_{\vee}\le r$, then using  Lemma \ref{l:intscale}(ii), \eqref{e:Green_0_2} and the condition that $\beta_2 \wedge 1 -\gamma_*<-1$ in the second inequality below, 
the scaling property of $\phi$ and  \eqref{e:boundary-function-c}  in the third, and \eqref{e:Sgreen_0} in the fourth, we get
\begin{align*}
&\int_{\Phi(r)}^{2\Phi(\diam(D))}  \frac{h(s,x,y)}{sV(x,\Phi^{-1}(s))\phi(1/s)}ds \le \frac{1}{V(x, r)\Phi(r)}\int_{\Phi(r)}^{2\Phi(\diam(D))}  \frac{h(s,x,y)}{\phi(1/s)}ds\\
& \le c_5 \frac{h(\Phi(r),x,y)\psi(r)}{V(x, r)} \le  \frac{c_6}{V(x, r)\Phi(r)}	\int_{\Phi(r)/2}^{\Phi(r)} \frac{h(s,x,y)}{\phi(1/s)}ds \le c_7[h](x,y) \frac{\psi(r)}{V(x, r)}.
\end{align*}
Note that  $g_0(x,y) \simeq \psi(r)/V(x,r)$ when $\delta_{\vee}\le r$. 
Thus by Proposition \ref{p:Sgreen} and \eqref{e:Sgreen_0}, we get $G_D(x,y) \simeq [h](x,y)\psi(r)/V(x,r) \simeq [h](x,y)g_0(x,y)$ when $\delta_{\vee}\le r$. This completes the proof for \eqref{e:Sgreen_20}.  \qed

For completeness, we record the Green function estimates when $C_0=1$, $\beta_1=\beta_2$ and $\gamma_*=\gamma=\beta_1 +1$. 

\begin{thm}\label{t:Sgreen_3}
	 Suppose that  $C_0=1$, $\beta_1=\beta_2$,  $h(t,x,y)$ is  regular and   (H2*) holds with $\gamma_*=\gamma=\beta_1 +1$ and (H2**) holds.    Suppose also 
that one of the following holds: 
(1) $ d_1>\alpha_2(\beta_2 \wedge 1)$ or (2) $d_1=d_2=\alpha_1\beta_1=\alpha_2\beta_2$ or (3) $d_2<\alpha_1\beta_1$. Then it holds that
	\begin{equation}\label{e:Sgreen_3}
		G_D(x,y) \simeq  h(\Phi(\rho(x,y)),x,y)  \left[ \frac{ \Phi(\rho(x,y))}{\Psi(\rho(x,y)} \frac{ \Phi(\rho(x,y))^{\gamma-1}}{V(x,\rho(x,y))} 
		\log \Big( e + \frac{\Phi(\rho(x,y))}{\Phi(\delta_\vee(x,y))}\Big) +g_0(x,y)\right].
	\end{equation} 
In particular, if $\Psi \simeq \Phi$, then
	\begin{equation}\label{e:Sgreen_30}
	G_D(x,y) \simeq  h(\Phi(\rho(x,y)),x,y) g_0(x,y) \log \Big( e + \frac{\Phi(\rho(x,y))}{\Phi(\delta_\vee(x,y))}\Big).  
\end{equation} 
\end{thm}
\pf 
Take $x,y \in D$ and let $r:=\rho(x, y)$ and $\delta_{\vee}=\delta_{\vee}(x,y)$. 
Using Lemma \ref{l:intscale}(i) (which is applicable due to  \eqref{e:Green_0_1}),
\eqref{e:phi_strong_stable}, the second comparability in \eqref{e:gamma=gamma0},  \eqref{e:boundary-function-c} and the scaling property of $\Phi$,  we get 
that,  if $\delta_{\vee} >2r$, then 
\begin{align*}
	\int_0^{\Phi(r)} \frac{h(s,x,y)}{\phi(1/s)}ds \simeq 	\int_0^{\Phi(r)} s^{\beta_1}h(s,x,y)ds \simeq \Phi(r)^\gamma h(\Phi(r),x,y)
\end{align*}
and if $\delta_{\vee} \le 2r$, then
\begin{align*}
	&	\int_0^{\Phi(r)} \frac{h(s,x,y)}{\phi(1/s)}ds  = 	\int_0^{2^{-1}\Phi(2^{-1}(\delta_{\vee}))} \frac{h(s,x,y)}{\phi(1/s)}ds + \int_{2^{-1}\Phi(2^{-1}(\delta_{\vee}))}^{\Phi(r)} \frac{h(s,x,y)}{\phi(1/s)}ds \\
	& \simeq \int_0^{2^{-1}\Phi(2^{-1}(\delta_{\vee}))} s^{\beta_1} h(s,x,y)ds +\Phi(r)^\gamma h(\Phi(r),x,y) \int_{2^{-1}\Phi(2^{-1}(\delta_{\vee}))}^{\Phi(r)} s^{\beta_1-\gamma }ds \\
	& \simeq  \Phi(\delta_{\vee})^\gamma h(\Phi(\delta_{\vee}),x,y) + \Phi(r)^\gamma h(\Phi(r),x,y) \log \Big( e + \frac{\Phi(r)}{\Phi(\delta_{\vee})}\Big)\\
	& \simeq \Phi(r)^\gamma h(\Phi(r),x,y) \log \Big( e + \frac{\Phi(r)}{\Phi(\delta_{\vee})}\Big).
\end{align*}
The last comparability above is valid by the second comparability in \eqref{e:gamma=gamma0}. 
Therefore, 
in either case,  it holds that
$$
\int_0^{\Phi(r)} \frac{h(s,x,y)}{\phi(1/s)}ds \simeq \Phi(r)^\gamma h(\Phi(r),x,y) \log \Big( e + \frac{\Phi(r)}{\Phi(\delta_{\vee})}\Big).
$$
Combining this with Proposition \ref{p:Sgreen} and Theorem \ref{t:Sgreen}, we get \eqref{e:Sgreen_3}. 

Now we also assume that $\Psi \simeq \Phi$. Then as in the proof of Theorem \ref{t:Sgreen_2}, one can check that in  \eqref{e:Sgreen_3}, 
if $\delta_{\vee} > r$, then  the second term dominates the first one, and if $\delta_{\vee} \le  r$, then  the first term dominates the second one. Combining with the facts that $\psi(r) \simeq \Phi(r)^{\gamma-1}$ for $0<r<\diam(D)$ since $\beta_1=\beta_2$ and that $g_0(x,y) \simeq \psi(r)/V(x,r)$ when $\delta_{\vee}\le r$, we obtain \eqref{e:Sgreen_30}.  \qed

%%%%%%%%%%%%%%%%%%%%%%%%%%%%%%%%%%%%%%%%%%%%%%%%%%%%%%%%%%%%%%%%%%%%%%%%%%%%%%%%%%%
%%%%%%%%%%%%%%%%%%%                      Parabolic Harnack inequality   and H\"older regularly                      %%%%%%%%%%%%%%%%%%%%%%%%
%%%%%%%%%%%%%%%%%%%%%%%%%%%%%%%%%%%%%%%%%%%%%%%%%%%%%%%%%%%%%%%%%%%%%%%%%%%%%%%%%%%

\section{Parabolic Harnack inequality   and H\"older regularity}\label{s:phi}

Throughout this section, we assume that
 $h(t,x,y)$ is a regular boundary function  and  that either (1) {\bf (Poly-$R_1$)}\ and 
\HKBh \ hold, or (2) {\bf (Poly-$\infty$)}  and 	\HKUh \ hold. 

For $x_0 \in D$ and $r>0$, 
let $\tau_{B(x_0,r)}:=\inf\{s>0:\, X_s\notin B(x_0,r)\}$ and $X^{B(x_0, r)}$ be the part process of $X$ in $B(x_0,r)$. Denote by  $q_{B(x_0,r)}(t,x,y)$ the heat kernel of $X^{B(x_0, r)}$. By the strong Markov property,
\begin{equation}\label{e:Dirichlet-heat-kernel}
q_{B(x_0,r)}(t,x,y)=q(t,x,y)-\E_x\big[q(t-\tau_{B(x_0,r)}, X_{\tau_{B(x_0,r)}},y); \,  \tau_{B(x_0,r)}<t\big]\, .
\end{equation}

Recall the definition of
$\psi$ in \eqref{e:defpsi}.

\begin{lemma}\label{l:lower-hke-killed}
There exist constants $C>0$ and $\epsilon\in (0,1/4)$ 
 such that for all $x_0\in D$ and $r \in (0, \delta_D(x_0))$,  
$$
q_{B(x_0,r)}(t,x,y)\ge \frac{C}{V(x_0,\psi^{-1}(t))} \quad \textrm{ for all } t\in (0,\psi(\epsilon r)] \textrm{ and }x,y\in  B\big(x_0, \epsilon \psi^{-1}(t)\big ).
$$
\end{lemma}
\pf Since the proofs are similar, we only give the proof 
in case (1).

Fix $x_0\in D$ and $r \in (0, \delta_D(x_0))$. 
Let $\epsilon \in  (0,1/4) $ 
to be chosen later.  Let $0<t \le \psi(\epsilon r)$ and $x,y\in B(x_0, \epsilon  \psi^{-1}(t))$. Clearly, $x,y \in B(x_0, \epsilon^2 r)$. Further,  
$\delta_\wedge(x, y)
>\delta_D(x_0)-\epsilon^2 r>r-r/4> \psi^{-1}(t)> 2\epsilon \psi^{-1}(t) > \rho(x,y).$ Therefore, we have
$
\delta_\wedge(x, y)
> \psi^{-1}(t)= \rho(x,y) \vee \psi^{-1}(t).
$
From Corollary \ref{c:smallint}(i) (with $T=\psi(\diam(D))$) and  \eqref{e:compare_ball}, it follows that there exist constants
$c_1,c_2>0$ (independent of $t,x_0, x,y$) such that
\begin{equation}\label{e:NDL-1}
q(t,x,y)\ge \frac{c_1}{V(x,\psi^{-1}(t))}\ge \frac{c_2}{V(x_0,\psi^{-1}(t))}\,.
\end{equation}

Let $z\in D\setminus B(x_0,r)$. Then since $t\le \psi(\epsilon r)$ and $\epsilon<1/4$,
\begin{equation*}
	\rho(z,y)\ge  \rho(z,x_0)-\rho(y,x_0)\ge (1-\epsilon^2) \rho(z,x_0) \ge (1-\epsilon^2)r \ge (2\epsilon)^{-1}\psi^{-1}(t).
\end{equation*}
Then according to Proposition \ref{p:hke-upper-rough}(i), \eqref{e:volume_doubling}  and \eqref{e:compare_ball},  since $h \le 1$, there exist constants $c_3,c_4,c_5>0$ (independent of $t,x_0,z,y$)  such that for every $0<s<t$,
\begin{align}\label{e:NDL-2}
q(s,z,y)&\le \frac{c_3s}{V(z, \rho(z,y))\psi(\rho(z,y))}\le \frac{c_3}{V(z, \rho(z,y))}\le \frac{c_4}{V(x_0, \rho(z,y))} \le \frac{c_5\epsilon^{\,d_1}}{V(x_0, \psi^{-1}(t))}.
\end{align}

Combining \eqref{e:Dirichlet-heat-kernel} with \eqref{e:NDL-1} and \eqref{e:NDL-2},  we get that
\begin{equation*}
	q_{B(x_0,r)}(t,x,y) \ge \frac{c_2 -c_5\epsilon^{d_1}}{V(x_0,\psi^{-1}(t))}.
\end{equation*}
Now we finish the proof by choosing $\epsilon=(c_2/(2c_5))^{1/d_1}$  so that $c_2-c_5 \epsilon^{d_1} = c_2/2$. \qed

\begin{remark}\label{r:lower-hke-killed}{\rm
By using \eqref{e:upper-scaling-psi} we may replace $\psi(\epsilon r)$ and $\epsilon \psi^{-1}(r)$ in the statement of Lemma \ref{l:lower-hke-killed} with $\epsilon\psi(r)$ and $\psi^{-1}(\epsilon r)$ respectively, cf.~\cite[p.3758]{CKW20}.
}
\end{remark}

\begin{lemma}\label{l:exit-time}
There exists a constant $C>1$ such that for all $x\in D$ and $r \in (0, \delta_D(x))$,
\begin{equation}\label{e:E-psi}
C^{-1}\psi(r)\le \E^x[\tau_{B(x,r)}]\le C \psi(r).
\end{equation}
\end{lemma}
\pf  Fix $x\in D$ and $r \in (0, \delta_D(x))$.
Let $\epsilon\in  (0,1/4)$ 
be as in the statement of Lemma \ref{l:lower-hke-killed}. Then by Lemma \ref{l:lower-hke-killed}, we have that
$$
q_{B(x,r)}(\psi(\epsilon r),x,y)\ge \frac{c_1}{V(x,\epsilon r)}, \quad y \in B(x, \epsilon^2 r).
$$ 
By \eqref{e:volume_doubling}, this implies that
$$
\P^x(\tau_{B(x,r)}> \psi(\epsilon r))
\ge \int_{B(x,\epsilon^2r)}q_{B(x,r)}(\psi(\epsilon r),x,y)\, dy\ge  \frac{c_1V(x,\epsilon^2r)}{V(x,\epsilon r)} \ge c_2.
$$
Hence, by Markov's inequality and \eqref{e:upper-scaling-psi}, we get that
$$
\E^x [\tau_{B(x,r)}] \ge \psi(\epsilon r)
\P^x(\tau_{B(x,r)}> \psi(\epsilon r)) 
\ge c_2  \psi(\epsilon r) \ge c_3\psi(r).
$$

To obtain the upper bound in \eqref{e:E-psi}, we first assume that \HKBh \ holds.  We claim that there exists a constant $A>1$ such that
\begin{equation}\label{e:upper-E-onestep}
\sup_{z \in B(x,r)}
\P^z(\tau_{B(x,r)} > \psi(Ar)) 
\le \frac{1}{2}.
\end{equation}
Indeed, according to Proposition \ref{p:hke-upper-rough}(i) and Theorem \ref{t:Slarge}, since $h \le 1$, there exists $c_4>1$ such that
\begin{equation}\label{e:bdd_upper_global}
q(t,z,y) \le c_4 \left(V(z,\psi^{-1}(t))^{-1} \1_{\{t \le 1\}} + e^{-\phi(\lambda_D)t} \1_{\{t > 1\}}\right), \quad z,y \in B(x,r).
\end{equation} 
 Further, by \eqref{e:volume_doubling}, there is  $c_5>1$ such that for all $z \in B(x,r)$,
\begin{equation}\label{e:bdd_upper_global_1}
V(z,c_5r) \ge V(x, (c_5-1)r) \ge 2c_4V(x,r).
\end{equation}
Let $A > c_5$ be a constant such that
\begin{equation}\label{e:bdd_upper_global_2}
\exp \big(\phi(\lambda_D)\psi(Ac_5^{-1}\psi^{-1}(1))  \big) \ge 2 m(D).
\end{equation} 
In case when $r\le c_5^{-1}\psi^{-1}(1)$, we see from \eqref{e:bdd_upper_global} and \eqref{e:bdd_upper_global_1} that for all $z \in B(x,r)$,
\begin{equation*}
\P^z(\tau_{B(x,r)} > \psi(Ar))  
\le \P^z(\tau_{B(x,r)} > \psi(c_5r)) \le \int_{B(x,r)}q(\psi(c_5r),z,y) dy  \le \frac{c_4V(x,r)}{V(z, c_5r)} \le \frac 12.
\end{equation*}
On the other hand, if $r>c_5^{-1}\psi^{-1}(1)$, then since $B(x,r) \subset B(x, \delta_D(x)) \subset D$, we get from \eqref{e:bdd_upper_global} and \eqref{e:bdd_upper_global_2} that for all $z \in B(x,r)$,
\begin{align*}
\P^z(\tau_{B(x,r)} > \psi(Ar)) 
& \le \P^z\big(\tau_{B(x,r)} > \psi(Ac_5^{-1}\psi^{-1}(1))\big) \le V(x,r)\exp \big(-\phi(\lambda_D)\psi(Ac_5^{-1}\psi^{-1}(1))  \big)\\
& \le m(D)\exp \big(-\phi(\lambda_D)\psi(Ac_5^{-1}\psi^{-1}(1))  \big) \le \frac{1}{2}.
\end{align*}
Hence, \eqref{e:upper-E-onestep} holds.

Now, by  \eqref{e:upper-E-onestep} and Markov's property, we see that for all $n \ge 2$,
\begin{align*}
&\sup_{z \in B(x,r)}\P^z(\tau_{B(x,r)} > n\psi(Ar)) = \sup_{z \in B(x,r)}\P^z\big(\tau_{B(x,r)} > n\psi(Ar), \, \tau_{B(x,r)} > \psi(Ar)\big) \\
& \le \sup_{z \in B(x,r)}\P^z\big( \P^{X_{\psi(Ar)}}(\tau_{B(x,r)} > (n-1)\psi(Ar)), \, \tau_{B(x,r)} > \psi(Ar)\big) \\
& \le \sup_{z \in B(x,r)} \P^z(\tau_{B(x,r)} > (n-1)\psi(Ar)) \sup_{z \in B(x,r)} \P^z(\tau_{B(x,r)} > \psi(Ar)) \\
&\le \cdots \le \Big( \sup_{z \in B(x,r)} \P^z(\tau_{B(x,r)} > \psi(Ar)) \Big)^n \le 2^{-n}.
\end{align*}
Therefore, we get from \eqref{e:upper-scaling-psi} that
\begin{align*}
\E^x[\tau_{B(x,r)}] &\le \sum_{n=1}^\infty n\psi(Ar) \P^x\big(\tau_{B(x,r)} \in
 ((n-1)\psi(Ar),  n \psi(Ar)]\big)\\
 & \le c_6A^{\alpha_2(\beta_2 \wedge 1)} \psi(r) \sum_{n=1}^\infty n 2^{-(n-1)} =  4c_6A^{\alpha_2(\beta_2 \wedge 1)} \psi(r).
\end{align*}

Similarly, by using Proposition \ref{p:hke-upper-rough}(ii), we can obtain the upper bound in \eqref{e:E-psi} when {\bf (Poly-$\infty$)}  and 	\HKUh \ hold. \qed

Recall that the jump kernel $J(x,y)$ is given in \eqref{e:jumping-J}.

 \begin{lemma}\label{l:Tail-estimate}
 	There exists a constant $C>1$ such that for all $x\in D$ and $r \in (0, \delta_D(x))$,
 	\begin{equation}\label{e:Tail-psi}
 \int_{D \setminus B(x,r)} J(x,y)dy \le \frac{C}{\psi(r)}.
 	\end{equation}
 \end{lemma}
\pf 
 By Theorem \ref{t:jump-estimate} and the fact that $h \le 1$, we have that
\begin{align*}
 \int_{D \setminus B(x,r)} J(x,y)dy& \le c_1 \left(  \int_{D \setminus B(x,r)}\frac{\int_0^{\Phi(\rho(x,y))} w(s)ds }{V(x,\rho(x,y))\Psi(\rho(x,y))}dy  +  \int_{D \setminus B(x,r)}  \frac{ w\big(\Phi(\rho(x,y))\big)}{V(x,\rho(x,y))} dy \right)\\
 &=:c_1(I_1+I_2).
\end{align*} 
By the first inequality in \eqref{e:phi-w}, since $\Psi \ge \Phi$, we see that
\begin{align*}
I_1 \le \int_{D \setminus B(x,r)} \frac{e\Phi(\rho(x,y))\phi(1/\Phi(\rho(x,y)))}{V(x,\rho(x,y)) \Psi(\rho(x,y))}dy \le  \int_{D \setminus B(x,r)} \frac{e\,dy}{V(x,\rho(x,y)) \psi(\rho(x,y))}.
\end{align*}
It follows from \eqref{e:H-w} that
$
I_2 \le \frac{e}{e-2} \int_{D \setminus B(x,r)} V(x,\rho(x,y))^{-1}\psi(\rho(x,y))^{-1}dy.
$
Hence, by \eqref{e:volume_doubling}, \eqref{e:upper-scaling-psi} and the proof of \cite[Lemma 2.1]{CKW16a}, we conclude that $I_1 + I_2 \le c_2/\psi(r)$. \qed

Let $Z := (V_s, X_s)_{s\ge 0}$ be the time-space process corresponding to $X$, where $V_s = V_0 − s$.
The augmented filtration of $Z$ will be denoted by $(\wt{\FF}_s)_{s\ge 0}$.
The law of the time-space process $s\mapsto Z_s$ starting from $(t, x)$ will be denoted by $\P^{(t,x)}$. For
every open subset $B$ of $[0,∞) ×D$, define $\tau_B = \inf\{s > 0:\, Z_s\notin B\}$ and $\sigma_B=\tau_{B^c}$.

Recall that a Borel measurable function $u:[0,\infty)\times D\to \R$ is parabolic (or caloric) on $(a,b)\times B(x_0,r)$ with respect to the process $X$ if for every relatively compact open set $U\subset (a,b)\times B(x_0,r)$ it holds that $u(t,x)=\E^{(t,x)}u(Z_{\tau_U})$ for all $(t,x)\in U$.

 We denote by $dt \otimes m$ the product  of the Lebesgue measure on $[0,\infty)$ and $m$ on $E$. 

\begin{lemma}\label{l:ckw-3-7} 
Let $\epsilon\in (0,1/4)$ be the constant from Lemma \ref{l:lower-hke-killed}. 
For every $\delta\in (0,\epsilon]$, there exists a constant $C_1>0$  such that for all $x\in D$, $r \in (0, \delta_D(x))$, 
 $t\ge \delta \psi(r)$, and any compact set $A\subset [t-\delta\psi(r), t-\delta\psi(r)/2]\times B(x, \psi^{-1}(\epsilon \delta \psi(r)/2))$, 
\begin{equation}\label{e:ckw-3-7-a}
\P^{(t,x)}(\sigma_A < \tau_{[t-\delta \psi(r), t]\times B(x,r)})\ge C_1 \frac{dt \otimes m (A)}{V(x,r) \psi(r)}. 
\end{equation}
\end{lemma}
\pf Let $\tau_r=\tau_{[t-\delta \psi(r), t]\times B(x,r)}$ and $A_s=\{y\in D:\,  (s,y)\in A\}$. For any $t,r>0$ and $x\in D$ such that $B(x,r)\subset D$,
\begin{align}
&\delta\psi(r)\P^{(t,x)}(\sigma_A <\tau_r)\ge \int_0^{\delta\psi(r)} \P^{(t,x)}\left(\int_0^{\tau_r}\1_A(t-s, X_s)ds > 0\right)du
\nn\\
&\ge \int_0^{\delta\psi(r)} \P^{(t,x)}\left(\int_0^{\tau_r}\1_A(t-s, X_s)ds > u\right)du
= \E^{(t,x)}\left[\int_0^{\tau_r} \1_A (t-s, X_s)ds\right]. \label{e:ckw-3-7-b}
\end{align}
For any $t\ge \delta \psi(r)$, 
\begin{align*}
&\E^{(t,x)}\left[\int_0^{\tau_r} \1_A (t-s, X_s)ds\right]\ge \int_{\delta \psi(r)/2}^{\delta\psi(r)} \P^{(t,x)}\left((t-s, X_s^{B(x,r)})\in A\right) ds \\
&= \int_{\delta \psi(r)/2}^{\delta\psi(r)}\P^x (X_s^{B(x,r)}\in A_{t-s})ds= \int_{\delta \psi(r)/2}^{\delta\psi(r)} ds \int_{A_{t-s}}q_{B(x,r)}(s,x,y)\, dy.
\end{align*}
Let $s\in [\delta\psi(r)/2, \delta\psi(r)]$ and $y\in B(x, \psi^{-1}(\epsilon\delta \psi(r)/2))$. Then $s\le \epsilon \psi(r)$ and $\psi^{-1}(\epsilon\delta \psi(r)/2)\le \psi^{-1}(\epsilon s)$ so that $y\in B(x, \psi^{-1}(\epsilon s))$. Hence, by \eqref{e:volume_doubling}, \eqref{e:upper-scaling-psi}, Lemma \ref{l:lower-hke-killed} and Remark \ref{r:lower-hke-killed},
$$
q_{B(x,r)}(s,x,y) \ge c_1 V(x, \psi^{-1}(s))^{-1} \ge c_2 V(x,r)^{-1}.
$$
Therefore
$$
\E^{(t,x)}\left[\int_0^{\tau_r} \1_A (t-s, X_s)ds\right]\ge \frac{c_2}{V(x,r)} \int_{\delta \psi(r)/2}^{\delta\psi(r)} ds \int_{A_{t-s}}dy=c_2\frac{dt \otimes m(A)}{V(x,r)}.
$$
Combining with \eqref{e:ckw-3-7-b}, we arrive at \eqref{e:ckw-3-7-a}. \qed

\begin{thm}\label{t:phr}
There exists a constant $\eta \in (0,1]$ such that for all $\delta \in (0,1)$, there exists a constant $C=C(\delta)>0$ so that for every $x_0 \in D$, $r \in (0, \delta_D(x_0))$,  $t_0\ge 0$, and any  function $u$ on $(0,\infty)\times D$ which is parabolic in $(t_0, t_0+ \psi(r))\times B(x_0, r)$ and bounded in $(t_0, t_0+ \psi(r))\times D$, we have
	\begin{equation}\label{e:phr}
		|u(s,x) - u(t,y)| \le C \Big(\frac{\psi^{-1}(|s-t|) + \rho(x,y)}{r}\Big)^\eta \esssup_{[t_0, t_0 + \psi(r)] \times D}|u|,
	\end{equation}
for every $s,t \in (t_0 + \psi(r)-\psi(\delta r), t_0 + \psi(r))$ and $x,y \in B(x_0, \delta r)$.
\end{thm}
\pf Using \eqref{e:volume_doubling}, \eqref{e:upper-scaling-psi} and Lemmas  \ref{l:exit-time},  \ref{l:Tail-estimate} and \ref{l:ckw-3-7}, the result can be proved 
using the same argument as in the proof of \cite[Theorem 4.14]{CK03} (see also the proof of \cite[Proposition 3.8]{CKW20}). We omit details here. 
\qed

\begin{lemma}\label{l:ckw-4-2}
	Let $\epsilon\in (0,1/4)$ be the constant from Lemma \ref{l:lower-hke-killed} and let $\delta\in (0,\epsilon/4)$ be such that $4\delta\psi(2r)\le \epsilon\psi(r)$ for all $r>0$. Then there exists a constant $C_2>0$ such that  for all $x_0 \in D$, $R \in (0, \delta_D(x_0))$, $r \in  (0, \psi^{-1}(\epsilon\delta \psi(R)/2)/2]$,
	 $\delta\psi(R)/2\le t-s \le 4\delta\psi(2R)$, 
	$x\in B(x_0, \psi^{-1}(\epsilon\delta \psi(R)/2) /2  )$,  and $z\in B(x_0, \psi^{-1}(\epsilon\delta \psi(R)/2) )$, 
	\begin{equation}\label{e:ckw-4-2-a}
	\P^{(t,z)}(\sigma_{ \{s\}\times B(x,r)}\le \tau_{[s,t]\times B(x_0,R)})\ge C_2\frac{V(x,r)}{V(x,R)}. 
	\end{equation}
\end{lemma}
\pf The left-hand side of \eqref{e:ckw-4-2-a} is equal to 
\begin{equation}\label{e:ckw-4-2-b}
	\P^z(X_{t-s}^{B(x_0, R)}\in B(x,r))=\int_{B(x,r)}q_{B(x_0,R)}(t-s,z,y)\, dy.
\end{equation}
Note that $t-s\le 4\delta\psi(2R)\le \epsilon \psi(R)$ by the choice of $\delta$. Next, if $z\in B(x,r)$, then $\rho(z,x_0)\le \rho(z,x)+\rho(x,x_0)\le r+\rho(x,x_0)\le \psi^{-1}(\epsilon\delta\psi(R)/2)\le \psi^{-1}(\epsilon(t-s))$, implying that $z\in B(x_0, \psi^{-1}(\epsilon(t-s)))$ (with the same conclusion for $y\in B(x,r)$). Thus it follows from  Lemma \ref{l:lower-hke-killed}, Remark \ref{r:lower-hke-killed} and \eqref{e:compare_ball} that $q_{B(x_0,R)}(t-s,z,y)\ge c_1 V(x_0,R)^{-1} \ge c_2V(x,R)^{-1}$. 
 By inserting in \eqref{e:ckw-4-2-b} we obtain \eqref{e:ckw-4-2-a}. \qed

 In the remainder of this section, we further assume that $h$ is of Harnack-type.

Suppose that $x,y,z\in D$ are such that $\rho(x,z)\le \rho(x,y)/2$. Then
\begin{equation}\label{e:xyz}
	\frac23 \rho(x,y)\le \rho(z,y)\le \frac32 \rho(x,y).
\end{equation}
As a consequence, by the scalings of $\Phi$ and $\Psi$, there exists $a>1$ such that
\begin{equation}\label{e:phipsi-a}
	a^{-1}\Phi(\rho(x,y))\le \Phi(\rho(z,y))\le a\Phi(\rho(x,y)), \quad a^{-1}\Psi(\rho(x,y))\le \Psi(\rho(z,y))\le a\Psi(\rho(x,y)).
\end{equation}

\begin{prop}\label{p:comp-J}
 Suppose that $h$ is of Harnack-type. Then there exists a constant $C>0$ such that for all $x,y,z\in D$ satisfying $\rho(x,z)\le (\rho(x,y)\wedge \delta_D(x))/2$, it holds that
	$$
	J(x,y)\le C J(z,y).
	$$
\end{prop}
\pf  
This follows 
from Theorem \ref{t:jump-estimate}, \eqref{e:H_3}, the scaling property of $w$, \eqref{e:boundary-function-c} and \eqref{e:phipsi-a}. 
\qed

\begin{cor}\label{c:ujs}  Suppose that $h$ is of Harnack-type. Then there exists a constant $C>0$ such that for all $x,y\in D$ and $0<r \le (\rho(x,y) \wedge \delta_D(x))/2$,
	it holds that
	$$
	J(x,y)\le \frac{C}{V(x,r)}\int_{B(x,r)}J(z,y)dz.
	$$
\end{cor}
\pf Let $x,y\in D$ and $r>0$ be as in the statement. If $z\in B(x,r)$, then $\rho(x,z)<r\le \rho(x,y)/2$. Therefore, by Proposition \ref{p:comp-J}, $J(x,y)\le c_1 J(z,y)$, whence the claim immediately follows. \qed

\begin{lemma}\label{l:ckw-4-1}  Suppose that $h$ is of Harnack-type.  Let $\epsilon\in (0,1/4)$ be the constant from Lemma \ref{l:lower-hke-killed}, 
	and $\theta\ge 1/2$. Further, let $0<\delta_0 <\epsilon$,
	 and $0<\delta_1<\delta_2<\delta_3<\delta_4$ be such that $(\delta_3-\delta_2)\psi(r)\ge \psi(\delta_0 r)$ 
and $\delta_4\psi(r)\le \psi(\epsilon r)$
for all   $r \in (0, \diam(D))$. 
For $x_0\in D$, $t_0\ge 0$ and $r \in (0, \delta_D(x))$, 
define
\begin{align*}
Q_1&=(t_0, t_0+\delta_4 \psi(r))\times B(x_0, \delta_0^2 r), &Q_2&=(t_0, t_0+\delta_4 \psi(r))\times B(x_0, r),\\
Q_3&=[t_0+\delta_1\psi(r), t_0+\delta_2\psi(r)]\times B(x_0, \delta_0^2 r/2), &Q_4&=[t_0+\delta_3\psi(r), t_0+\delta_4\psi(r)]\times B(x_0, \delta_0^2 r/2).
\end{align*}
Let $f:(t_0, \infty)\times D\to [0,\infty)$ be bounded and supported in $(t_0, \infty)\times (D\setminus B(x_0, (1+\theta)r))$. Then there exists a constant $C_2>0$ such that
$$
\E^{(t_1, y_1)}f(Z_{\tau_{Q_1}})\le C_2 \E^{(t_2, y_2)}f(Z_{\tau_{Q_2}})\quad \textrm{for all }(t_1,y_1)\in Q_3, (t_2,y_2)\in Q_4.
$$
\end{lemma}
\pf Without loss of generality we may assume that $t_0=0$. 
  Fix $x_0 \in D$. For $s>0$,
we set $B_s=B(x_0,s)$.  Let  $(t_1,y_1) \in Q_3$ and $(t_2,y_2)\in Q_4$. By the L\'evy system formula for the time-space process $Z$ we have
\begin{eqnarray}
\lefteqn{ \E^{(t_2, y_2)}f(Z_{\tau_{Q_2}})=\E^{(t_2,y_2)}f(t_2-(\tau_{B_r}\wedge t_2), X_{\tau_{B_r}\wedge t_2}) }\nonumber \\
 &=& \E^{(t_2,y_2)}\left[\int_0^{t_2}\1_{t\le \tau_{B_r}}\, dt \int_{D\setminus B(x_0, (1+\theta)r)}f(t_2-t, v)J(X_t,v)\, dv\right] \nonumber \\
 &=&\int_0^{t_2} dt \int_{D\setminus B(x_0, (1+\theta)r)} f(t_2-t,v)  \E^{(t_2, y_2)}\left[\1_{t\le \tau_{B_r}}J(X_t,v)\right]dv \nonumber \\
 &=&\int_0^{t_2} ds \int_{D\setminus B(x_0, (1+\theta)r)} f(s,v)  \E^{(t_2, y_2)}\left[\1_{t_2-s\le \tau_{B_r}}J(X_{t_2-s},v)\right]dv \nonumber \\
 &=&\int_0^{t_2} ds \int_{D\setminus B(x_0, (1+\theta)r)} f(s,v)dv\int_{B_r}q_{B_r}(t_2-s, y_2,z)J(z,v)dz \label{e:ckw-4-1-a}\\
 &\ge & \int_0^{t_1} ds \int_{D\setminus B(x_0, (1+\theta)r)} f(s,v)dv \int_{B_{\delta_0^2 r}}q_{B_r}(t_2-s, y_2,z)J(z,v)dz. \nonumber
\end{eqnarray}
For $s\in [0,t_1]$ it holds that $\psi(\delta_0 r)\le (\delta_3-\delta_2)\psi(r) \le t_2-t_1\le t_2-s \le  \delta_4\psi( r)
\le \psi(\epsilon r)$, hence $\delta_0^2 r\le \epsilon \delta_0 r\le \epsilon \psi^{-1}(t_2-s)$. Therefore,  for any $z \in B_{\delta_0^2r}$, by Lemma \ref{l:lower-hke-killed} and \eqref{e:compare_ball}, 
$q_{B_r}(t_2-s, y_2,z)\ge c_1 V(x_0,\psi^{-1}(t_2-s))^{-1}\ge c_1 V(x_0, \epsilon r)^{-1}$. 
We conclude that
\begin{equation}\label{e:ckw-4-1-b}
 \E^{(t_2, y_2)}f(Z_{\tau_{Q_2}}) \ge \frac{c_1}{V(x_0,\epsilon r)}\int_0^{t_1} ds \int_{D\setminus B(x_0, (1+\theta)r)} f(s,v)dv \int_{B_{\delta_0^2 r}}J(z,v)dz.
\end{equation}
Now, by using the L\'evy system formula again, similar to \eqref{e:ckw-4-1-a} we obtain that
\begin{eqnarray*}
\lefteqn{ \E^{(t_1, y_1)}f(Z_{\tau_{Q_1}})=\int_0^{t_1} ds \int_{D\setminus B(x_0, (1+\theta)r)} f(s,v)dv\int_{B_{\delta_0^2 r}}q_{B_{\delta_0^2 r}}(t_1-s, y_1,z)J(z,v)dz}\\
&=&\int_0^{t_1} ds \int_{B_{\delta_0^2 r}}q_{B_{\delta_0^2 r}}(t_1-s, y_1,z)dz  \int_{D\setminus B(x_0, (1+\theta)r)} f(s,v)J(z,v) dv\\
&=& \int_0^{t_1} ds\left(\int_{B_{\delta_0^2 r}\setminus B_{3\delta_0^2 r/4}}q_{B_{\delta_0^2 r}}(t_1-s, y_1,z)dz  \int_{D\setminus B(x_0, (1+\theta)r)} f(s,v)J(z,v) dv\right.\\
& &\left. +\int_{B_{3\delta_0^2 r/4}}q_{B_{\delta_0^2 r}}(t_1-s, y_1,z)dz  \int_{D\setminus B(x_0, (1+\theta)r)} f(s,v)J(z,v) dv\right)\,=:\,\int_0^{t_1} ds(I_1+I_2).
\end{eqnarray*}
Let $z\in B_{\delta_0^2 r}\setminus B_{3\delta_0^2 r/4}$. Since $y_1\in B_{\delta_0^2 r/2}$, we have that  $ 2\rho(x_0,y_1)  \le \delta_0^2 r/4  \le \rho(y_1,z)\le 3\delta_0^2 r/2$, which implies by \eqref{e:upper-scaling-psi} that 
 $\psi(\rho(y_1,z))\ge c_2 \psi(r) \ge c_2\delta_2^{-1}t_1$. 
 Hence by 
 Proposition \ref{p:hke-upper-rough}, \eqref{e:volume_doubling} and \eqref{e:compare_ball}, 
we get that for any $s>0$,
\begin{align*}
q_{B_{\delta_0^2r}}(t_1-s,y_1,z) &\le q(t_1-s,y_1,z) \le \frac{c_3t_1}{V(y_1, \rho(y_1,z))\psi(\rho(y_1,z))}\le \frac{c_2^{-1} c_3 \delta_2}{V(y_1, \rho(y_1,z))}\\
&\le \frac{c_4}{V(x_0, \rho(y_1,z))} \le \frac{c_4}{V(x_0, \delta_0^2r/4)}\le \frac{c_5}{V(x_0, \epsilon r)}.
\end{align*}
 Therefore,
$$
I_1 \le \frac{c_5}{V(x_0,\epsilon r)}\int_{B_{\delta_0^2 r}\setminus B_{3\delta_0^2 r/4}}dz  \int_{D\setminus B(x_0, (1+\theta)r)} f(s,v)J(z,v) dv.
$$
 Let $z\in B_{3\delta_0^2 r/4}$. Then by  \eqref{e:compare_ball}, we have
\begin{equation}\label{e:ujs_volume}
V(z, \delta_0^2r/4)  \ge c_6V(x_0, \epsilon r).
\end{equation}
We also have that $\delta_D(z)\ge r-3\delta_0^2 r/4 \ge (1-\delta_0^2)r$, and for $v\in D\setminus  B(x_0, (1+\theta)r)$,
$$
\rho(z,v)\ge \rho(x_0, v)-\rho(x_0,z)\ge (1+\theta)r-\frac{3\delta_0^2 r}{4}\ge \frac{3r}{2}-\frac{3\delta_0^2 r}{4} \ge  (1-\delta_0^2)r. 
$$
Thus, for any $w\in B(z, \delta_0^2 r/4)$,  since $\delta_0^2<2/5$,
$$
\frac{1}{2}\left(\rho(z,v)\wedge \delta_D(z)\right)\ge \frac{1}{2}(1-\delta_0^2)r\ge \frac{\delta_0^2 r}{4}\ge \rho(z,w).
$$
From Proposition \ref{p:comp-J}, we get that $J(z,v)\le c_7 J(w,v)$. Therefore, by \eqref{e:ujs_volume},
\begin{eqnarray*}
I_2 &=& \int_{B_{3\delta_0^2 r/4}}q_{B_{\delta_0^2 r}}(t_1-s, y_1,z)dz\int_{D\setminus B(x_0, (1+\theta)r)} f(s,v)J(z,v) dv\\
& \le & \int_{B_{3\delta_0^2 r/4}}q_{B_{\delta_0^2 r}}(t_1-s, y_1,z)dz \left( \frac{c_7}{V(z,\delta_0^2r/4)}\int_{B(z,\delta_0^2 r/4)}J(w,v)dw \int_{D\setminus B(x_0, (1+\theta)r)} f(s,v) dv \right) \\
&\le & \frac{c_6^{-1}c_7}{V(x_0,\epsilon r)}\int_{B(z,\delta_0^2 r/4)}dw \int_{D\setminus B(x_0, (1+\theta)r)} f(s,v) J(w,v) dv.
\end{eqnarray*}
It follows that 
\begin{eqnarray*}
\E^{(t_1, y_1)}f(Z_{\tau_{Q_1}})&=&\int_0^{t_1}ds(I_1+I_2) \le \frac{c_8}{V(x_0,\epsilon r)}\int_0^{t_1}ds \int_{B_{\delta_0^2 r}}dz  \int_{D\setminus B(x_0, (1+\theta)r)} f(s,v)J(z,v) dv\\
&=&\frac{c_8}{V(x_0,\epsilon r)}\int_0^{t_1} ds \int_{D\setminus B(x_0, (1+\theta)r)} f(s,v)dv \int_{B_{\delta_0^2 r}}J(z,v)dz \stackrel{\eqref{e:ckw-4-1-b}}{\le } c_9 \E^{(t_2, y_2)}f(Z_{\tau_{Q_2}}).
\end{eqnarray*}
\qed

\begin{thm}\label{t:phi}
 Suppose that $h$ is of Harnack-type.  Then there exist constants $\delta>0$, $C>1$ and $K\ge 1$ such that for all $t_0\ge 0$, $x_0\in D$ and $R\in (0, R_1)$ with $B(x_0,CR)\subset D$, and any non-negative function $u$ on $(0,\infty)\times D$ which is parabolic on $Q:=(t_0, t_0+4\delta \psi(CR))\times B(x_0, CR)$, we have
\begin{equation}\label{e:phi-1}
\sup_{(t_1,y_1)\in Q_{-}}u(t_1,y_1)\le K \inf_{(t_2, y_2)\in Q_{+}}u(t_2, y_2),
\end{equation}
where $Q_-=[t_0+\delta\psi(CR), t_0+2\delta\psi(CR)]\times B(x_0, R)$ and $Q_-=[t_0+3\delta\psi(CR), t_0+4\delta\psi(CR)]\times B(x_0, R)$.
\end{thm}
\pf Using \eqref{e:volume_doubling}, \eqref{e:upper-scaling-psi}, Lemmas \ref{l:lower-hke-killed},   \ref{l:exit-time},  \ref{l:ckw-3-7}, \ref{l:ckw-4-2}, \ref{l:ckw-4-1} and Corollary \ref{c:ujs},  the result can be proved 
using the same argument as in the proof of \cite[Lemma 5.3]{CKK09} (see also the proof of \cite[Lemma 4.1]{CKW20}). We omit details here. 
\qed

%%%%%%%%%%%%%%%%%%%%%%%%%%%%%%%%%%%%%%%%%%%%%%%%%%%%%%%%%%%%%%%%%%%%%%%%%%%%%%%%%%%
%%%%%%%%%%%%%%%%%%%%%%%%%%                         Examples %%%%%%%%%%%%%%%%%%%%%%%%%%%%
%%%%%%%%%%%%%%%%%%%%%%%%%%%%%%%%%%%%%%%%%%%%%%%%%%%%%%%%%%%%%%%%%%%%%%%%%%%%%%%%%%%
\section{Examples}\label{s:examples}
 Recall the  definitions  of $h_{p,q}(t,x,y)$  from \eqref{e:def-hp}, and  $\sB_h(t,x,y)$ from \eqref{e:defsB}. We remind the reader that $h_{p,q}(t,x,y)$ is quite typical and it is the most important boundary function. 
Recall that  
$\psi(r)=1/\phi(1/\Phi(r))$, 
	$\delta_{\vee}(x,y)= \delta_D(x) \vee \delta_D(y)$
and	$\delta_{\wedge}(x,y)= \delta_D(x) \wedge \delta_D(y).$
For simplicity, we will use $\delta(x)$ and $\delta(y)$ instead of  $\delta_D(x)$ and  $\delta_D(y)$, respectively.

We let
\begin{align}\label{e:time-regularization}
	\delta^t(x)&:=\delta_D(x)\vee \psi^{-1}(t),\nn\\
	\delta^t_\vee(x,y)&:=\delta^t(x) \vee \delta^t(y)=\delta_\vee(x,y) \vee \psi^{-1}(t),\nn\\
	\delta^t_\wedge(x,y)&:=\delta^t(x) \wedge \delta^t(y)=\delta_\wedge(x,y) \vee \psi^{-1}(t).
\end{align}

The following lemma provides 
 a 
list of estimates of  $\sB_{h_{p,q}}(t,x,y)$ 
 depending on the relationship between 
the parameters $p,q, \beta_1, \beta_2$. 
The list is not exhaustive, but it suffices for our purpose. The proof of the lemma is rather technical and consists 
of carefully estimating the integral defining  $\sB_{h_{p,q}}(t,x,y)$. 
  The factorization	\eqref{e:standard-factorization} below is inspired by \eqref{e:on} and \eqref{e:simple-off0}. See also \eqref{e:heatkernel-estimates} below.

\begin{lemma}\label{l:asym-Bp}
	Let $q\ge p\ge0$, $p+q>0$. Suppose that {\bf (Poly-$R_1$)}\ holds with $\beta_1, \beta_2\in(0,1)$. 
	Then 
\begin{equation}\label{e:standard-factorization}
\sB_{h_{p,q}}(t,x,y)= \frac{\Phi(\rho(x,y))}{\psi(\rho(x,y))} \Big(1 \wedge \frac{\Phi(\delta_D(x))}{\phi^{-1}(1/t)^{-1}} \Big)^{p}\Big(1 \wedge \frac{\Phi(\delta_D(y))}{\phi^{-1}(1/t)^{-1}} \Big)^{q} A_{p,q}(t,x,y),
\end{equation} 
where $A_{p,q}(t,x,y)$ satisfies the following   estimates for all $x,y \in D$ and  $0<t \le  \psi(\rho(x,y))$ 
 such that $\Phi(\rho(x,y))<R_1/8$. 	
	\smallskip
	
	\noindent (i) If $\beta_2<1-p-q$, then
	\begin{equation*}
		A_{p,q}(t,x,y) \simeq \Big(1 \wedge \frac{\Phi(\delta^t(y))}{\Phi(\rho(x,y))} \Big)^{q-p}\Big(1 \wedge \frac{\Phi(\delta^t_\wedge(x,y))}{\Phi(\rho(x,y))} \Big)^{p} \Big(1 \wedge \frac{\Phi(\delta^t_\vee(x,y))}{\Phi(\rho(x,y))} \Big)^{p}.
	\end{equation*}

\noindent (ii) If $1-p-q<\beta_1 \le \beta_2<1- q$, then
\begin{align*}
	A_{p,q}(t,x,y)&\simeq \Big(1 \wedge \frac{\Phi(\delta^t(y))}{\Phi(\rho(x,y))} \Big)^{q-p} \Big(1 \wedge \frac{\Phi(\delta^t_\wedge(x,y))}{\Phi(\rho(x,y))} \Big)^{p} \Big(1 \wedge \frac{\Phi(\delta^t_\vee(x,y))}{\Phi(\rho(x,y))} \Big)^{1-q} \Big(1\wedge \frac{\psi(\delta^t_\vee(x,y))}{\psi(\rho(x,y))} \Big)^{-1} .
\end{align*}

\noindent (iii) If $1- q<\beta_1 \le \beta_2<1- p $, then
\begin{align*}
	A_{p,q}(t,x,y)&\simeq   \Big(1 \wedge \frac{\Phi(\delta^t(y))}{\Phi(\rho(x,y))} \Big)^{1-p}\Big(1 \wedge \frac{\Phi(\delta^t_\wedge(x,y))}{\Phi(\rho(x,y))} \Big)^{p}  \Big(1\wedge \frac{\psi(\delta^t(y))}{\psi(\rho(x,y))} \Big)^{-1}.
\end{align*}

\noindent (iv) If $1-p<\beta_1 \le \beta_2<1$, then
\begin{equation*}
	A_{p,q}(t,x,y)\simeq \Big(1 \wedge \frac{\Phi(\delta^t_\wedge(x,y))}{\Phi(\rho(x,y))} \Big)  \Big(1\wedge \frac{\psi(\delta^t_\wedge(x,y))}{\psi(\rho(x,y))} \Big)^{-1}.
\end{equation*}

\noindent (v) If $\beta_1=\beta_2=1-p-q$ and $p>0$, then
\begin{equation*}
	A_{p,q}(t,x,y) \simeq \Big(1 \wedge \frac{\Phi(\delta^t(y))}{\Phi(\rho(x,y))} \Big)^{q-p}  \Big(1 \wedge \frac{\Phi(\delta^t_\wedge(x,y))}{\Phi(\rho(x,y))} \Big)^{p} \Big(1 \wedge \frac{\Phi(\delta^t_\vee(x,y))}{\Phi(\rho(x,y))} \Big)^{p}  \log\Big( e+ \frac{\Phi(\rho(x,y)) }{ \Phi(\delta^t_{\vee}(x,y)) }\Big). 	\end{equation*}

	\noindent (vi) If $\beta_1=\beta_2=1-p-q$ and $p=0$, then
	\begin{align*}
		A_{p,q}(t,x,y) &\simeq \Big(1 \wedge \frac{\Phi(\delta^t(y))}{\Phi(\rho(x,y))} \Big)^q \log \Big(e+\frac{\Phi(\rho(x,y))}{\Phi(\delta^t(y))}\Big).
	\end{align*}

\noindent (vii) If $\beta_1=\beta_2=1-q$ and $q=p$, then
\begin{align*}
	A_{p,q}(t,x,y) &\simeq \Big(1 \wedge \frac{\Phi(\delta^t_\wedge(x,y))}{\Phi(\rho(x,y))} \Big)^{p} \log \Big( e+ \frac{\Phi(\rho(x,y)) \wedge \Phi(\delta^t_{\vee}(x,y))}{\Phi(\delta^t_\wedge(x,y))}\Big).
\end{align*}

	\noindent (viii) If $\beta_1=\beta_2=1- q$ and $q>p>0$, then
	\begin{align*}
		A_{p,q}(t,x,y)&\simeq   \Big(1 \wedge \frac{\Phi(\delta^t(y))}{\Phi(\rho(x,y))} \Big)^{q-p}\Big(1 \wedge \frac{\Phi(\delta^t_\wedge(x,y))}{\Phi(\rho(x,y))} \Big)^{p}  \log \Big( e+ \frac{\Phi(\rho(x,y))\wedge \Phi(\delta^t(x))}{\Phi(\delta^t(y))}\Big).
	\end{align*}

	\noindent (ix) If $\beta_1=\beta_2=1-p$ and $q>p$, then
	\begin{equation*}
		A_{p,q}(t,x,y) \simeq \Big(1 \wedge \frac{\Phi(\delta^t_\wedge(x,y))}{\Phi(\rho(x,y))} \Big)^{p}  
		\log\Big(e+\frac{\Phi(\rho(x,y)) \wedge \Phi(\delta^t(y))}{\Phi(\delta^t(x))}\Big).
	\end{equation*}
\end{lemma}
\pf Fix $x,y \in D$ such that $r:=\rho(x,y)<\Phi^{-1}(R_1/8)$, and  $t \in (0, \psi(r)]$. 
Let $\delta_\wedge:=\delta_\wedge(x,y)$ and $\delta_\vee:=\delta_\vee(x,y)$, 
and note that $\delta_{\vee}\le r+\delta_{\wedge}$.  We also note that  since $\beta_2<1$, by  \cite[Lemma 2.6, Proposition 2.9]{Mi} and  \eqref{e:phi-upper-scaling},
\begin{equation}\label{e:1-w-phi}
	w(s) \simeq \phi(1/s), \quad s\in (0,R_1/2),
\end{equation}  
which is equivalent to that $w(\Phi(s)) \simeq \psi(s)^{-1}$ for $s \in (0, \Phi^{-1}(R_1/2))$. 

 Define 
\begin{equation}\label{e:standard-factorization-alt}
A_{p,q}(t,x,y):=\frac{\psi(r)}{\Phi(r)}\frac{\sB_{h_{p,q}}(t,x,y)}{h_{p,q}(\phi^{-1}(1/t)^{-1},x,y)}.
\end{equation}
 Observe that for any $s>\phi^{-1}(1/t)^{-1}$,
\begin{align}\label{e:equiv-form}
	\frac{h_{p,q}(s,x,y)}{h_{p,q}(\phi^{-1}(1/t)^{-1},x,y)} &= \Big(1 \wedge \frac{\Phi(\delta^t(x))}{s} \Big)^p \Big(1 \wedge \frac{\Phi(\delta^t(y))}{s} \Big)^q\nn\\
	&= \Big(1 \wedge \frac{\Phi(\delta^t(y))}{s} \Big)^{q-p} \Big(1 \wedge \frac{\Phi(\delta^t_\wedge(x,y))}{s} \Big)^p  \Big(1 \wedge \frac{\Phi(\delta^t_\vee(x,y))}{s} \Big)^p. 
\end{align}
Hence, by the definition of $\sB_{ h_{p,q}}(t,x,y)$  and \eqref{e:standard-factorization-alt},  we get
\begin{align*}
	A_{p,q}(t,x,y)=\frac{\psi(r)}{\Phi(r)}\int_{2\phi^{-1}(1/t)^{-1}}^{4\Phi(r)} \Big(1 \wedge \frac{\Phi(\delta^t(y))}{s} \Big)^{q-p} \Big(1 \wedge \frac{\Phi(\delta^t_\wedge(x,y))}{s} \Big)^p  \Big(1 \wedge \frac{\Phi(\delta^t_\vee(x,y))}{s} \Big)^p w(s)ds.
\end{align*}

If $\delta_{\vee} \ge 2r$, then $\delta_\wedge \ge \delta_\vee -r \ge r$ and hence we get from  Lemma \ref{l:intscale}(i), {\bf (Poly-$R_1$)} 
 (by using $\beta_2<1$) 
and \eqref{e:1-w-phi} that
\begin{equation*}
	A_{p,q}(t,x,y)\simeq \frac{\psi(r)}{\Phi(r)}\int_{2\phi^{-1}(1/t)^{-1}}^{4\Phi(r)}  w(s)ds \simeq \psi(r) w(4\Phi(r)) \simeq 1.
\end{equation*}
 On the other hand, we see that  in every of the cases (i)-(ix), the right-hand side of comparability relation for
$A_{p,q}(t,x,y)$ is comparable to $1$. Hence the assertion of the lemma is valid when $\delta_\vee \ge 2r$.

Suppose  now  that $\delta_\vee<2r$. 
 Since $\psi^{-1}(t)<r$, we get $\delta^t_\vee(x,y)<2r$, hence by 
 the scaling property \eqref{e:psiphi} of $\Phi$, we have
\begin{equation}\label{e:delta<r}
	3\Phi(a) \wedge 4\Phi(r) \simeq \Phi(a) \quad \text{for} \;\; a\in \{\delta^t(y), \, \delta^t_\wedge(x,y),  \,\delta^t_\vee(x,y)\}.
\end{equation} 

\smallskip

\noindent (i) 
 The desired comparability 
follows  from Lemma \ref{l:intscale}(i), {\bf (Poly-$R_1$)} and \eqref{e:1-w-phi}.

\noindent (ii) By  Lemma \ref{l:intscale}(i)-(ii), {\bf (Poly-$R_1$)} and \eqref{e:1-w-phi}, since $q-p+p+\beta_2<1<p+q+\beta_1$,
\begin{align}\label{e:allow-1}
	A_{p,q}(t,x,y) &\simeq \frac{\psi(r)}{\Phi(r)}\int_{2\phi^{-1}(1/t)^{-1}}^{3\Phi(\delta^t_\vee(x,y)) \wedge 4\Phi(r)} \Big(1 \wedge \frac{\Phi(\delta^t(y))}{s} \Big)^{q-p} \Big(1 \wedge \frac{\Phi(\delta^t_\wedge(x,y))}{s} \Big)^p   w(s)ds\nn\\
	& \quad  + \frac{\psi(r)}{\Phi(r)}\Phi(\delta^t(y))^{q-p}\Phi(\delta^t_\wedge(x,y))^p \Phi(\delta^t_\vee(x,y))^p\int_{3\Phi(\delta^t_\vee(x,y)) \wedge 4\Phi(r)}^{4\Phi(r)}  s^{-p-q}w(s)ds \nn\\
	&\simeq \frac{\psi(r)}{\Phi(r)}  \Big(1 \wedge \frac{\Phi(\delta^t(y))}{\Phi(\delta^t_\vee(x,y))} \Big)^{q-p} \Big(1 \wedge \frac{\Phi(\delta^t_\wedge(x,y))}{
		\Phi(\delta^t_\vee(x,y))} \Big)^p   \Phi(\delta^t_\vee(x,y)) w(\Phi(\delta^t_\vee(x,y))) \nn\\
	&\simeq \frac{\Phi(\delta^t(y))^{q-p} \Phi(\delta^t_\wedge(x,y))^{p}}{ \Phi(r)^q} \, \Big(\frac{\Phi(r)}{ \Phi(\delta^t_\vee(x,y))} \Big)^q \, \frac{\Phi(\delta^t_\vee(x,y))}{\psi(\delta^t_\vee(x,y))} \frac{\psi(r)}{\Phi(r)}.
\end{align}
 We used \eqref{e:delta<r} in the second comparability above. 
Since $\Phi$ and $\psi$ are increasing and satisfy scaling properties, the desired comparability holds.

\noindent (iii)  If $\delta^t_{\vee}(x,y) =\delta^t(y)$, then \eqref{e:allow-1} holds by Lemma \ref{l:intscale}(i)-(ii), {\bf (Poly-$R_1$)}, \eqref{e:1-w-phi} and \eqref{e:delta<r} since  $p+\beta_2<1<p+q+\beta_1$. 
 The desired comparability then follows immediately from \eqref{e:allow-1}. 
If $\delta^t_{\wedge}(x,y) =\delta^t(y)$, then 
  we get, by Lemma \ref{l:intscale}(i)-(ii), {\bf (Poly-$R_1$)}, \eqref{e:1-w-phi}, \eqref{e:delta<r} and the assumption  $q+\beta_1>1$, that 
\begin{align*}
	A_{p,q}(t,x,y)&\simeq \frac{\psi(r)}{\Phi(r)}\int_{2\phi^{-1}(1/t)^{-1}}^{3\Phi(\delta^t(y)) \wedge 4\Phi(r)}  w(s)ds+ \frac{\psi(r)}{\Phi(r)} \Phi(\delta^t(y))^q \int_{3\Phi(\delta^t(y)) \wedge 4\Phi(r)}^{4\Phi(r)}   \Big(1 \wedge \frac{\Phi(\delta^t_\vee(x,y))}{s} \Big)^p s^{-q} w(s)ds\nn\\
	&\simeq  \frac{\Phi(\delta^t(y))}{\psi(\delta^t(y))} \frac{\psi(r)}{\Phi(r)} =\Big( \frac{\Phi(\delta^t(y))}{\Phi(r)}\Big)^{1-p} \Big( \frac{\Phi(\delta^t_\wedge(x,y))}{\Phi(r)}\Big)^{p}  \frac{\psi(r)}{\psi(\delta^t(y))}. 
\end{align*}

\noindent (iv)  Since $p+\beta_1>1$, we get from  Lemma \ref{l:intscale}(i)-(ii), {\bf (Poly-$R_1$)}, \eqref{e:1-w-phi} and \eqref{e:delta<r} that
\begin{align*}
A_{p,q}(t,x,y) &\simeq \frac{\psi(r)}{\Phi(r)}\int_{2\phi^{-1}(1/t)^{-1}}^{3\Phi(\delta^t_\wedge(x,y)) \wedge 4\Phi(r)}  w(s)ds\nn\\
&\quad +  \frac{\psi(r)}{\Phi(r)} \Phi(\delta^t_\wedge(x,y))^p\int_{3\Phi(\delta^t_\wedge(x,y)) \wedge 4\Phi(r)}^{4\Phi(r)} \Big(1 \wedge \frac{\Phi(\delta^t(y))}{s} \Big)^{q-p}  \Big(1 \wedge \frac{\Phi(\delta^t_\vee(x,y))}{s} \Big)^p s^{-p}w(s)ds\nn\\
&\simeq \frac{\Phi(\delta^t_\wedge(x,y))}{\psi(\delta^t_\wedge(x,y))} \frac{\psi(r)}{\Phi(r)}.
\end{align*}

\noindent (v) Note that $w(s) \simeq s^{p+q-1}$ and $\psi(s) \simeq \Phi(s)^{1-p-q}$ for $s \in (0,R_1)$  in this case. Thus, we get  from  Lemma \ref{l:intscale}(i) and \eqref{e:delta<r} that
\begin{align*}
	A_{p,q}(t,x,y) &\simeq \frac{\psi(r)}{\Phi(r)}\int_{2\phi^{-1}(1/t)^{-1}}^{3\Phi(\delta^t_\vee(x,y)) \wedge 4\Phi(r)}  \big(s \wedge \Phi(\delta^t(y)) \big)^{q-p} \big(s \wedge \Phi(\delta^t_\wedge(x,y)) \big)^p    s^{p-1}ds\\
	&\quad + \frac{\psi(r)}{\Phi(r)}\Phi(\delta^t(y))^{q-p}\Phi(\delta^t_\wedge(x,y))^p\Phi(\delta^t_\vee(x,y))^p\int_{3\Phi(\delta^t_\vee(x,y)) \wedge 4\Phi(r)}^{4\Phi(r)}  s^{-1}ds\\
	&\simeq \frac{\Phi(\delta^t(y))^{q-p}\,\Phi(\delta^t_\wedge(x,y))^p\,\Phi(\delta^t_\vee(x,y))^p }{\Phi(r)^{p+q}} \log\Big( e+ \frac{\Phi(r) }{ \Phi(\delta^t_{\vee}(x,y)) }\Big).
\end{align*}

\noindent (vi) Note that $w(s) \simeq s^{q-1}$ and $\psi(s) \simeq \Phi(s)^{1-q}$ for $s \in (0,R_1/2)$  in this case. We get  from  \eqref{e:delta<r} that
\begin{align*}
	A_{p,q}(t,x,y) &\simeq \frac{\psi(r)}{\Phi(r)}\int_{2\phi^{-1}(1/t)^{-1}}^{3\Phi(\delta^t(y)) \wedge 4\Phi(r)}   s^{q-1}ds +  \frac{\psi(r)}{\Phi(r)}  \Phi(\delta^t(y))^q\int_{3\Phi(\delta^t(y)) \wedge 4\Phi(r)}^{4\Phi(r)}  s^{-1}ds\nn\\
	&\simeq   \frac{\Phi(\delta^t(y))^{q}}{\Phi(r)^q}\log\Big( e+ \frac{\Phi(r) }{ \Phi(\delta^t(y)) }\Big).
\end{align*}

\noindent (vii) Note that $w(s) \simeq s^{p-1}$ and $\psi(s) \simeq \Phi(s)^{1-p}$ for $s \in (0,R_1)$  in this case. By \eqref{e:delta<r}, we obtain
\begin{align*}
	A_{p,q}(t,x,y) &\simeq \frac{\psi(r)}{\Phi(r)}\int_{2\phi^{-1}(1/t)^{-1}}^{3\Phi(\delta^t_\wedge(x,y)) \wedge 4\Phi(r)}   s^{p-1}ds+ \frac{\psi(r)}{\Phi(r)} \Phi(\delta^t_\wedge(x,y))^p \int_{3\Phi(\delta^t_\wedge(x,y)) \wedge 4\Phi(r)}^{3\Phi(\delta^t_\vee(x,y)) \wedge 4\Phi(r)}  s^{-1}ds\\
	&\quad + \frac{\psi(r)}{\Phi(r)} \Phi(\delta^t_\wedge(x,y))^p \Phi(\delta^t_\vee(x,y))^p \int_{3\Phi(\delta^t_\vee(x,y)) \wedge 4\Phi(r)}^{4\Phi(r)} s^{-1-p}ds\\
	& \simeq\frac{\Phi(\delta^t_\wedge(x,y))^p}{\Phi(r)^p}  \log\Big( e+ \frac{\Phi(\delta^t_\vee(x,y)) }{ \Phi(\delta^t_\wedge(x,y)) }\Big).
\end{align*}

\noindent (viii) Note that $w(s) \simeq s^{q-1}$ and $\psi(s) \simeq \Phi(s)^{1-q}$ for $s \in (0,R_1)$  in this case. If $\delta^t_\vee(x,y)=\delta^t(y)$, then we get  from  Lemma \ref{l:intscale}(i) and \eqref{e:delta<r} that
\begin{align*}
	A_{p,q}(t,x,y) &\simeq \frac{\psi(r)}{\Phi(r)}\int_{2\phi^{-1}(1/t)^{-1}}^{3\Phi(\delta^t(y)) \wedge 4\Phi(r)}    \big(s \wedge \Phi(\delta^t(x)) \big)^p   s^{q-p-1}ds\\
	&\quad + \frac{\psi(r)}{\Phi(r)}\Phi(\delta^t(x))^p\Phi(\delta^t(y))^q
	\int_{3\Phi(\delta^t(y)) \wedge 4\Phi(r)}^{4\Phi(r)}  s^{-p-1}ds\\
	&\simeq \frac{\Phi(\delta^t(x))^p\,\Phi(\delta^t(y))^{q-p} }{\Phi(r)^{q}} =  \frac{\Phi(\delta^t(y))^{q-p} \,\Phi(\delta^t_\wedge(x,y))^p }{\Phi(r)^{q}}.
\end{align*}
 If $\delta^t_\wedge(x,y)=\delta^t(y)$, then we get  from  \eqref{e:delta<r} that
\begin{align*}
	A_{p,q}(t,x,y) &\simeq \frac{\psi(r)}{\Phi(r)}\int_{2\phi^{-1}(1/t)^{-1}}^{3\Phi(\delta^t(y)) \wedge 4\Phi(r)}   s^{q-1}ds+ \frac{\psi(r)}{\Phi(r)} \Phi(\delta^t(y))^q \int_{3\Phi(\delta^t(y)) \wedge 4\Phi(r)}^{3\Phi(\delta^t(x)) \wedge 4\Phi(r)}  s^{-1}ds\\
	&\quad + \frac{\psi(r)}{\Phi(r)} \Phi(\delta^t(x))^p \Phi(\delta^t(y))^q \int_{3\Phi(\delta^t(x)) \wedge 4\Phi(r)}^{4\Phi(r)} s^{-1-p}ds\\
	& \simeq\frac{\Phi(\delta^t(y))^q}{\Phi(r)^q}  \log\Big( e+ \frac{\Phi(\delta^t(x)) }{ \Phi(\delta^t(y)) }\Big) = \frac{\Phi(\delta^t(y))^{q-p}\Phi(\delta^t_\wedge(x,y))^p}{\Phi(r)^q}  \log\Big( e+ \frac{\Phi(\delta^t(x)) }{ \Phi(\delta^t(y)) }\Big).
\end{align*}

\noindent (ix) Note that $w(s) \simeq s^{p-1}$ and $\psi(s) \simeq \Phi(s)^{1-p}$ for $s \in (0,R_1)$  in this case. If $\delta^t_\vee(x,y)=\delta^t(y)$, then by \eqref{e:delta<r}, it holds that
\begin{align*}
	A_{p,q}(t,x,y) &\simeq \frac{\psi(r)}{\Phi(r)}\int_{2\phi^{-1}(1/t)^{-1}}^{3\Phi(\delta^t(x)) \wedge 4\Phi(r)}   s^{p-1}ds+ \frac{\psi(r)}{\Phi(r)} \Phi(\delta^t(x))^p \int_{3\Phi(\delta^t(x)) \wedge 4\Phi(r)}^{3\Phi(\delta^t(y)) \wedge 4\Phi(r)}  s^{-1}ds\\
	&\quad + \frac{\psi(r)}{\Phi(r)} \Phi(\delta^t(x))^p \Phi(\delta^t(y))^q \int_{3\Phi(\delta^t(y)) \wedge 4\Phi(r)}^{4\Phi(r)} s^{-1-q}ds\\
	& \simeq\frac{\Phi(\delta^t(x))^p}{\Phi(r)^p}  \log\Big( e+ \frac{\Phi(\delta^t(y)) }{ \Phi(\delta^t(x)) }\Big) =\frac{\Phi(\delta^t_\wedge(x,y))^p}{\Phi(r)^p}  \log\Big( e+ \frac{\Phi(\delta^t(y)) }{ \Phi(\delta^t(x)) }\Big).
\end{align*}
If $\delta^t_\wedge(x,y)=\delta^t(y)$, then we get  from Lemma \ref{l:intscale}(ii) and \eqref{e:delta<r} that
\begin{align*}
	A_{p,q}(t,x,y) &\simeq \frac{\psi(r)}{\Phi(r)}\int_{2\phi^{-1}(1/t)^{-1}}^{3\Phi(\delta^t(y)) \wedge 4\Phi(r)}   s^{p-1}ds+ \frac{\psi(r)}{\Phi(r)} \Phi(\delta^t(y))^q \int_{3\Phi(\delta^t(y)) \wedge 4\Phi(r)}^{ 4\Phi(r)} \Big(1 \wedge \frac{\Phi(\delta^t(x))}{s} \Big)^{p} s^{p-q-1}ds\\
	& \simeq\frac{\Phi(\delta^t(y))^p}{\Phi(r)^p} = \frac{\Phi(\delta^t_\wedge(x,y))^p}{\Phi(r)^p}.
\end{align*}

 \qed

\begin{example}\label{ex:basic}
	{\rm Let $d,\alpha>0$, $\beta \in (0,1)$ and
		 $p,q \ge 0$ such that $p+q>0$. 
		  Suppose that 	for every $r_0 \ge 1$,  there are comparability constants such that
		\begin{equation}\label{e:d_set}
		V(x,r) \simeq r^d, \quad x \in D, \; 0<r<r_0.
		\end{equation}
	Let $Y^D$ be a Hunt process in $D$ and $S=(S_t)_{t\ge 0}$ be 
an independent driftless subordinator with Laplace exponent $\phi$.
Suppose that the tail $w$ of the L\'evy measure of $S$ satisfies
    \begin{equation}\label{e:beta_stable_like}
    w(r) \simeq r^{-\beta}, \quad 0<r<r_1,
    \end{equation}
for some $r_1>0$. Suppose  that  the heat kernel $p_D(t,x,y)$ of $Y^D$ satisfies either
${\bf HK^{h_{p,q}}_B}$ or  ${\bf HK^{h_{p,q}}_U}$ 
with $\Phi(r)=\Psi(r)=r^\alpha$  where the boundary 
function $h_{p,q}$ is defined as \eqref{e:def-hp}. When  ${\bf HK^{h_{p,q}}_U}$ 
is satisfied, we also assume that \eqref{e:d_set} and \eqref{e:beta_stable_like} hold  for all $r>0$.  See	Example \ref{ex:HK} for concrete examples of $Y^D$. By switching the roles of $x$ and $y$ if needed, without loss of generality, we assume that $q \ge p$.

Let $q(t,x,y)$, $J(x,y)$ and  $G_D(x,y)$ be the  heat kernel, the jump kernel and the Green function   of the subordinate process $X_t:=Y^D_{S_t}$ respectively. 
Using our theorems in 
Sections \ref{s:key-estimates}  and \ref{s:green},
and Lemma	\ref{l:asym-Bp},
we get explicit estimates on  $q(t,x,y)$, $J(x,y)$ and $G_D(x,y)$. 
We list them in terms of the range of 
$p+q$,
similar to the format of the  Green function estimates for Dirichlet forms degenerate at the boundary in \cite{KSV21}.

		 In particular, by putting 
		 $p=q=1/2$, we get  Theorem \ref{t:special}.
		
		\smallskip
		
	 We first give the Green function estimates. Define
	 \begin{equation*}
	 \mathfrak	g(x,y):= \Big(1 \wedge \frac{\delta(x)}{\rho(x,y)} \Big)^{\alpha p} \Big(1 \wedge \frac{\delta(y)}{\rho(x,y)} \Big)^{\alpha q} \times 
	 \begin{cases}
 \rho(x,y)^{\alpha\beta-d} 	
	 		,&  d>\alpha\beta,\\[4pt]
		\displaystyle	 \log\Big( e +  \frac{\delta_\vee(x,y)}{\rho(x,y)}\Big), &  d=\alpha\beta, \\[7pt]
			 	[\rho(x,y) \vee \delta_\vee(x,y)]^{\alpha\beta-d}
	 		, & d<\alpha\beta.
	 \end{cases}
	 \end{equation*}
When $C_0=0$, by  Theorem \ref{t:Sgreen} and  Example \ref{ex:H2_opposite},  for all $x,y \in D$, if $d>\alpha(\beta-p-q)$, then 
\begin{equation}\label{e:newrszv1}
			G_D(x,y)\simeq \mathfrak g(x,y)
\end{equation}
	and if $d \le \alpha(\beta-p-q)$, then 
\begin{equation}\label{e:newrszv2}
		G_D(x,y)\simeq
		\begin{cases}
		\delta(x)^{\alpha p} \delta(y)^{\alpha q}\log\left( e +  \frac{\diam(D)}{\rho(x,y) \vee \delta_\vee(x,y)}\right), &  d=\alpha(\beta-p-q)
		\text{ and ${\bf HK^{h_{p,q}}_B}$ holds},\\[6pt]
			\delta(x)^{\alpha p} \delta(y)^{\alpha q}, &  d<\alpha(\beta-p-q)\text{ and ${\bf HK^{h_{p,q}}_B}$ holds}, \\[4pt]
			\infty, &d\le \alpha(\beta-p-q)  \text{ and ${\bf HK^{h_{p,q}}_U}$ holds}.	
		\end{cases}
		\end{equation}
				Now assume that $C_0=1$.	If $p+q<\beta+1$, then using Theorem \ref{t:Sgreen} and  Example \ref{ex:H2_opposite} again, we see that 
\eqref{e:newrszv1} and 	\eqref{e:newrszv2} also hold.
		If $p+q=\beta+1$ and 
		$q<\beta+1$ 
(so that (H2**) holds, cf.~Example \ref{ex:H2**}),  
		then by Theorem \ref{t:Sgreen_3}, for all $x,y \in D$,
		\begin{align}\label{e:newrszv3}
			G_D(x,y)\simeq  \mathfrak g(x,y)\log\Big( e +  \frac{\rho(x,y)}{\delta_\vee(x,y)}\Big).
		\end{align}
	 If $p+q>\beta+1$  and 
	 $q<\beta+1$ 
(so, again, (H2**) holds), 
	 then by Theorem \ref{t:Sgreen_2} and \eqref{e:bracket_hpq},  for all $x,y \in D$,
		\begin{align}\label{e:newrszv4}
			G_D(x,y)&\simeq \Big(1 \wedge \frac{\delta_\vee(x,y)}{\rho(x,y)} \Big)^{-\alpha (p+q-\beta-1)}\mathfrak g(x,y).
		\end{align}
The unusual form of the estimates in \eqref{e:newrszv3}-\eqref{e:newrszv4} 
should be compared with similar estimates of the Green function obtained in a different context in \cite[Theorem 1.1 (2),(3)]{KSV21}. Such estimates lead to anomalous boundary behavior of the corresponding Green potential, cf.~\cite{AGV}.

	Next, we get heat kernel estimates from Lemma
	\ref{l:asym-Bp},  Corollary \ref{c:interior-small} and Theorem \ref{t:Slarge}. By the definition \eqref{e:time-regularization}, since $\psi(r) \simeq r^{\alpha\beta}$ for $r \in (0,1]$, we have for all  $(t,x,y) \in (0,1] \times D \times D$,
	\begin{equation}\label{e:time-regularization-1}
		\delta^t(x)\simeq \delta_D(x)\vee t^{1/(\alpha\beta)}, \qquad 
		\delta^t_\vee(x,y)\simeq \delta_\vee(x,y) \vee t^{1/(\alpha\beta)},\qquad 
		\delta^t_\wedge(x,y)\simeq \delta_\wedge(x,y) \vee t^{1/(\alpha\beta)}.
	\end{equation}
Moreover, when ${\bf HK^{h_{p}}_U}$ holds, since $\psi(r) \simeq r^{\alpha\beta}$ for all $r>0$ in this case, \eqref{e:time-regularization-1} holds for all $(t,x,y) \in (0,\infty) \times D \times D$.

	\medskip
	
	\textit{(Small time estimates)} 
	The following estimates hold for all $(t,x,y) \in (0,1] \times D \times D$ :
	\begin{align}\label{e:heatkernel-estimates}
			q(t,x,y)&\simeq \Big(1 \wedge \frac{\delta(x)}{t^{1/(\alpha\beta)}} \Big)^{\alpha p} \Big(1 \wedge \frac{\delta(y)}{t^{1/(\alpha\beta)}} \Big)^{\alpha q} B_{p,q}(t,x,y ;C_0) \left( t^{-d/(\alpha\beta)} \wedge \frac{t  }{\rho(x,y)^{d+\alpha\beta}}\right)\nn\\
			 &\simeq \Big(1 \wedge \frac{\delta(x)}{t^{1/(\alpha\beta)}} \Big)^{\alpha p} \Big(1 \wedge \frac{\delta(y)}{t^{1/(\alpha\beta)}} \Big)^{\alpha q}  \left( t^{-d/(\alpha\beta)} \wedge \frac{t  B_{p,q}(t,x,y ;C_0)}{\rho(x,y)^{d+\alpha\beta}}\right),
	\end{align}
	where
	\begin{equation}\label{e:explicit-boundary-1}
	B_{p,q}(t,x,y ;0):= \Big(1\wedge \frac{\delta^t_\wedge(x,y)}{\rho(x,y)} \Big)^{\alpha p}  \Big(1\wedge \frac{\delta^t_\vee(x,y)}{\rho(x,y)} \Big)^{\alpha p}  \Big(1\wedge \frac{\delta^t(y)}{\rho(x,y)} \Big)^{\alpha (q-p)}
	\end{equation}
	and
	\begin{align}\label{e:explicit-boundary-2}
		&B_{p,q}(t,x,y;1):= \nn\\
		&\begin{cases} \displaystyle
			\Big( 1 \wedge \frac{\delta^t_\wedge(x,y)}{\rho(x,y)}\Big)^{\alpha(1-\beta)} ,& q \ge p >1-\beta,\\[9pt]
			\displaystyle\Big(1 \wedge \frac{\delta^t_\wedge(x,y)}{\rho(x,y)} \Big)^{\alpha (1-\beta) }
			\log \Big(e+\frac{ \rho(x,y) \wedge  \delta^t(y)}{ \delta^t(x) }\Big), 	
			& q> 1-\beta = p, \\[9pt]
			\displaystyle\Big(1 \wedge \frac{\delta^t_\wedge(x,y)}{\rho(x,y)} \Big)^{\alpha p }\Big(1 \wedge \frac{\delta^t(y)}{\rho(x,y)} \Big)^{\alpha (1-\beta-p) }, & q>1-\beta>p, \\[9pt]
				\displaystyle\Big(1 \wedge \frac{\delta^t_\wedge(x,y)}{\rho(x,y)} \Big)^{\alpha (1-\beta) }
			\log \Big(e+\frac{ \rho(x,y) \wedge  \delta^t_\vee(x,y)}{ \delta^t_\wedge(x,y) }\Big), 	
			& q=1-\beta = p, \\[9pt]
			\displaystyle\Big(1 \wedge \frac{\delta^t_\wedge(x,y)}{\rho(x,y)} \Big)^{\alpha p } \Big(1 \wedge \frac{\delta^t(y)}{\rho(x,y)} \Big)^{\alpha(1-\beta-p)}   
			\log \Big(e+\frac{\delta^t(x) \wedge \rho(x,y)}{ \delta^t(y) }\Big), 
			& 1-\beta=q>p>0, \\[9pt]
			\displaystyle\Big(1 \wedge \frac{\delta^t_\wedge(x,y)}{\rho(x,y)} \Big)^{\alpha p } \Big(1 \wedge \frac{\delta^t_\vee(x,y)}{\rho(x,y)} \Big)^{\alpha(1-\beta-q)  }  \Big(1 \wedge \frac{\delta^t(y)}{\rho(x,y)} \Big)^{\alpha (q-p)}  , &1-\beta>q>1-\beta-p,\\[9pt]
			\displaystyle\Big(1 \wedge \frac{\delta^t_\wedge(x,y)}{\rho(x,y)} \Big)^{\alpha p } \Big(1 \wedge \frac{\delta^t_\vee(x,y)}{\rho(x,y)} \Big)^{\alpha p }\Big(1 \wedge \frac{\delta^t(y)}{\rho(x,y)} \Big)^{\alpha (q -p)}	\log \left(  e+\frac{  \rho(x,y)}{ \delta^t_\vee(x,y) } \right),   &1-\beta>q=1-\beta-p,\\[9pt]
			\displaystyle \Big(1 \wedge \frac{\delta^t(y)}{\rho(x,y)} \Big)^{\alpha q}	\log \left(  e  + \frac{\rho(x,y)}{ \delta^t(y)}\right),   &q=1-\beta=1-\beta-p,\\[9pt]
			\displaystyle \Big(1 \wedge \frac{\delta^t_\wedge(x,y)}{\rho(x,y)} \Big)^{\alpha p } \Big(1 \wedge \frac{\delta^t_\vee(x,y)}{\rho(x,y)} \Big)^{\alpha p} \Big(1 \wedge \frac{\delta^t(y)}{\rho(x,y)} \Big)^{\alpha (q -p)},   &q<1-\beta-p.
		\end{cases}
	\end{align}
		\textit{(Large time estimates)} 
If  ${\bf HK^{h_{p}}_B}$ holds, then  for all $(t,x,y) \in [1,\infty) \times D \times D$, 
	$$
	q(t,x,y) \simeq   	e^{-t \phi(\lambda_D) }\delta(x)^{\alpha p} \delta(y)^{\alpha q},
	$$
	and if ${\bf HK^{h_{p}}_U}$ holds, then \eqref{e:heatkernel-estimates} holds  for all $(t,x,y) \in (0,\infty) \times D \times D$. 
	
	\medskip

	Lastly,   we give estimates on the jump kernel $J(x,y)$.
	By  Theorem \ref{t:jump-estimate} and the fact that $\sB_h^*(x,y) \simeq \sB_h(0,x,y)$ for $x,y \in D$, we deduce from \eqref{e:heatkernel-estimates} that for any $q\ge p \ge 0$, 
		
	$$
	J(x,y) \simeq \frac{B_{p,q}(0,x,y;C_0)}{\rho(x,y)^{d+\alpha\beta}}, \quad x,y \in D.
	$$

}

\end{example}

\begin{example}
	{\rm
		Let $D:=\{x\in \R^d: x_d>0\}$ be the upper half space in  $\R^d$ and  
		$q\in [\alpha-1, \alpha)\cap (0, \alpha)$.
		We recall 
		the process $Y^D$ from Example \ref{ex:HK} (b-4),  
		 which corresponds to the Feynman-Kac semigroup of 
		the part process $Z^D$, in $D$, of the reflected  isotropic 
		$\alpha$-stable process in $\overline{D}$  
		via the multiplicative functional 
		$\exp(-\int^t_0C(d, \alpha, q)(Z^D_s)_d^{-\alpha}ds)$, 
		where the positive constants $C(d, \alpha, q)$ is defined on
		\cite[p.~233]{CKSV}.
		It is easy to see that $Y^D$ satisfies the scaling property and horizontal translation
		invariance, more precisely, for any $\lambda>0$, the transition density $p_D$ of
		$Y^D$ satisfies
		$$
		p_D(\lambda^\alpha, \lambda x, \lambda y)=p_D(t, x, y)\lambda^{-d}, \quad
		t>0, x, y\in D
		$$
		and
		$$
		p_D(t, x+(\widetilde{z}, 0), y+(\widetilde{z}, 0))=p_D(t, x, y), \quad
		t>0, x, y\in D, \widetilde{z}\in \R^{d-1}.
		$$
 Let $S=(S_t)_{t\ge 0}$ be a $\beta$-stable subordinator independent of the process 
  $Y^D$, 
 $\beta\in (0, 1)$. Then the 
		process $X_t:=Y^D_{S_t}$ falls into the framework of the present paper and thus we
		can get sharp two-sided estimates on the jump kernel, heat kernel and Green 
		functions of $X$. By using the scaling property and  horizontal translation
		invariance of $p_D$ above, we can show that the killing function $\kappa(x)$
		of $X$ is given by
		$$
		\kappa(x)=Cx_d^{-\alpha\beta}, \quad x\in D
		$$
		for some constant $C\in (0, \infty)$.  
		One can check,  by following arguments at the end of \cite[Section 2]{KSV200},  that the jump kernel of the subordinate process $X_t:=Y^D_{S_t}$ satisfies assumptions in \cite[(A1)-(A4)]{KSV21}.
				
		Moreover, by comparing the Green function estimates  in 
		\cite[Theorem 1.1]{KSV21}
		with the Green function function estimates in Example \ref{ex:basic}, 
		one can see that the value
		of the constant in the critical killing potential is indeed related to the power of the decay
		correctly.
		Thus, 
		instead of computing the constant $C$ of 
		the killing function $\kappa(x)$, 
		we see that the exponent $p$ in  \cite{KSV21} should be $q$ (and the constant $\alpha$ in \cite{KSV21} is equal to $\alpha \beta$ in the present case) and we can use   \cite[Theorems 1.2 and 1.3]{KSV21} directly.
		Therefore, by checking the range of $p=q/\alpha$ in the jump
		kernel estimates in Example \ref{ex:basic} and \cite[Theorems 1.2 and 1.3]{KSV21}, 
		from  \cite[Theorems 1.2 and 1.3]{KSV21} we obtain the following corollary. See \cite[Theorem 1.2]{KSV21} for the precise statement of the scale-invariant boundary Harnack principle.

		\begin{cor}\label{c:BHP}
			Suppose  $d > \alpha $ and
			 $X$ 
			 is the  process defined above. Let $\beta_{1/2}:=\beta \vee 1/2$.
			If $q\in [\alpha-1, \alpha)\cap ((\alpha\beta-1)_+, \alpha\beta_{1/2})$ , then
			the scale-invariant boundary Harnack principle is valid for $X$.  
			If $ q\in [\alpha\beta_{1/2} \vee (\alpha-1), \alpha)$, then
			the non-scale-invariant boundary Harnack principle is not valid for $X$.
		\end{cor} 
			Note that constants $(\beta_1, \beta_2)$
		in \cite{KSV21} are 
		$(\alpha-\alpha \beta, 0)$  for $\alpha(1-\beta) \le q<\alpha$,
		$(q, \alpha-q-\alpha \beta)$  for $\alpha(1-\beta)/2 < q<\alpha(1-\beta) $ and
	$(q, q)$  for $0 \le q \le \alpha(1-\beta)/2$.

	}
\end{example}

The next example illustrates that any of the four terms in \eqref{e:general-off} can not dominates  all the other terms in general.

\begin{example}
	{\rm
Let  $0<\beta<1<L$  
and $S$ be a subordinator 
with tail L\'evy measure $w$ satisfying
\begin{equation*}	w(r) \simeq r^{-\beta} \wedge r^{-L}, \quad r>0.\end{equation*}
 Let $0<\alpha\le K$.  Suppose that the heat kernel $p_D(t,x,y)$ of $Y^D$ enjoys the estimate ${\bf HK^{h_{1/2}}_U}$ with
\begin{equation*}
	\Phi(r)=r^{\alpha}, \quad\;\; \Psi(r)=r^{\alpha} \vee r^{K}, \quad r>0.
\end{equation*} 
Examples of such $Y^D$ can be found in \cite[Theorem 1.4]{KM18} where Dirichlet heat kernel estimates for a large class of subordinate Brownian motions are treated. 
Recall that $\psi(r):=\phi(1/\Phi(r))^{-1}$. Note that {\bf (Poly-$\infty$)} holds, and 
by \eqref{e:phi-w}, 
\begin{align}\label{e:polydecay}
	\phi(1/r)\simeq &\begin{cases}
		r^{-\beta}, &\mbox{if } \; r \le 1\\
		r^{-1}, &\mbox{if } \; r > 1,
	\end{cases}
	&\psi(r) \simeq &\begin{cases}
		r^{\alpha\beta}, &\mbox{if } \; r \le 1\\
		r^{\alpha}, &\mbox{if } \; r > 1,
	\end{cases}\nn\\
\phi^{-1}(1/r)^{-1}\simeq&\begin{cases}
	r^{1/\beta}, &\mbox{if } \; r \le 1\\
	r, &\mbox{if } \; r > 1,
\end{cases}  &\psi^{-1}(r)\simeq&\begin{cases}
	r^{1/(\alpha\beta)}, &\mbox{if } \; r \le 1\\
	r^{1/\alpha}, &\mbox{if } \; r > 1.
\end{cases}
\end{align}
Let $q(t,x,y)$ be the heat kernel of the subordinate process $X_t:=Y^D_{S_t}$. We get global estimates on $q(t,x,y)$ from Theorem
\ref{t:largegeneral} and Lemma \ref{l:asym-Bp}.
Write the terms on the right-hand side of  \eqref{e:general-off} respectively as
\begin{align*}
 A_1&=\frac{C_0t\sB_h(t,x,y)}{V(x,\rho(x,y)) \Psi(\rho(x,y))}, 
& A_2&=\frac{C_0\phi^{-1}(1/t)^{-1}h(\phi^{-1}(1/t)^{-1},x,y)}{V(x,\rho(x,y)) \Psi(\rho(x,y))}, \\
 A_3&=\frac{h(\phi^{-1}(1/t)^{-1},x,y)}{V(x,\psi^{-1}(t))} \exp\Big(-c\frac{\rho(x,y)^2}{\psi^{-1}(t)^2}\Big), &A_4&=\frac{t h(\Phi(\rho(x,y)), x, y)  w(\Phi(\rho(x,y)))}{V(x,\rho(x,y))}.
\end{align*}

\medskip

\textit{(Small time estimates)} Note that for any $R_1>0$,  
$w$ satisfies {\bf (Poly-$R_1$)} with both the upper and the lower index equal to 
 $\beta 
\in (0,1)$.  Using this fact, it is easy to see that for all $(t,x,y) \in (0,1/\phi(4)] \times D \times D$ with $\rho(x,y)\ge t^{1/(\alpha\beta)}$, $A_4$ dominates $A_3$  and   $A_1$ dominates $A_2$. Thus, the following  estimates hold for all $(t,x,y) \in (0,1/\phi(4)] \times D \times D$.

\smallskip

(i) If $\rho(x,y)\le \phi(4)^{-1/(\alpha\beta)}$, then 
$$
q(t,x,y) \simeq  \Big(1 \wedge \frac{\delta(x)}{t^{1/(\alpha\beta)}} \Big)^{\alpha/2} \Big(1 \wedge \frac{\delta(y)}{t^{1/(\alpha\beta)}} \Big)^{\alpha/2} B_{1/2,1/2}(t,x,y ;C_0) \left( \frac{1}{V(x, t^{1/(\alpha\beta)})} \wedge \frac{t  }{V(x,\rho(x,y))\rho(x,y)^{\alpha\beta}}\right),
$$
where the function $B_{1/2,1/2}(t,x,y;C_0)$ is defined as  in  \eqref{e:explicit-boundary-1} and \eqref{e:explicit-boundary-2}.

(ii) Let $\rho(x,y)>\phi(4)^{-1/(\alpha\beta)}$. Since $L>1$, we see that for all $u>1$, 
\begin{equation}\label{e:ex_mixedBh} \int_u^{\infty}h_{1/2}(s,x,y)s^{-L}ds \le h_{1/2}(u,x,y) \int_u^{\infty}s^{-L}ds \le \frac{u^{1-L}}{L-1}h_{1/2}(u,x,y) .
\end{equation}
Thus, we get from \eqref{e:intscale_lower},  \eqref{e:equiv-form} and \eqref{e:polydecay} that
\begin{align*}
\sB_{h_{1/2,1/2}}(t,x,y)&\simeq \int_{2\phi^{-1}(1/t)^{-1}}^{1} h_{1/2}(s,x,y)s^{-\beta}ds+\int_{1}^{ 4\rho(x,y)^{\alpha}} h_{1/2}(s,x,y)s^{-L}ds\\
&\simeq \int_{2\phi^{-1}(1/t)^{-1}}^{1} h_{1/2}(s,x,y)s^{-\beta}ds + h_{1/2}(1,x,y) \\
&\simeq \int_{2\phi^{-1}(1/t)^{-1}}^{1} h_{1/2}(s,x,y)s^{-\beta}ds\\
&= h_{1/2}(2\phi^{-1}(1/t)^{-1},x,y) \int_{2\phi^{-1}(1/t)^{-1}}^{1}   \Big(1 \wedge \frac{\delta^t_\wedge(x,y)^{\alpha}}{s} \Big)^{1/2} \Big(1 \wedge \frac{\delta^t_\vee(x,y)^{\alpha}}{s} \Big)^{1/2} s^{-\beta}ds.
\end{align*}
In the end, by similar arguments to that given in the proof of Lemma \ref{l:asym-Bp}, using \eqref{e:polydecay}, we conclude that
\begin{align*}
	&q(t,x,y) \simeq  \Big(1 \wedge \frac{\delta(x)}{t^{1/(\alpha\beta)}} \Big)^{\frac\alpha2} \Big(1 \wedge \frac{\delta(y)}{t^{1/(\alpha\beta)}} \Big)^{\frac\alpha2} \\
	&\quad \times \bigg[\Big( 1 \wedge \frac{\delta_\wedge(x,y) \vee t^{1/(\alpha\beta)}}{\rho(x,y)} \Big)^{\frac\alpha2} \Big( 1 \wedge \frac{\delta_\vee(x,y) \vee t^{1/(\alpha\beta)}}{\rho(x,y)} \Big)^{\frac\alpha2} \frac{t}{V(x,\rho(x,y))\rho(x,y)^{\alpha L}} +  \frac{C_0t F_\beta(t,x,y)}{V(x,\rho(x,y))\rho(x,y)^{K}}\bigg],
\end{align*}
where
\begin{align*} 
F_\beta(t,x,y):= \begin{cases}
		\big(1 \wedge (\delta_\wedge(x,y) \vee t^{1/(\alpha\beta)})\big)^{\alpha/2}  \,	\big(1 \wedge (\delta_\vee(x,y) \vee t^{1/(\alpha\beta)})\big)^{\alpha/2-\alpha\beta}, &\mbox{if }\, \beta<1/2,\\[2pt]
	\displaystyle\big(1 \wedge (\delta_\wedge(x,y) \vee t^{1/(\alpha\beta)})\big)^{\alpha/2}  \log \Big(e + \frac{1 \wedge (\delta_\vee(x,y) \vee t^{1/(\alpha\beta)})}{1 \wedge (\delta_\wedge(x,y) \vee t^{1/(\alpha\beta)})}\Big), &\mbox{if } \beta=1/2,\\[8pt]
	\big(1 \wedge (\delta_\wedge(x,y) \vee t^{1/(\alpha\beta)})\big)^{\alpha-\alpha\beta}, &\mbox{if } \beta>1/2.
	\end{cases}
\end{align*}
Therefore, when $C_0=1$, since $K$ or $L$ may be arbitrarily large number, one should keep both terms $A_1$ and $A_4$ for off-diagonal estimates in general. 

\begin{comment}
\noindent (i) If $\rho(x,y)<t^{1/(\alpha_1\beta_1)}$, then
\begin{align*}
	q(t,x,y) \simeq &\Big(1 \wedge \frac{\delta(x)}{t^{1/(\alpha_1\beta_1)}} \Big)^{\alpha_1/2} \Big(1 \wedge \frac{\delta(y)}{t^{1/(\alpha_1\beta_1)}} \Big)^{\alpha_1/2} \frac{1}{V(x, t^{1/(\alpha_1\beta_1)} )}.
\end{align*}

\noindent (ii) Let $C_0=0$ and $\rho(x,y) \ge t^{1/(\alpha_1\beta_1)}$.
It is easy to see that $A_4$ dominates $A_3$. Hence,  we have that
\begin{align*}
	q(t,x,y) \simeq &  \Big(1 \wedge \frac{\delta(x)}{\rho(x,y)} \Big)^{\alpha_1/2} \Big(1 \wedge \frac{\delta(y)}{\rho(x,y)} \Big)^{\alpha_1/2} \frac{t}{V(x, \rho(x,y)) [\rho(x,y)^{\alpha_1\beta_1} \vee \rho(x,y)^{\alpha_1\beta_2}]}.
\end{align*}

\noindent (iii) Let $C_0=1$ and $t^{1/(\alpha_1\beta_1)} \le \rho(x,y) <1$.
First observe that $t w(\phi^{-1}(1/t)^{-1})\simeq 1$. Thus it follows from  \eqref{e:intscale_lower}, {\bf (Poly-$R_1$)} and \eqref{e:boundary-function-c} 
that $A_1$ dominates $A_2$ and $A_4$ (hence also $A_3$). 
Further, note that for any $R_1>0$,  
$w$ satisfies {\bf (Poly-$R_1$)} with both the upper and the lower index equal to $\beta_1\in (0,1)$. Therefore we 
see that  \eqref{e:subordinate_alpha_beta} holds with $\alpha=\alpha_1$, $\beta=\beta_1$ and $p=1/2$, after multiplying $\rho(x,y)^dV(x,\rho(x,y))^{-1}$ in each case.

\end{comment}

\smallskip

 \textit{(Large time estimates)} By \eqref{e:ex_mixedBh},  we see that  for all $(t,x,y) \in [1/\phi(4),\infty) \times D \times D$ with $\rho(x,y)\ge t^{1/\alpha}$, $A_2$ dominates $A_1$ since $t\phi^{-1}(1/t)^{-L} \simeq t^{1-L} \le \phi(4)^{L-1}$ in this case. Thus, the following estimates hold for all $(t,x,y) \in [1/\phi(4),\infty) \times D \times D$. 
 
 \smallskip  (i) If $\rho(x,y) < t^{1/\alpha}$, then
 \begin{align}\label{e:ex_mixed_large_on}
 	q(t,x,y) \simeq \Big(1 \wedge \frac{\delta(x)}{t^{1/\alpha}} \Big)^{\alpha/2} \Big(1 \wedge \frac{\delta(y)}{t^{1/\alpha}} \Big)^{\alpha/2} \frac{1}{V(x, t^{1/\alpha} )}.
 \end{align}
 
 (ii) If $\rho(x,y) \ge t^{1/\alpha}$, then 
 \begin{align}\label{e:ex_mixed_large_off}
 	q(t,x,y) &\asymp    
	\Big(1 \wedge \frac{\delta(x)}{t^{1/\alpha}} \Big)^{\alpha/2} \Big(1 \wedge \frac{\delta(y)}{t^{1/\alpha}} \Big)^{\alpha/2} \bigg[\frac{C_0t}{V(x,\rho(x,y)) \rho(x,y)^{K}} +\frac{1}{V(x, t^{1/\alpha} )} \exp\Big(-c \frac{\rho(x,y)^2}{t^{2/\alpha}} \Big) \nn\\
 	& \qquad  +    \Big(1 \wedge \frac{\delta(x)\vee t^{1/\alpha}}{\rho(x,y)} \Big)^{\alpha/2} \Big(1 \wedge \frac{\delta(y) \vee t^{1/\alpha}}{\rho(x,y)} \Big)^{\alpha/2} \frac{t}{V(x,\rho(x,y))\rho(x,y)^{\alpha L}}\bigg].
 \end{align}
  In particular, one should keep all terms $A_2$, $A_3$ and $A_4$ for off-diagonal estimates in general. 
 
 Note that, in case of $C_0=1$ and $L \ge K/\alpha$, we see from  \eqref{e:ex_mixed_large_on} and \eqref{e:ex_mixed_large_off} that
 \begin{equation*}
 	q(t,x,y) \asymp p_D(ct,x,y), \quad t \ge 1/\phi(4), \; x,y \in D.
 \end{equation*}

}
\end{example}

\vspace{.1in}
\textbf{Acknowledgment}:
 We thank the referee for the comments on the first version of this paper. The third-named 
author thanks Professor Habib Maagli for pointing a mistake in
Theorem 3.10 of \cite{So}, which led to a mistake in Theorem 4.6 there.
Part of the research for this paper was done while the third-named author was visiting
Jiangsu Normal University, where he was partially supported by  a grant from
the National Natural Science Foundation of China (11931004) and by
the Priority Academic Program Development of Jiangsu Higher Education Institutions.

%%%%%%%%%%%%%%%%%%%%%%%%%%%%%%%%%%%%%%%%%%%%%%%%%%%%%%%%%%%%%%%%%%%%%%%%%%%%%%%%%%%
%%%%%%%%%%%%%%%%%%%%%%%%%%%%%%%%%%%%%%%%%%%%%%%%%%%%%%%%%%%%%%%%%%%%%%%%%%%%%%%%%


\vspace{.1in}

\smallskip
\begin{thebibliography}{99}
	



\bibitem{AGV}
N. Abatangelo, D. Gomez-Castro and J. L V\'azquez. 
Singular boundary behaviour and large solutions for fractional elliptic equations. arXiv: 1910.00366


\bibitem{BGR}
K. Bogdan, T. Grzywny, and M. Ryznar. {Heat kernel
  estimates for the fractional {L}aplacian with {D}irichlet conditions}. Ann.
  Probab. \textbf{38} (2010), 1901--1923. 
  
  
  \bibitem{BGR3}
K. Bogdan, T. Grzywny, and M. Ryznar. {Dirichlet heat kernel for unimodal {L}\'evy processes}.
  Stochastic Process. Appl. \textbf{124} (2014),  3612--3650.


\bibitem{BGR2}
K. Bogdan, T. Grzywny, and M. Ryznar. {Barriers, exit time and survival probability for unimodal
  {L}\'evy processes}. Probab. Theory Related Fields \textbf{162} (2015),
   155--198. 


\bibitem{BSV} M. Bonforte, Y. Sire and J. L. V\'azquez.
Existence, uniqueness and asymptotic behaviour 
for fractional porous medium equations on bounded domains. Discrete Contin. Dyn. Syst. \textbf{35} (2015), 5725--5767. 



\bibitem{Bo84} N. Bouleau. Quelques résultats probabilistes sur la subordination au sens de Bochner.
Seminar on potential theory, Paris, No. 7, 54–81,
Lecture Notes in Math., 1061, Springer, Berlin, 1984.


\bibitem{CK} Z.-Q. Chen and  P. Kim.
Global Dirichlet heat kernel estimates for symmetric L\'{e}vy processes in half-space. Acta Appl. Math. \textbf{146} (2016), 113--143

	
\bibitem{CKKW} Z.-Q. Chen, P. Kim, T. Kumagai and J. Wang.
Heat kernel estimates for time fractional equations. Forum Math.  \textbf{30}  (2018),  1163--1192.


\bibitem{CKS-jems}
Z.-Q. Chen, P. Kim and R. Song.
Heat kernel estimates for the  Dirichlet fractional Laplacian. 
J. Eur. Math. Soc. \textbf{12} (2010), 1307--1329.

\bibitem{CKS-ptrf} 
Z.-Q. Chen, P. Kim and R. Song.
Two-sided heat kernel estimates for censored stable-like processes. 
Probab. Theory Relat. Fields, \textbf{146} (2010), 361--399.


\bibitem{CKS50}
Z.-Q. Chen, P. Kim and R. Song. {Heat kernel estimates for {$\Delta+\Delta^{\alpha/2}$} in
  {$C^{1,1}$} open sets}, J. Lond. Math. Soc. (2) \textbf{84} (2011), no.~1,
  58--80.

\bibitem{CKS5}
Z.-Q. Chen, P. Kim and R. Song. {Dirichlet heat kernel estimates for fractional {L}aplacian with
  gradient perturbation}, Ann. Probab. \textbf{40} (2012), 2483--2538.


\bibitem{CKS-ejp}
Z.-Q. Chen, P. Kim and R. Song.
Global heat kernel estimates for $\Delta+\Delta^{\alpha/2}$ in half-space-like domains. Elect. J. Probab. \textbf{17} (2012), No. 32, 1--32.


\bibitem{CKS} Z.-Q. Chen, P. Kim and R. Song.
Dirichlet heat kernel estimates for rotationally symmetric Lévy processes. Proc. Lond. Math. Soc. \textbf{109} (2014), 90--120.



\bibitem{CKS16} Z.-Q. Chen, P. Kim and R. Song.
Dirichlet Heat Kernel Estimates for Subordinate Brownian Motions with Gaussian Components. J. Reine Angew. Math. \textbf{711} (2016), 111--138.



\bibitem{CKK09} Z.-Q. Chen, P. Kim and T. Kumagai. On heat kernel estimates and parabolic Harnack
inequality for jump processes on metric measure spaces. Acta Math. Sin. (Engl. Ser.) 25
(2009), 1067–1086.


\bibitem{CK03} Z.-Q. Chen and T. Kumagai.
Heat kernel estimates for stable-like processes on d-sets.
Stochastic Process. Appl. \textbf{108}  (2003), no. 1, 27–62.

\bibitem{CKW16a} Z.-Q. Chen, T. Kumagai, and J. Wang.
Stability of heat kernel estimates for symmetric non-local {D}irichlet forms.
Mem. Amer. Math. Soc. \textbf{271} (2021), no. 1330.


\bibitem{CKW20} Z.-Q. Chen, T. Kumagai, and J. Wang. Stability of parabolic Harnack inequalities for symmetric non-local Dirichlet forms.  J. Eur. Math. Soc. \textbf{22} (2020), 3747--3803.


\bibitem{CS} Z.-Q. Chen and R. Song. 
Intrinsic ultracontractivity and conditional gauge for symmetric stable processes. 
J. Funct. Anal., \textbf{150} (1997), 204--239.

\bibitem{CT} Z.-Q. Chen and J. Tokle. 
Global heat kernel estimates for fractional Laplacians in unbounded open sets. 
Probab. Theory Relat. Fields, \textbf{149} (2011), 373--395.



	
\bibitem{CK1} S. Cho and P. Kim.
Estimates on the tail probabilities of subordinators and applications to general time fractional equations. 	
Stochastic Process. Appl. \textbf{130} (2020), 4392--4443.



\bibitem{CK20} S. Cho and P. Kim.
Estimates on transition densities of subordinators with jumping density decaying in mixed polynomial orders. 
Stochastic Process. Appl. \textbf{139} (2021), 229–279.


\bibitem{CKP}   S. Cho, P. Kim. and H. Park. 
Two-sided estimates on Dirichlet heat kernels for time-dependent parabolic operators with singular drifts in $C^{1,\alpha}$-domains. 
J. Differential Equations, \textbf{252} (2012), 1101--1145.



\bibitem{CKSV} S. Cho, P. Kim, R. Song, R. and Z. Vondra\v{c}ek.
Factorization and estimates of Dirichlet heat kernels for non-local operators with critical killings. 
J. Math. Pures Appl. \textbf{143} (2020), 208--256.


\bibitem{CZ} K. L. Chung and Z. X. Zhao.
From Brownian motion to Schrödinger’s equation. 
Springer-Verlag, Berlin, 1995.

\bibitem{Da} E. B. Davies.
Heat kernels and spectral theory. Cambridge University Press, Cambridge, 1989.


\bibitem{DD} J. Davila and L. Dupaigne.
Hardy-type inequalities. 
J. Eur. Math. Soc.  \textbf{6} (2004), 335–-365.

\bibitem{FMT}
S. Filippas, L. Moschini and A.Tertikas.
Sharp two-sided heat kernel estimates for critical Schrödinger operators on bounded domains. Comm. Math. Phys. \textbf{273} (2007),  237–-281.


\bibitem{GS} P. Gyrya and L. Saloff-Coste.
Neumann and Dirichlet heat kernels in inner uniform domains.
Ast\'erisque, \textbf{336} (2011), viii+144 pp.


  \bibitem{GKK} 
T. Grzywny, K.-Y. Kim and P. Kim.
{Estimates of Dirichlet heat kernel for symmetric Markov processes.}  
Stochastic Process. Appl. \textbf{130} (2020), no. 1, 431--470.




\bibitem{JP} N. C. Jain and W. E. Pruitt.
Lower tail probability estimates for subordinators and nondecreasing random walks. Ann. Probab.  \textbf{15} (1987),  75--101. 

\bibitem{JW}
T. Jakubowski and J. Wang.
Heat kernel estimates of fractional Schr\"odinger operators with negative Hardy potential. 
Potential Anal. \textbf{53} (2020),  997–1024.

\bibitem{KS} K. Kaleta and P. Sztonyk.
Estimates of transition densities and their derivatives for jump L\'evy processes. J. Math. Anal. Appl. \textbf{431} (2015), 260--282.

\bibitem{K}
K.-Y. Kim. {Global heat kernel estimates for symmetric {M}arkov
  processes dominated by stable-like processes in exterior {$C^{1,\eta}$} open
  sets}. Potential Anal. \textbf{43} (2015),  127--148.

\bibitem{KK}
K.-Y. Kim and P. Kim. {Two-sided estimates for the transition
  densities of symmetric {M}arkov processes dominated by stable-like processes
  in {$C^{1,\eta}$} open sets}. Stochastic Process. Appl. \textbf{124} (2014),
  3055--3083.  Erratum to "Two-sided estimates for the transition densities of symmetric Markov processes dominated by stable-like processes in {$C^{1,\eta}$} open sets'' [Stochastic Process. Appl. 124 (2014) 3055--3083]. Stochastic Process. Appl. \textbf{129} (2019),  2228--2230.
  
  

\bibitem{KM14} P. Kim and A. Mimica.
Green function estimates for subordinate Brownian motions: Stable and beyond.
Trans. Amer. Math. Soc. \textbf{366} (2014), No. 8, 4383--4422.



\bibitem{KM18} P. Kim and A. Mimica.
Estimates of Dirichlet heat kernels for subordinate Brownian motions. Electron. J. Probab. \textbf{23} (2018), No. 64, 45 pp. 


\bibitem{KSo}
P. Kim and R. Song. {Dirichlet heat kernel estimates for stable
  processes with singular drift in unbounded {$C^{1,1}$} open sets}. Potential
  Anal. \textbf{41} (2014), 555--581.

\bibitem{KSV-spa14} 
P. Kim, R. Song and Z. Vondra\v{c}ek. 
Global uniform boundary Hanack principle with explicit decay rate and its application. Stoch. Proc. Appl. \textbf{124} (2014), 235--267.


\bibitem{KSV19}
P. Kim, R. Song and  Z. Vondra\v{c}ek.
Potential theory of subordinate killed Brownian motions. 
Trans. Amer. Math. Soc. \textbf{371} (2019), 3917--3969.


\bibitem{KSV20}
P. Kim, R. Song and  Z. Vondra\v{c}ek.
On the boundary theory of subordinate killed L\'evy processes. 
 Potential  Anal.  \textbf{53} (2020), 131--181.

\bibitem{KSV200}
P. Kim, R. Song and  Z. Vondra\v{c}ek.
On potential theory of Markov processes with jump kernels decaying at the boundary. 
To appear in  Potential Anal. 
DOI: 10.1007/s11118-021-09947-8 (2021).

\bibitem{KSV21}
P. Kim, R. Song and  Z. Vondra\v{c}ek.
Sharp two-sided Green function estimates for Dirichlet forms degenerate at the boundary. 
	arXiv:2011.00234v3  (2021).
	
	
	\bibitem{KSe20}
	D. Kinzebulatov and Yu. A. Semënov. 
Fractional Kolmogorov operator and desingularizing weights. 
	 arXiv:2005.11199v4 (2020). 
	
	
	\bibitem{MNS18}
G. Metafune, L. Negro and C. Spina.
Sharp kernel estimates for elliptic operators with
second-order discontinuous coefficients. J. Evol. Equ. \textbf{18} (2018), 467--514.



\bibitem{Mi} A. Mimica.
Heat kernel estimates for subordinate Brownian motions. Proc. Lond. Math. Soc. (3) \textbf{113} (2016), 627--648.


\bibitem{Miy} A. Miyake. The subordination of Lévy system for Markov processes.
Proc. Japan Acad. \textbf{45} (1969), 601--604.

\bibitem{So} R. Song.
 Sharp bounds on the density, Green function and jumping function of subordinate killed BM.
 Probab. Theory Related Fields, Vol. \textbf{128}  (2004), No. 4, 606--628.



\bibitem{SV08} R. Song and Z.  Vondra\v{c}ek. On the relationship between subordinate killed and killed subordinate processes. Electron. Commun. Probab. \textbf{13} (2008), 325--336.

\bibitem{Var}
N. Th. Varopoulos. Gaussian estimates in Lipschitz domains. 
Can. J. Math.  \textbf{55} (2003) 401--431.



\bibitem{Zhao} Z. Zhao.
Green function for Schrödinger operator and conditional Feynman-Kac
gauge. J. Math. Anal. Appl. \textbf{116} (1986), 309–334.


\bibitem{Zh} Q. S. Zhang. 
The boundary behavior of heat kernels of Dirichlet Laplacians. J. Diff.
Eqs. \textbf{182} (2002), 416--430.



\end{thebibliography}
\end{document}